\documentclass[12pt,twoside, reqno]{amsart}

\usepackage{amssymb}
\usepackage{amsmath}
\usepackage{mathrsfs}
\usepackage{amsthm}
\usepackage{amsfonts}
\usepackage{txfonts}
\usepackage{enumerate}
\usepackage{hyperref}
\hypersetup{hidelinks}
\usepackage{graphicx}
\usepackage{subfigure}
\usepackage{float}
\usepackage{latexsym,amssymb}
\usepackage{color}
\usepackage{graphicx}
\usepackage{pgf,tikz}
\usepackage{mathrsfs}
\usepackage{fancyhdr}
\usetikzlibrary{arrows}

\allowdisplaybreaks
\usepackage{indentfirst}
\usepackage{setspace}
\usepackage{geometry}
\newtheorem{theorem}{Theorem}[section]
\newtheorem{proposition}[theorem]{Proposition}
\newtheorem{corollary}[theorem]{Corollary}
\newtheorem{lemma}[theorem]{Lemma}
\newtheorem{remark}[theorem]{Remark}
\newtheorem{definition}[theorem]{Definition}

\numberwithin{equation}{section} 

\numberwithin{equation}{section}

\newcommand{\R}{\mathbb{R}}

\newcommand{\ben}{\begin{eqnarray*}}
\newcommand{\enn}{\end{eqnarray*}}
\newcommand{\pa}{\partial}

\newcommand{\g}{\gamma}

\newcommand{\ve}{\varepsilon}

\newcommand{\al}{\alpha}
\newcommand{\la}{\lambda}

\newcommand{\ol}{\overline}

\newcommand{\half}{\frac{1}{2}}

\newcommand{\na}{\nabla}
\newcommand{\be}{\begin{equation}}
\newcommand{\ee}{\end{equation}}
\newcommand{\ba}{\begin{aligned}}
\newcommand{\ea}{\end{aligned}}
\newcommand{\lf}{\left}
\newcommand{\rt}{\right}

\newcommand{\wt}{\tilde{w}}

\def\9{{\infty}}
\def\a{{\alpha}}
\def\b{{\beta}}
\def\g{{\gamma}}

\def\lbb{{\lambda}}

\def\calo{{\mathcal{O}}}
\def\calp{{\mathcal{P}}}

\def\bbp{{\mathbb{P}}}
\def\bbr{{\mathbb{R}}}
\def\ve{{\varepsilon}}
\def\vf{{\varphi}}

\def\wt{\widetilde}
\def\wh{\widehat}
\def\ol{\overline}
\def\({\left(}
\def\){\right)}
\def\<{\langle}
\def\>{\rangle}
\def\scri{{\mathscr{I}}}
\def\cale{{\mathcal{E}}}

\def\cald{{\mathcal{D}}}
\def\calh{{\mathcal{H}}}
\def\calk{{\mathcal{K}}}

\begin{document}

\title
[Multi-bubble Bourgain-Wang solutions]{Multi-bubble Bourgain-Wang solutions to nonlinear Schr\"odinger equation}

\author{Michael R\"ockner}
\address{Fakult\"at f\"ur Mathematik,
Universit\"at Bielefeld, D-33501 Bielefeld, Germany}
\email{roeckner@math.uni-bielefeld.de}
\thanks{}

\author{Yiming Su}
\address{Department of mathematics,
Zhejiang University of Technology, 310014 Zhejiang, China}
\email{yimingsu@zjut.edu.cn}
\thanks{}

\author{Deng Zhang}
\address{School of mathematical sciences,
Shanghai Jiao Tong University, 200240 Shanghai, China }
\email{dzhang@sjtu.edu.cn}
\thanks{}

\subjclass[2010]{35B44, 35B40, 35Q55}

\keywords{Bourgain-Wang solutions, $L^2$-critical, multi-bubbles, nonlinear Schr\"odinger equation, \quad
non-pure multi-solitons}

\date{}


\begin{abstract}
  We consider a general class of focusing $L^2$-critical nonlinear Schr\"odinger equations
  with lower order perturbations,
  for which the pseudo-conformal symmetry and the conservation law of energy are absent.
  In dimensions one and two,
  we construct Bourgain-Wang type solutions
  concentrating at $K$ distinct singularities,
  $1\leq K<\9$,
  and prove that they are unique if the asymptotic behavior is within the order $(T-t)^{4+}$,
  for $t$ close to the blow-up time $T$.
  These results apply to the canonical nonlinear Schr\"odinger equations
  and, through the pseudo-conformal transform,
  in particular yield the existence and conditional uniqueness of
  non-pure multi-solitons.
  Furthermore, through a Doss-Sussman type transform,
  these results also apply to stochastic nonlinear Schr\"odinger equations,
  where the driving noise is taken in the sense of controlled rough path.
\end{abstract}

\maketitle
\begin{spacing}{0.95}
\tableofcontents
\end{spacing}

\section{Introduction and main results}    \label{Sec-Intro}

\subsection{Introduction}

We consider a general class of focusing $L^2$-critical nonlinear Schr\"{o}dinger equations
with lower order perturbations
\begin{align} \label{equa-NLS-perturb}
    i\partial_t v +\Delta v + a_1 \cdot \nabla v + a_0 v +|v|^{\frac{4}{d}}v  =0
\end{align}
on $\bbr^d$, where $d=1,2$, the coefficients of lower order perturbations are of form
\begin{align}
 & a_1(t,x)= 2 i \sum\limits_{l=1}^N \na \phi_l(x)h_l(t),   \label{a1-loworder} \\
 & a_0(t,x)= - \sum\limits_{j=1}^d \(\sum\limits_{l=1}^N \partial_j \phi_l(x)h_l(t)\)^2
          + i \sum\limits_{l=1}^N \Delta \phi_l(x) h_l(t), \label{a0-loworder}
\end{align}
and $\phi_l \in C_b^\9(\bbr^d, \bbr)$,
$h_l \in C(\bbr^+; \bbr)$, $1\leq l\leq N$.

Equation \eqref{equa-NLS-perturb} is mainly motivated by two canonical models.
The first is the nonlinear Schr\"odinger equation (NLS),
corresponding to the case without lower order perturbations, i.e.,
\begin{align} \label{equa-NLS}
    i\partial_t v +\Delta v +|v|^{\frac{4}{d}}v =0.
\end{align}
NLS is a canonical equation of major importance
in continuum mechanics, plasma physics and optics (\cite{DNPZ92}).
In particular, for the cubic nonlinearity in the critical dimension two,
the phenomenon of mass concentration near collapse gives a rigorous basis
to the physical concept of ``strong collapse'' (\cite{SS99}).
For more physical interpretations we refer to \cite{DNPZ92,GS11,SS99}.

Another important model
is the stochastic nonlinear Schr\"odinger equation (SNLS)
\begin{align} \label{equa-SNLS}
    & idX + \Delta X dt + |X|^{\frac 4d} X dt = - i \mu X dt +i X dW(t),
\end{align}
where
$W$ is a Wiener process of form
$$W(t,x)=\sum_{l=1}^N i\phi_l(x)B_l(t),\ \ x\in \bbr^d,\ t\geq 0,$$
$\{\phi_l\} \subseteq C_b^\9(\bbr^d, \bbr)$,
$\{B_l\}$ are standard $N$-dimensional real valued Brownian motions
on a normal stochastic basis $(\Omega, \mathscr{F}, \{\mathscr{F}_t\}, \bbp)$,
and
$\mu= \frac 12 \sum_{l=1}^N  \phi_l^2.$
The last term $XdW(t)$ in \eqref{equa-SNLS} is taken in the sense of controlled rough path
(see Definition \ref{def-X-rough} below).
The key relationship is that,
through the Doss-Sussman type transformation
$v:= e^{-W} X$,
$v$ satisfies  equation \eqref{equa-NLS-perturb}
with the   functions
$\{h_l\}$ being exactly the Brownian motions $\{B_l\}$.

The physical significance of SNLS is well known.
One significant model
arises from molecular aggregates
with thermal fluctuations,
where the multiplicative noise corresponds to
scattering of exciton by phonons,
due to thermal vibrations of the molecules.
In particular, for the cubic nonlinearity in dimension two,
the noise effect on the coherence of the ground state solitary solution
was studied in \cite{BCIR94,BCIRG95}.
The case of quintic nonlinearity
in the critical one dimensional case was studied in \cite{RGBC95}.
We also refer to \cite{BG09} for applications to open quantum systems.

It is known that,
equation \eqref{equa-NLS-perturb} is locally well-posed in the space $H^1$,
see, e.g., \cite{C03} for the NLS,
and \cite{BD03,BM14,BRZ16} for the SNLS.

The long time behavior of solutions
is, however, more delicate.
An important role here is played by the ground state,
which is a positive radial solution to the elliptic equation
\begin{align} \label{equa-Q}
    \Delta Q - Q + Q^{1+\frac{4}{d}} =0.
\end{align}
It is known that (see \cite[Theorem 8.1.1]{C03})
$Q$ is smooth and decays at infinity exponentially fast,
i.e., there exist $C, \delta>0$ such that for any multi-index $|\upsilon|\leq 3$,
\be\label{Q-decay}
|\partial_x^\upsilon Q(x)|\leq C e^{-\delta |x|}, \ \ x\in \bbr^d.
\ee

More importantly,
the mass of the ground state is the threshold of global well-posedness and blow-up.
As a matter of fact,
in the NLS case,
solutions with subcritical mass (i.e., $\|v\|_{L^2} < \|Q\|_{L^2}$)
exist globally and even scatter at infinity, \cite{D15,W83}.
In contrast to that,
in the critical mass regime,
two important dynamics are exhibited:
the non dispersive solitary wave
\begin{align}  \label{W-soliton-intro}
  W(t,x):=w^{-\frac d2}Q \(\frac{x- c t}{w} \)e^{i(\half c \cdot x-\frac{1}{4}|c|^2t+w^{-2}t+\vartheta)},
\end{align}
and the pseudo-conformal blow-up solutions
\begin{align}  \label{S-blowup-intro}
    S_T(t,x)=(w(T-t))^{-\frac d2}Q \(\frac{x-x^*}{w(T-t)}\)
             e^{ - \frac i 4 \frac{|x-x^*|^2}{T-t} + \frac{i}{w^2(T-t)} + i\vartheta},
\end{align}
where $w>0$, $c,x^* \in \bbr^d$ and $\vartheta \in \bbr$.
Both dynamics are closely related to each other in
the pseudo-conformal space
$\Sigma :=\{u\in H^1: xu\in L^2\}$,
through the
{\it pseudo-conformal transform}
\begin{align} \label{pseu-conf-transf}
   S_T(t,x) = \mathcal{C}_T(W)(t,x):= \frac{1}{(T-t)^{\frac d2}} W \(\frac{1}{T-t}, \frac{x}{T-t}\) e^{-i\frac{|x|^2}{4(T-t)}}, \ \ t\not =T,\ x^* =c.
\end{align}
Note that, $S_T$ blows up at time $T$,
and $x^*$ is the singularity
corresponding to the velocity $c$ of $W$.
A remarkable result in the seminal paper by Merle \cite{M93}
is that,
the pseudo-conformal blow-up solution is the unique critical mass blow-up solution
to $L^2$-critical NLS, up to symmetries of the equation.

In the small supercritical mass regime,
two different kinds of blow-up solutions to NLS are exhibited.
The first one is the Bourgain-Wang solution
behaving asymptotically as a sum of a
singular profile $S_T$ and a regular profile $z$, i.e.,
\begin{align}
   v(t) -S_T(t) -z(t) \to 0,\ \ as\ t\ \to T.
\end{align}
Note that, $v$ blows up at time $T$
with the pseudo-conformal speed
$$
   \|\na v(t)\|_{L^2} \sim (T-t)^{-1}.
$$
This kind of solutions was first constructed in the pioneering work
by Bourgain and Wang \cite{BW97} in dimensions $d=1,2$.
It was then extended by Krieger and Schlag \cite{KS09} to prove the existence of
a large set of initial data close to the ground state resulting
in pseudo-conformal speed blow-up solutions in dimension $d=1$,
this set is a codimension one stable manifold in the measurable category.
Moreover,
the instability of such solutions
was proved in the work by Merle, Rapha\"el and Szeftel \cite{MRS13},
which shows that Bourgain-Wang solutions lie on the boundary
of two  $H^1$ open sets of global scattering solutions
and loglog blow-up solutions.
We also would like to refer to \cite{KS06,S09}
for the stable manifolds for the supercritical NLS,
and \cite{B12} for the center-stable manifold
for the $\dot{H}^\frac 12$-critical cubic NLS in dimension three.

Another important kind of blow-up solutions
is of loglog blow-up rate
$$
\|\na v(t)\|_{L^2} \sim
   \((T-t)^{-1}{\log|\log (T-t)|}\)^{\frac 12}.
$$
Unlike Bourgain-Wang solutions, these solutions
are stable under $H^1$ perturbations.
In this respect,
we refer to the pioneering work by Perelman {\cite{P01}}
and  a series of works of Merle and Rapha\"el \cite{MR03,MR04,MR05.2,MR06}.

In the even larger mass regime,
the construction of multi-bubble blow-up solutions
was initiated by Merle \cite{M90},
which behave like a sum of $K$ pseudo-conformal blow-up solutions,
$1\leq K<\9$.
Through the pseudo-conformal transform,
this also yields the existence of multi-solitons, \cite{M90}.
Multi-bubble blow-up solutions with loglog speed have recently been
constructed by Fan \cite{F17}.

For general blow-up solutions to $L^2$-critical NLS,
it is conjectured that
the mass of blow-up solutions is quantized at each singularity
and the remaining part of solutions converges strongly to a residue away from the singularities,
see the {\it mass quantization conjecture} in \cite{MR05},
see also \cite{BW97}.

Let us also mention that,
according to the famous {\it soliton resolution conjecture},
global solutions to a nonlinear dispersive equation are expected to decompose at large time
as a sum of solitons plus a scattering remainder.
We refer to \cite{CKLS18,DJKM17,DKM13,DKM19}
and  references therein for the important progress for the energy critical wave equation.
For the NLS, except for the integrable one dimensional case, this conjecture is still open.
A series of (pure) multi-solitons
(i.e., solutions behaving as a sum of solitons without dispersive part)
have been constructed  for the NLS,
see e.g.  \cite{CF20,CL11,CMM11,LeLP15,LeT14,MM06,M90}.
See \cite{K-M-R} for the construction of
two soliton solutions for the subcritical Hartree equation.
For  the gKdV equations,
we also refer to \cite{Ma05,Co11} for the existence and classification of multi-solitions,
and \cite{C06,C07}
for the construction of solutions behaving as a sum of solitons and of a linear term.

Hence,  a natural question to ask is whether non-pure multi-solitons (including a dispersive part)
can be constructed for the NLS,
which, to the best of our knowledge,
seems not to have been done in literature.
See, e.g., the recent lecture notes of Cazenave \cite{C20}.

The two conjectured long time dynamics are indeed the main motivations
of the present work.

Furthermore,
in the stochastic case,
a remarkable result proved by
de Bouard and Debussche \cite{BD02,BD05} is that,
stochastic solutions can blow up at any short time with positive probability
in the $L^2$-supercritical case.
Several numerical experiments have been also made to investigate
the dynamics of stochastic blow-up solutions,
see, e.g., \cite{BDM02,DL02,DL02.2,MRRY20,MRY20}.

One major challenge in the stochastic case is that,
in contrast to NLS,
the classical pseudo-conformal symmetry is lost due to the input of noise.
Moreover,
the energy of solutions is no longer conserved,
which makes it more difficult to
understand the global behavior in the stochastic $L^2$-supercritical case,
see \cite{MRRY20,MRY20} for the numerical tracking of energy.

Recently, the quantitative construction of
critical mass stochastic blow-up solutions to \eqref{equa-SNLS}
is obtained in \cite{SZ19},
the proof there relies mainly on the modulation method developed in the work by Rapha\"el and Szeftel \cite{RS11}
and also on the rescaling approach in \cite{BRZ16,BRZ16.1,BRZ18,HRZ18,Z17,Z19}.
This  also
yields the threshold of the mass of the ground state for the global well-posedness
and blow-up in the stochastic case.
Later, stochastic blow-up solutions with  loglog speed
have been constructed in \cite{FSZ20}.
Furthermore,
multi-bubble blow-up solutions to \eqref{equa-SNLS},
behaving as a sum of pseudo-conformal blow-up solutions,
were constructed and proved to be unique if the asymptotic behavior
is of the order $(T-t)^{3+}$, \cite{SZ20}.
The conditional uniqueness result has been further used in the very recent work \cite{CSZ21}
to enlarge the energy class for the uniqueness of
both multi-bubble solutions and multi-solitons,
particularly in the low asymptotical regime
with the orders $\calo(T-t)^{0+}$ and $s^{-2-}$,
respectively,
where $t$ is close to $T$ and $s$ is large.

In the present work,
we study the
Bourgain-Wang type solutions,
concentrating at multiple points, in the large mass regime
for both equations \eqref{equa-NLS} and \eqref{equa-SNLS}
in a uniform manner.

More precisely,
in both dimensions one and two,
we construct  multi-bubble Bourgain-Wang solutions to \eqref{equa-NLS-perturb},
which  behave asymptotically as a sum of pseudo-conformal blow-up solutions
and a regular profile, i.e.,
for $t$ close to $T$,
\begin{align} \label{v-Sk-H1-intro}
\|v(t) -\sum_{k=1}^KS_k(t) - z(t) \|_{L^2}
+ (T-t) \|\na v(t) - \na \sum_{k=1}^KS_k(t) - \na z(t) \|_{L^2}
\leq C(T-t)^{\frac 12 (\kappa-1)},
\end{align}
where $z$ is the regular profile propagating along the flow generated by equation \eqref{equa-NLS-perturb}
with $z(T)=z^*$,
$\{S_k\}$ are the pseudo-conformal blow-up solutions as in \eqref{S-blowup-intro}
with distinct singularities,
and the exponent $\kappa\ (\geq 3)$ is closely related to the flatness  at singularities of
both the spatial functions $\{\phi_l\}$ and the residue $z^*$.
Moreover,
we prove that the multi-bubble Bourgain-Wang solutions are unique
if their asymptotical behavior is within the order  $(T-t)^{4+}$.

This provides examples of the conjectured mass quantization phenomena
for both the $L^2$-critical NLS and SNLS.
Furthermore, in the NLS case,
through the pseudo-conformal transform,
the existence and conditional uniqueness of non-pure multi-solutions
are also obtained,
which behave asymptotically as a sum of solitons with distinct velocities
plus a dispersive part.
To the best of our knowledge,
this provides the first  examples of non-pure multi-solitons to
the $L^2$-critical NLS,
predicted by the soliton resolution conjecture.
Let us also mention that,
the uniqueness holds in the energy class of solutions
with decay rate $t^{-5-}$, where $t$ is large enough,
which is larger than the class of exponential convergence
in which (pure) multi-solitons naturally lie. \\

{\bf Notations.}
For any $x=(x_1,\cdots,x_d) \in \bbr^d$
and any multi-index $\upsilon=(\upsilon_1,\cdots, \upsilon_d)$,
let
$|\upsilon|= \sum_{j=1}^d \upsilon_j$,
$\<x\>=(1+|x|^2)^{1/2}$,
$\partial_x^\upsilon=\partial_{x_1}^{\upsilon_1}\cdots \partial_{x_d}^{\upsilon_d}$,
and
$\<\na\>=(I-\Delta)^{1/2}$.

We use the standard Sobolev spaces
$H^{s,p}(\bbr^d)$, $s\in \bbr, 1\leq p\leq\9$.
In particular,
$L^p := H^{0,p}(\bbr^d)$ is
the space of $p$-integrable (complex-valued) functions,
$L^2$ denotes the Hilbert space   with the inner product
$\<v,w\> =\int_{\bbr^d} v(x) \ol w(x)dx$,
and $H^s:= H^{s,2}$.
Let $\Sigma$ denote the pseudo-conformal space,
i.e., $\Sigma:=\{u\in H^1, x u\in L^2\}$.
The local smoothing space is defined by
$L^2(I;H^\a_{\beta})=\{u\in \mathscr{S}': \int_{I} \int \<x\>^{2\beta}|\<\na\>^{\a} u(t,x)|^2  dxdt <\9 \}$,
$\a, \beta \in \mathbb{R}$.
Let $C_c^\9$ be the space of all compactly supported smooth functions on $\bbr^d$.

We also use the notation $\dot{g} = \frac{d}{dt}g$ for any $C^1$
function $g$ on $\bbr$.
For any H\"older continuous function $g\in C^\a(I)$,  $\a>0$, $I\subseteq \bbr^+$,
let $\delta g_{st} := g(t)-g(s)$, $s,t\in I$,
and $\|g\|_{\a, I} := \sup_{s,t\in I,s\not =t} \frac{|\delta g_{st}|}{|s-t|^\a}$.
As $t\to T$ or $t\to \9$,
$f(t)=\calo(g(t))$
means that $|f(t)/g(t)|$ stays bounded,
and $f(t)=o(g(t))$ means that $|f(t)/g(t)|$ converges to zero.

Throughout this paper,
the positive constants $C$ and $\delta$ may change from line to line.

\subsection{Formulation of main results}
Let $K\in \mathbb{N}^+$ and
$\{x_k\}_{k=1}^K$ denote the distinct blow-up points in $\bbr^d$.
We assume that the spatial functions $\{\phi_l\}$ in the noise
and the residue $z^*$ satisfy the following hypotheses:

\begin{enumerate}
   \item[(H1)] {\it Asymptotical flatness}:
  For any multi-index $\upsilon \not = 0$ and $1\leq l\leq N$,
\begin{align} \label{decay}
   \lim_{|x|\to \9} \<x\>^2 |\partial_x^\upsilon  \phi_l(x)| =0.
\end{align}
  {\it Flatness at singularities}:
  Let $\upsilon_*\in \mathbb{N}^+$.
  For every $1\leq l\leq N$ and multi-index $|\upsilon |\leq \upsilon_*$,
\begin{align} \label{degeneracy}
   \partial_x^\upsilon  \phi_l(x_k)=0,\ \ 1\leq l\leq N,\ 1\leq k\leq K.
\end{align}
  \item[(H2)] {\it Smallness}:
  Let $\a^*$ be a positive (small) constant, $m\in \mathbb{N}^+$.
  Let   $z^*$ satisfy
  \begin{align}
    & \|z^*\|_{H^{2m+2+d}}  \leq \a^*,  \label{z*-H-small} \\
    &  \|\<x\> z^*\|_{H^1}  \leq \a^*. \label{z*-Sigma-small}
  \end{align}
  {\it Flatness at singularities}:
  For any multi-index $|\upsilon|\leq 2m$,
  \begin{align}
    & \partial_x^\upsilon z^*(x_k) =0, \ \ 1\leq k\leq K. \label{z*-flat}
  \end{align}
\end{enumerate}

\begin{remark}
On one hand,
the asymptotical flatness condition \eqref{decay} ensures
the Strichartz and local smoothing estimates for the Laplacian with lower order perturbations,
which guarantees the local solvability of equation \eqref{equa-NLS-perturb},
see \cite{BRZ16,MMT08,Z17}.
On the other hand,
the flatness at singularities \eqref{degeneracy} and \eqref{z*-flat}
permit to construct blow-up solutions,
which reflect the local nature of the singularities.
\end{remark}

The main result of this paper is formulated in the following theorem.

\begin{theorem} \label{Thm-BW-RNLS}
Consider equation \eqref{equa-NLS-perturb} with $d=1,2$.
Let $K\in \mathbb{N}^+$, $T\in \bbr_{>0}$, $\{\vartheta_k\}_{k=1}^K \subseteq \bbr$.
Assume that $\{\phi_l\}_{l=1}^N$ and $z^*$ satisfy
Hypotheses $(H1)$  and $(H2)$, respectively, with
$\upsilon_*\geq 5$, $m\geq 3$ if $d=2$ and $m\geq 4$ if $d=1$.

Then, for any distinct points $\{x_k\}_{k=1}^K\subseteq \bbr^d$, $w>0$
(resp. any $\{w_k\}_{k=1}^K \subseteq \bbr_{>0}$),
there exists $\ve^* >0$ small enough such that
for any $\a^*,\ve \in (0,\ve^*)$ and
for any $\{w_k\}_{k=1}^K \subseteq \bbr_{>0}$ with $|w_k - w|\leq \ve$, $1\leq k\leq K$
(resp. any $\{x_k\}_{k=1}^K\subseteq \bbr^d$ with $|x_k-x_j|\geq \ve^{-1}$, $j\not =k$),
the following holds:

$(i)$ {\it Existence.}
There exists a solution $v$ to \eqref{equa-NLS-perturb}
satisfying that for $t$ close to $T$,
\begin{align} \label{v-Sj-L2-RNLS}
\|v(t) -\sum_{k=1}^KS_k(t) - z(t) \|_{L^2}\leq C(T-t)^{\frac {1}{2}(\kappa -1)},
\end{align}
and
\begin{align} \label{v-Sj-Sigma-RNLS}
   \|v(t) -\sum_{k=1}^KS_k(t) - z(t) \|_{\Sigma}\leq C(T-t)^{\frac {1}{2}(\kappa-3)},
\end{align}
where  $\kappa:= (m+\frac d2-1)\wedge (\upsilon_*-2)$,
$C>0$,
$\{S_k\}$ are the pseudo-conformal blow-up solutions of the form
\begin{align}  \label{Sj-blowup}
S_{k}(t,x)=(w_k(T-t))^{-\frac d2}Q \(\frac{x-x_k}{w_k(T-t)}\)
             e^{ - \frac i 4 \frac{|x-x_k|^2}{T-t} + \frac{i}{w_k^2(T-t)} + i\vartheta_k}.
\end{align}
and $z$ is the unique solution of the equation
\be    \label{equa-z}
\left\{ \begin{aligned}
 & i \partial_t z + \Delta z + a_1 \cdot \nabla z+a_0 z +|z|^{\frac 4d} z =0,   \\
    & z(T) = z^*,
\end{aligned}\right.
\ee
where the coefficients $a_1, a_0$ are given by \eqref{a1-loworder} and \eqref{a0-loworder}, respectively.

$(ii)$ {\it Conditional uniqueness.}
Assume in addition that $m\geq 10, \upsilon_*\geq 12$.
Then,
for any small $\zeta >0$,
there exists a unique solution $v$ to \eqref{equa-NLS-perturb} satisfying
that for $t$ close to $T$,
\begin{align} \label{v-Sj-Uniq-RNLS}
\|v(t) -\sum_{k=1}^KS_k(t) - z(t) \|_{L^2}
+ (T-t)\|\na v(t) -\na \sum_{k=1}^KS_k(t) - \na z(t) \|_{L^2}\leq C(T-t)^{4+\zeta}.
\end{align}
\end{theorem}

\begin{remark}
$(i)$. Theorem \ref{Thm-BW-RNLS} mainly treats two cases of  singularities
$\{x_k\}$ and frequencies $\{w_k\}$:

{\rm \bf Case (I)}:
$\{x_k\}_{k=1}^K$ are arbitrary distinct points in $\bbr^d$,
and $\{w_k\}_{k=1}^K (\subseteq \bbr_{>0})$ satisfy that
for some $w>0$, $|w_k - w| \leq \ve$ for every $1\leq k \leq K$.

{\rm \bf Case (II)}:
$\{w_k\}_{k=1}^K$ are arbitrary points in $\bbr_{>0}$,
and $\{x_k\}_{k=1}^K (\subseteq \bbr^d)$ satisfy that
$|x_j-x_k| \geq \ve^{-1}$ for any $1\leq j \neq k \leq K$.

{\rm Cases (I)} and {\rm (II)} roughly mean
certain decoupling between the profiles.
In particular,
{\rm Case (I)} allows the arbitrariness of singularities
when the frequencies are the same.
In the special single bubble case,
both the singularity  and the frequency  can be arbitrary.
Unlike in {\rm Case (II)},
the arbitrariness of singularities in {\rm Case (I)}
is mainly due to the conservation law of mass,
which gives a rapid exponential decay of the sum of the localized masses.

$(ii)$.
The decay order in \eqref{v-Sj-L2-RNLS} and \eqref{v-Sj-Sigma-RNLS}
is closely related to the flatness of $\{\phi_l\}$ and $z^*$ at the singularities.
For $\kappa \geq 4$,
the asymptotics hold in the more regular $H^\frac 32$ space.

$(iii)$. It is important that the regular profile $z$ propagates along the flow generated by equation \eqref{equa-NLS-perturb}.
This fact permits to control the energy,
particularly in the absence of the conservation law of energy,
and to gain one more smallness of the remainder to fulfill the bootstrap arguments in the construction.
The solvability of equation \eqref{equa-z} can be guaranteed by
the smallness of $z^*$ in the Sobolev space and
the Strichartz and local smoothing estimates for the Laplacian with
lower order perturbations (see, e.g., \cite{BRZ14,BRZ16,MMT08,Z17}).

$(iv)$.
The conditional uniqueness  reflects certain rigidity of the flow around
multi-bubble pseudo-conformal blow-up solutions and the regular profile.
It was first proved by Merle, Rapha\"el and Szeftel \cite{MRS13}
in the single bubble case (i.e., $K=1$)
to ensure the continuity of the one-parameter curve in the instability result, \cite{MRS13}.
See also \cite{KK20} for Chern-Simons-Schr\"odinger equations
and \cite{SZ20} for the SNLS case.
It would be very interesting to prove the uniqueness in the low asymptotic regime,
e.g. $(T-t)^{0+}$, as in the very recent work \cite{CSZ21}.
The main challenge here lies in the
linear terms of the remainder in the control of localized mass and energy,
which destroy the upgradation procedure in \cite{CSZ21}.
Let us also mention that,
the proof of the conditional uniqueness in Theorem \ref{Thm-BW-RNLS}
is mainly due to the monotonicity of the generalized energy,
which, fortunately, is stable under the effect of the regular profile.
\end{remark}

The applications to the NLS and SNLS cases are presented below.  \\

\paragraph{{\bf Application 1: The NLS case.}}
One main outcome of  Theorem \ref{Thm-BW-RNLS} is the following theorem
concerning multi-bubble Bourgain-Wang solutions to $L^2$-critical NLS.

\begin{theorem} (Multi-bubble Bourgain-Wang solutions to NLS)  \label{Thm-BW-NLS}
Consider equation \eqref{equa-NLS} with $d=1,2$.
Let $K\in \mathbb{N}^+$, $T\in \bbr_{>0}$, $\{\vartheta_k\}_{k=1}^K \subseteq \bbr$.
Assume that $z^* \in H^{2m+2+d}$ satisfying Hypothesis $(H2)$
with $m\geq 3$ if $d=2$ and $m\geq 4$ if $d=1$.

Then, for any distinct points $\{x_k\}_{k=1}^K \subseteq \bbr^d$, $w>0$
(resp. any $\{w_k\}_{k=1}^K \subseteq \bbr_{>0}$),
there exists $\ve^* >0$ small enough such that
for any $\a^*,\ve \in (0,\ve^*)$ and
for any $\{w_k\}_{k=1}^K \subseteq \bbr_{>0}$ with $|w_k - w|\leq \ve$, $1\leq k\leq K$
(resp. any $\{x_k\}_{k=1}^K \subseteq \bbr^d$ with $|x_k-x_j|\geq \ve^{-1}$, $j\not =k$),
there exists a solution $v$ to \eqref{equa-NLS}
satisfying the asymptotics \eqref{v-Sj-L2-RNLS} and \eqref{v-Sj-Sigma-RNLS},
where the regular profile $z$ is the unique solution of equation
\be     \label{equa-z-NLS}
\left\{ \begin{aligned}
  & i \partial_t z + \Delta z+|z|^{\frac 4d} z =0,   \\
    & z(T) = z^*.
\end{aligned}\right.
\ee
Moreover, if in addition $m\geq 10$,
then for any arbitrarily small $\zeta>0$,
there exists a unique solution to \eqref{equa-NLS} satisfying
the asymptotic \eqref{v-Sj-Uniq-RNLS}.
\end{theorem}

\begin{remark}
Theorem \ref{Thm-BW-NLS} gives  examples
for the conjectured mass quantization in \cite{MR05}.
Actually, by the asymptotical behavior \eqref{v-Sj-L2-RNLS},
for any $R>0$,
\begin{align*}
     v(t) \to z^*\ \ in\ L^2\(\bbr^d - \bigcup\limits_{k=1}^K B(x_k, R)\),
\end{align*}
and
\begin{align*}
    |v(t)|^2 \rightharpoonup \sum\limits_{k=1}^K \|Q\|_{L^2}^2 \delta_{x_k} + |z^*|^2,\ \ as\ t\to T.
\end{align*}
Hence, the solutions concentrate the mass $\|Q\|_{L^2}^2$ at each singularity
and the remaining part converges to a regular residue $z^*$.
\end{remark}

The next result is concerned with the non-pure multi-solitons to $L^2$-critical NLS,
thanks to the pseudo-conformal transform which connects blow-up solutions and solitons.

\begin{theorem} (Non-pure multi-solitons to NLS) \label{Thm-nonpure-soliton-NLS}
Consider equation \eqref{equa-NLS} with $d=1,2$.
Let $K\in \mathbb{N}^+$,  $\{\vartheta_k\}_{k=1}^K \subseteq \bbr$.
Assume that $z^* \in H^{2m+2+d}$ satisfying Hypothesis $(H2)$
with $m\geq 6$.

Then, for any distinct speeds $\{c_k\}_{k=1}^K \subseteq \bbr^d$, $w>0$
(resp. any $\{w_k\}_{k=1}^K  \subseteq \bbr_{>0}$),
there exists $\ve^* >0$ small enough such that
for any $\a^*,\ve \in (0,\ve^*)$ and
for any $\{w_k\}_{k=1}^K  \subseteq \bbr_{>0}$ with $|w_k - w|\leq \ve$, $1\leq k\leq K$
(resp. any $\{c_k\}_{k=1}^K \subseteq \bbr^d$ with $|c_j-c_k|\geq \ve^{-1}$, $j\not =k$),
the following holds:

$(i)$ Existence.
There exists a solution $u$ to \eqref{equa-NLS} satisfying
\begin{align}   \label{u-Wk-NLS}
   \|u(t) - \sum\limits_{k=1}^K W_k(t)-\widetilde{z}(t)\|_{\Sigma} \leq Ct^{-\frac 12 \kappa + \frac 52},\ \ for\ t\ large\ enough,
\end{align}
where $\kappa = m+ \frac d2-1$, $C>0$,
$\{W_k\}$ are the solitary waves to \eqref{equa-NLS} of form
\begin{align}  \label{Wj-soliton}
W_k(t,x) =w_k^{-\frac d2}Q \(\frac{x-c_k t}{w_k} \)e^{i(\half c_k\cdot x-\frac{1}{4}|c_k|^2t+ w_k^{-2}t+\vartheta_{k})},
\end{align}
and $\wt z$ corresponds to the regular part $z$ for \eqref{equa-z-NLS}
through the inverse of the pseudo-conformal transform:
\begin{align}
   \wt z(t,x)  = \mathcal{C}_T^{-1} z (t,x )= t^{-\frac d2} z\(T-\frac 1t, \frac x t\) e^{i\frac{|x|^2}{4t}}.
\end{align}

$(ii)$ Conditional uniqueness.
If in addition $m\geq 16$,
then for any arbitrarily small $\zeta>0$,
there exists a unique non-pure multi-soliton $u$ to \eqref{equa-NLS}
satisfying
\begin{align}
    \|u(t) - \sum\limits_{k=1}^K W_k(t)-\widetilde{z}(t)\|_{\Sigma} \leq Ct^{-5-\zeta},\ \ for\ t\ large\ enough.
\end{align}

\end{theorem}

\begin{remark}
$(i)$. It is known (\cite{D15}) that in the subcritical mass regime
$\|z^*\|_{L^2} < \|Q\|_{L^2}$,
the solution $z$ to \eqref{equa-z-NLS} scatters both forward and backward in time, i.e.,
$\|z\|_{L^{2+\frac{4}{d}}(\bbr\times \bbr^d)} <\9$.
Since the pseudo-conformal transform leaves the $L^2$-critical NLS and the
$L^{2+\frac{4}{d}}(\bbr\times \bbr^d)$-norm invariant,
$\wt z$ also scatters both forward and backward in time with small data $\|z^*\|_{L^2} \leq \a^* <<1$.
Hence, by the asymptotics \eqref{u-Wk-NLS},
the constructed solution
behaves as a sum of solitons plus a dispersive part.
In particular, Theorem \ref{Thm-nonpure-soliton-NLS}
provides new examples of non-pure multi-solitons
to $L^2$-critical NLS,
predicted by the soliton resolution conjecture.

$(ii)$. It is also interesting to see that,
the uniqueness of non-pure multi-solitons holds in the energy class
of solutions with decay rate $t^{-5-}$,
which is much larger than the class of exponential convergence
in which multi-solitons naturally lie (see, e.g, \cite{LeLP15,LeT14}).
We also refer to \cite{CF20,CSZ21}
for this kind of uniqueness in the case of pure multi-solitons
to the $L^2$-critical NLS.
It remains still open to prove the uniqueness or classification of even pure multi-solitons
for the NLS, as done for the gKdV equations in \cite{Ma05,Co11}.

$(iii)$. The relationship between the exponent $m$ and the decay orders in
Theorems \ref{Thm-BW-NLS} and \ref{Thm-nonpure-soliton-NLS}
can be seen from the following estimates:
for $v:= \mathcal{C}_T u$,
\begin{align}
  & \|u(t)\|_{\Sigma} \leq C t \|v(T-\frac 1t)\|_{\Sigma}, \\
  & \|v(t)\|_{\Sigma} \leq  \frac{C}{T-t} \|u(\frac{1}{T-t})\|_{\Sigma}.
\end{align}
\end{remark}

\paragraph{{\bf Application 2: The SNLS case.}}

Another important outcome of Theorem \ref{Thm-BW-RNLS} is in stochastic case.
Let us present the  precise definition of solutions to equation \eqref{equa-SNLS}
in the controlled rough path sense.
For more details of the theory of (controlled) rough paths,
we  refer the interested readers to the monograph \cite{FH14} and \cite{G04}.

\begin{definition} \label{def-X-rough}
We say that $X$ is a solution to \eqref{equa-SNLS} on $[0,\tau^*)$,
where $\tau^*\in (0,\9]$ is a random variable,
if $\bbp$-a.s. for any $\vf\in C_c^\9$,
$t \mapsto \<X(t), \vf\>$ is continuous on $[0,\tau^*)$
and for any $0<s<t<\tau^*$
\begin{align*}
   \<X(t)-X(s), \vf\>
   - \int_s^t \<i X, \Delta\vf\>  + \<i|X|^{\frac 4d} X, \vf\>  - \< \mu X, \vf\> dr
   = \sum\limits_{k=1}^N \int_s^t \<i\phi_k X, \vf\> dB_k(r).
\end{align*}
Here,
the integral $\int_s^t \<i\phi_k X, \vf\> d B_k(r)$
is taken in the sense of controlled rough paths
with respect to the rough paths $(B, \mathbb{B})$,
where $\mathbb{B}=(\mathbb{B}_{jk})$, $\mathbb{B}_{jk,st}:= \int_s^t \delta B_{j,sr} dB_k(r)$
with the integration taken in the sense of It\^o and
$\delta B_{j,st} = B_j(t) - B_j(s)$.
That is,
$\<i\phi_k X, \vf\> \in C^\a([s,t])$,
\begin{align} \label{phikX-st}
   \delta (\<i\phi_k X, \vf\>)_{st}
   = -\sum\limits_{j=1}^N \<\phi_j\phi_k X(s), \vf\> \delta B_{j,st}
     + \delta R_{k,st},
\end{align}
and
$\|\<\phi_j\phi_k X, \vf\> \|_{\a, [s,t]} <\9, \ \
   \|R_k\|_{2\a, [s,t]} <\9
$, where $\frac 13 <\a <\frac 12$.
\end{definition}

The important fact is that,
via the Doss-Sussman type transform
\begin{align}
  v = e^{-W} X,
\end{align}
the $H^1$ solvability of equations \eqref{equa-NLS-perturb} and \eqref{equa-SNLS}
is equivalent, see Theorem 2.10 in \cite{SZ19}.
Thus, by virtue of Theorem \ref{Thm-BW-RNLS},
we have the following result for the $L^2$-critical SNLS.

\begin{theorem} (Multi-bubble Bourgain-Wang solutions to SNLS)   \label{Thm-BW-SNLS}
Consider  \eqref{equa-SNLS} with $d=1,2$.
Let $K\in \mathbb{N}^+$, $\{\vartheta_k\}_{k=1}^K \subseteq \bbr$.
Assume that $\{\phi_l\}_{l=1}^N$ and $z^*$ satisfy
Hypotheses $(H1)$  and $(H2)$, respectively, with
$\upsilon_*\geq 5$,  $m\geq 3$ if $d=2$ and $m\geq 4$ if $d=1$.

Then, for $\bbp$-a.e $\omega\in \Omega$ and
for any distinct points $\{x_k\}_{k=1}^K \subseteq \bbr^d$, $w>0$
(resp. any $\{w_k\}_{k=1}^K \subseteq \bbr_{>0}$),
there exists $\ve^*(\omega) >0$ small enough such that
for any $\a^*,\ve \in (0,\ve^*)$
and any $\{w_k\}_{k=1}^K \subseteq \bbr_{>0}$
with $|w_k - w|\leq \ve$, $1\leq k\leq K$
(resp. any $\{x_k\}_{k=1}^K \subseteq \bbr^d$
with $|x_k-x_j|\geq \ve^{-1}$, $j\not =k$),
the following holds:

There exists $\tau^*(\omega)$ small enough
such that for any $T\in (0,\tau^*(\omega))$,
there exists a solution $X$ to \eqref{equa-SNLS}
satisfying for $t$ close to $T$,
\begin{align} \label{X-Sj-BW-L2-SNLS}
\|e^{-W(t,\omega)} X(t,\omega)-\sum_{k=1}^KS_k(t) - z(t) \|_{L^2}\leq C(T-t)^{\frac 12 (\kappa-1)},
\end{align}
and
\begin{align} \label{X-Sj-BW-SNLS}
\|e^{-W(t,\omega)} X(t,\omega)-\sum_{k=1}^KS_k(t) - z(t) \|_{\Sigma}\leq C(T-t)^{\frac 12 (\kappa-3)},
\end{align}
where $\kappa:= (m+\frac d2 -1) \wedge (\upsilon_*-2)$, $C>0$,
$\{S_k\}$ are the pseudo-conformal blow-up solutions as in \eqref{Sj-blowup},
and $z$ solves equation \eqref{equa-z}.

Moreover, if in addition $m\geq 10$ and $\upsilon_*\geq 12$,
then for any arbitrarily small $\zeta>0$ there exists a unique solution $X$ to \eqref{equa-SNLS} such that
\begin{align} \label{X-Sj-Uniq-SNLS}
\|e^{-W(t,\omega)} X(t,\omega)-\sum_{k=1}^KS_k(t) - z(t) \|_{\Sigma}\leq C(T-t)^{4+\zeta},\ \ for\ t\ close\ to\ T.
\end{align}
\end{theorem}

\begin{remark}
The blow-up time  $T\in (0,\tau^*)$ is chosen to be sufficiently small in Theorem \ref{Thm-BW-SNLS}
because the Brownian motions start moving at time zero.
\end{remark}

\paragraph{{\bf Sketch of Proof.}}
The strategy of proof relies mainly on the modulation method developed in the works \cite{RS11,MRS13}
and on  the multi-bubble analysis in  \cite{M90,CSZ21,SZ20}.

The modulation method in \cite{RS11} is very robust to handle the critical mass blow-up
even in the absence of pseudo-conformal symmetry.
It in particular enables us to treat equation \eqref{equa-NLS-perturb} with lower order perturbations
(or, the stochastic equation \eqref{equa-SNLS}).
Moreover,
as exhibited in \cite{MRS13},
it also permits to
construct Bourgain-Wang solutions
as the limit of both the scattering and loglog blow-up solutions,
rather than by the fixed point arguments in \cite{BW97}.
This inspires us to construct multi-bubble Bourgain-Wang solutions
by using compactness arguments in the modulation framework,
involving the backward integration from the singularity.

More precisely,
we first decompose the approximating solution into three profiles
\begin{align} \label{u-dec-intro}
    v(t,x)=U(t,x)+z(t,x) +R(t,x),
\end{align}
where $U,z,R$ are the blow-up profile,
the regular profile and the remainder, respectively,
which satisfy suitable orthogonality conditions
corresponding to the generalized null space of
the linearized operators around the ground state.
See Theorem \ref{Thm-dec-un} for the detailed statements.

Then, the localization analysis  in \cite{CSZ21,SZ20} and flatness conditions permit to
reduce the analysis to an almost critical mass regime,
in which more dynamical tools developed by \cite{RS11} can be employed.
One crucial ingredient here is
the monotonicity of the generalized energy adapted to the multi-bubble case,
which enables us to derive a uniform backwards control of the remainder.
Hence, the desired blow-up solutions can be constructed
by using compactness arguments as in \cite{RS11,SZ20}.

Let us mention that,
different types of interactions are exhibited in the multi-bubble case:
\begin{enumerate}
  \item[$(i)$] Interactions between different blow-up profiles $U_j$ and $U_k$, $j\not= k$.
  This kind of interaction is of exponentially small order (i.e., $e^{-\frac{\delta}{T-t}}$),
  due to the rapid decay of the ground state
  and the distinction of singularities.
  It was  treated in the pioneering work by Merle \cite{M90}
  in the construction of multi-bubble pseudo-conformal solutions to $L^2$-critical NLS.
  See also \cite{SZ20} for the recent treatments in the stochastic case.

  \item[$(ii)$] Interactions between different localized remainders $R_j$ and $R_k$, $j\not= k$.
  Unlike the previous interactions,
  the remainders are of low polynomial type decay orders.
  There is only little knowledge about remainders in the geometrical decomposition.
  Extra cancellations and decays have to be explored from the
  related localization and cut-off functions,
  e.g.,
  to derive the monotonicity of the generalized energy and
  the coercivity of energy,  \cite{SZ20}.

  \item[$(iii)$] Interactions between the blow-up profile $U$ and the remainder $R$.
  The typical interaction of this kind is the localized mass
  $M_k$ defined in \eqref{Mj-def} below.
  It creates no difficulty in the single bubble case,
  as it is of second order $\calo(\|R\|^2_{L^2})$ thanks to the conservation law of mass.
  However,  the conservation law of mass fails for each localized mass in the multi-bubble case.
  Hence, more delicate analysis has to be performed
  to gain enough temporal regularity, \cite{CSZ21,SZ20}.

  \item[$(iv)$] Interactions between the remainder $R$ and the regular profile $z$.
  This kind of interactions are acceptable in the construction procedure,
  as it is at least of the order $\calo(\|R\|_{L^2})$,
  which suffices for the bootstrap arguments.

  \item[$(v)$] Interactions between the blow-up profile $U$ and the regular profile $z$.
  These interactions are treated by using the flatness condition \eqref{z*-flat}.
  It is not difficult in the NLS case,
  as one may use Taylor's expansion and differentiate equation \eqref{equa-NLS} enough times
  to get high temporal and spatial regularity,  \cite{BW97,MRS13}.
  However, this argument is not applicable in the SNLS case,
  since the coefficients $a_1,a_0$ contain
  the rough paths of  Brownian motions of merely temporal regularity $C^{\frac 12-}_t$.
  The key observation here is that,
  when interacting with the blow-up profile $U$,
  the spatial size $|x-x_k|$ is comparable to the temporal size $T-t$.
  This comparability between space and time
  permits to gain high temporal regularity from the spatial regularity of the residue $z^*$,
  and leads to an inductive expansion of solutions
  for which the continuity of coefficients suffices.
\end{enumerate}

Another major difficulty is the failure of
the conservation law of energy for the solutions to \eqref{equa-NLS-perturb}.

Actually,
unlike in the radial case in \cite{MRS13},
two new modulation parameters $\a_k$ and $\beta_k$
need to be introduced in the multi-bubble case,
due to the distinct singularities
and the non-radialness of solutions.
The control of these new parameters   requires certain coercivity type control of energy,
which however is no longer conserved.
It might be tempting to use the variation control of energy as in \cite{SZ19,SZ20}.
However, in the Bourgain-Wang regime under consideration,
extra terms such as $\|z\|_{H^1}$ appear in the evolution formula of energy,
which, unfortunately, give no temporal regularity
and thus is far from sufficient to close the bootstrap arguments in the construction.

The key point here is,  that the temporal regularity can be gained
after subtracting the energy evolutions of
the solutions and the regular profile.
This leads us to introduce the evolution equation \eqref{equa-z} for the regular profile,
rather than the usual NLS.
Similar structural consideration will  also be used in the
controls of localized mass and
of the remainder in the pseudo-conformal space.

The remainder of this paper is organized as follows.
Section \ref{Sec-Gem-Mod} contains the geometrical decomposition
and preliminary estimates for the modulation equations
and different profiles.
Section \ref{Sec-Mass-Energy} is devoted to the
controls of the localized mass, energy,
and to the curial monotonicity of the generalized energy adapted to the multi-bubble case.
Then,
in Section \ref{Sec-Const-BW}
we mainly construct the multi-bubble Bourgain-Wang solutions to \eqref{equa-NLS-perturb}.
The conditional uniqueness result  in Theorem \ref{Thm-BW-RNLS} is then proved in Section \ref{Sec-Uniq-BW}.
At last, preliminaries concerning the linearized operators,
expansion of the nonlinearity
and the technical proof for modulation equations are collected in
Appendix, i.e., Section \ref{Sec-App}.

\section{Geometrical decomposition}  \label{Sec-Gem-Mod}

\subsection{Geometrical decomposition}  \label{Subsec-Geom-Decomp}
For each $1\leq k\leq K$,
define the modulation parameters by
$\mathcal{P}_k:=(\lbb_k,\alpha_k,\beta_k,\gamma_k,\theta_k) \in \mathbb{Y}:= \bbr\times \bbr^d \times \bbr^d \times \bbr \times \bbr$,
where $\lbb_k, \gamma_k, \theta_k \in \bbr$,
$\a_k, \beta_k \in \bbr^d$.
Set
$\calp:= (\calp_1, \cdots, \calp_K) \in \mathbb{Y}^K$.

Given any $K$ distinct blow-up points $\{x_k\}$,
set $P_{k}:=   |\lbb_{k}|+|\a_{k}-x_k| + |\beta_{k}| + |\g_{k}|$,
$1\leq k\leq K$,
and
$P:= \sum_{k=1}^K P_k$.
Let $S(t,x) = \sum_{k=1}^K S_k(t,x)$,
where $\{S_k\}$ are given by \eqref{Sj-blowup}.

\begin{theorem} (Geometrical decomposition)   \label{Thm-dec-un}
Given $T\in \bbr_{>0}$.
Assume that $v\in C([\wt t, T_*]; H^1)$ solves \eqref{equa-NLS-perturb}
and $v(T_*)=S(T_*) + z(T_*)$,
where $T_*<T$.
Then,
for $\a^*$ sufficiently small and for $T_*$ close to $T$,
there exist $t^*<T_*$
and
unique modulation parameters
$\mathcal{P}
\in C^1((t^*,T_*); \mathbb{Y}^K)$,
such that
$u$ admits the geometrical decomposition
\begin{align} \label{u-dec}
    v(t,x)=U(t,x)+z(t,x) +R(t,x),\ \ t\in [t^*, T_*],\ x\in \bbr^d,
\end{align}
where the main blow-up profile
\begin{align}  \label{U-Uj}
    U(t,x)
    =   \sum_{k=1}^{K} U_k(t,x),
\end{align}
with
\begin{align} \label{Uj-Qj-Q}
   U_k(t,x) = \lbb_k(t)^{-\frac d2} Q_{k}(t,\frac{x-\a_k(t)}{\lbb_k(t)}) e^{i\theta_k(t)},  \ \
   Q_{k}(t,y) = Q(y) e^{i\(\beta_k(t)\cdot y - \frac 14 \g_k(t)|y|^2\)},
\end{align}
the regular profile $z$ solves equation \eqref{equa-z},
$R(T_*)=0$,
and the modulation parameters satisfy
\begin{align} \label{PjT*}
   \calp_k(T_*)=(w_{k}(T-T_*), x_k, 0, w^2_k(T-T_*), w_k^{-2} (T-T_*)^{-1} + \vartheta_k), \ \ 1\leq k\leq K.
\end{align}
Moreover,
for each $1\leq k\leq K$, the following orthogonality conditions hold on $[t^*,T_*] $:
\be\ba\label{ortho-cond-Rn-wn}
&{\rm Re}\int (x-\a_{k}) U_{k}(t)\ol{R}(t)dx=0,\ \
{\rm Re} \int |x-\a_{k}|^2 U_{k}(t) \ol{R}(t)dx=0,\\
&{\rm Im}\int \nabla U_{k}(t) \ol{R}(t)dx=0,\ \
{\rm Im}\int \Lambda_k U_{k}(t) \ol{R}(t)dx=0,\ \
{\rm Im}\int \varrho_{k}(t) \ol{R}(t)dx=0,
\ea\ee
where $\Lambda_k = \frac d2 I_d + (x-\a_k) \cdot \na$,
and
\begin{align}
   & \varrho_{k}(t,x)= \lbb(t)_{k}^{-\frac d2}\rho_{k}(t,\frac{x-\a_{k}(t) }{\lbb_k(t) }) e^{i\theta_{k}(t) }
 \ \ with\   \  \rho_{k}(t,y) := \rho(y)^{i\(\beta_{k}(t)\cdot y - \frac 14 \g_{k}(t) |y|^2\)},  \label{rhon}
\end{align}
and $\rho$ is given by \eqref{def-rho}.
\end{theorem}

Theorem \ref{Thm-dec-un} is mainly based on the implicit function theorem.
The case $z^*=0$ is proved in \cite{SZ20}.
Since the smallness condition of $z$ still keeps the non-degeneracy of
the determinant of Jacobian matrix,
the arguments in \cite{SZ20} are also applicable here.
For simplicity, the proof is omitted.

\subsection{Modulation equations} \label{Subsec-Mod-Equa}
Let $\dot{g}:= \frac{d}{dt}g$  for any $C^1$ function $g$.
For each $1\leq k\leq K$,
define \emph{the vector of modulation equations} by
\begin{align} \label{Mod-def}
   Mod_{k}:= |\lambda_{k}\dot{\lambda}_{k}+\gamma_{k}|+|\lambda_{k}^2\dot{\gamma}_{k}+\gamma_{k}^2|
   +|\lambda_{k}\dot{\alpha}_{k}-2\beta_{k}|
                +|\lambda_{k}^2\dot{\beta}_{k}+\gamma_{k}\beta_{k}|
              +|\lambda_{k}^2\dot{\theta}_{k}-1-|\beta_{k}|^2|.
\end{align}
Set $Mod:=\sum_{k=1}^{K}Mod_{k}$.

The modulation equations mainly characterize the dynamics of geometrical parameters.
The main estimate
is contained in Theorem \ref{Thm-Mod-bdd} below.

\begin{theorem} (Control of modulation equations)  \label{Thm-Mod-bdd}
Assume that $u$ admits the geometrical decomposition \eqref{u-dec} on $[t^*,T_*] \subseteq [0,T)$
with the modulation parameters $\calp=(\lbb, \a, \beta, \g, \theta)\in \mathbb{Y}^{K}$.
Assume additionally that for some $C_1,C_2>0$,
\begin{align}
   C_1(T-t)\leq \lbb_k \leq C_2(T-t),\ \ t\in [t^*, T_*], \ 1\leq k\leq K.
\end{align}
Then,  for $t^*$ close to $T$,
there exists $C>0$ such that
for any $ t\in[t^*,T_*]$,
\begin{align} \label{Mod-bdd}
Mod\leq
C\( \sum_{k=1}^{K} |M_k| +  P^2 D + D^2
  + \a^* (T-t)^{m+1+\frac d2} + P^{\upsilon_*+1} \),
\end{align}
where $\upsilon_*$ is the index of flatness in \eqref{degeneracy},
$M_k$ is the localized mass
\begin{align}
   M_k  = 2 {\rm Re} \<R_k, U_k\> + \int |R|^2 \Phi_k dx,
\end{align}
and
$D$ is the important quantity to measure the size of remainder, defined by
\begin{align}   \label{D-def}
    D:= {\|R\|_{L^2}} + (T-t) \|\na R\|_{L^2}.
\end{align}

Moreover, we have the improved estimate
\begin{align} \label{lbb-ga-mod}
   |\lambda_{k}\dot{\lambda}_{k}+\gamma_{k}|
\leq& C\(    P^2 D +  D^2
              + \a^* (T-t)^{m+1+\frac d2}  + P^{\upsilon_*+1}   \).
\end{align}
\end{theorem}

\begin{remark}
By Lemma \ref{Lem-Mod-boot} below,
we shall see that
estimate \eqref{lbb-ga-mod} gains one more fact $T-t$ than \eqref{Mod-bdd},
which is important in the derivation of the monotonicity of generalized energy $\scri$.
\end{remark}

Theorem \ref{Thm-Mod-bdd} can be proved analogously as in \cite{SZ19,SZ20},
and hence is postponed to the Appendix for simplicity.

\subsection{Estimates of profiles}

We collect in this subsection the estimates of three profiles in the above
geometrical decomposition \eqref{u-dec},
which will be frequently used in the squeal.  \\

\paragraph{\bf The blow-up profile $U$.}
Let us first see that,
by the explicit formula \eqref{Uj-Qj-Q},
$U_k$ satisfies the equation
\begin{align}   \label{equa-Ut}
 i\partial_tU_{k}+&\Delta U_{k}+|U_{k}|^{\frac{4}{d}}U_{k}
 = \psi_k
 = \frac{e^{i\theta}}{\lbb_k^{2+\frac d2}} \Psi_k(t,\frac{x-\alpha_{k}}{\lambda_{k}}),
\end{align}
where  $1\leq k\leq K$, and
\begin{align} \label{Psik-def}
   \Psi_k =&   -(\lambda_{k}^2\dot{\theta}_{k}-1-|\beta_{k}|^2)Q_{k}
      -(\lambda_{k}^2\dot{\beta}_{k}+\gamma_{k}\beta_{k})\cdot yQ_{k}
    +\frac 14 (\lambda_{k}^2\dot{\gamma}_{k}+\gamma_{k}^2) |y|^2 Q_{k}  \nonumber \\
    &  -i(\lambda_{k}\dot{\alpha}_{k}-2\beta_{k})\cdot \nabla Q_{k}
     -i(\lambda_{k}\dot{\lambda}_{k}+\gamma_{k})\Lambda Q_{k}.
\end{align}

\begin{lemma}   \label{Lem-U}
Suppose that $P=\calo(1)$ and
$\lbb_k\geq C(T-t)$, $C>0$.
Then, for any $p\geq 2$,
there exists $C>0$
such that  for all $t\in [t^*, T_*]$, $1\leq k\leq K$,
\begin{align}
     \|U(t)\|^p_{L^p} \leq C (T-t)^{-d(\frac p2 -1)}.  \label{U-Lp}
\end{align}
\end{lemma}

{\bf Proof.}
Estimate  \eqref{U-Lp} follows from
the Gagliardo-Nirenberg inequality   that
for any $2\leq p<\9$,
\begin{align}  \label{G-N}
\|g\|_{L^p}\leq C \|g\|_{L^2}^{1 - d (\frac 12-\frac 1p)} \|\nabla g\|_{L^2}^{d(\frac 12-\frac 1p)}, \ \ \forall g\in H^1,
\end{align}
and the estimates
\begin{align} \label{Uj-L2-H1}
   \|U_k(t)\|_{L^2}= \|Q\|_{L^2}, \ \
   \|\nabla U_k(t)\|_{L^2} = \lbb_k^{-1}\|\na Q_k\|_{L^2}
   \leq  C (T-t)^{-1}.
\end{align}
\hfill $\square$

Because the blow-up profile $U_k$ is almost localized around $\frac{x-\a_k}{\lbb_k}$
and the singularities are separated from each other,
the interactions between different blow-up profiles are exponentially small.
Lemma \ref{Lem-decoup-U} below is a slight modification of  \cite[Lemma 3.1]{SZ20}.

\begin{lemma}(Interactions between blow-up profiles)  \label{Lem-decoup-U}
Let $0<t^*<T_*<\9$. For $1\leq k\leq K$,
set
\begin{align} \label{Gj-gj-g}
    G_k(t,x)
    := \lbb_k^{-\frac d2} g_k(t,\frac{x-\a_k}{\lbb_k}) e^{i\theta_k}, \ \
    with\ \ g_{k}(t,y) := g(y) e^{i(\beta_{k}(t) \cdot y - \frac 14\g_{k}(t) |y|^2)},
\end{align}
where $g \in C_b^2(\bbr^d)$
decays exponentially fast at infinity
\begin{align*}
   |\partial^\upsilon g(y)| \leq C e^{-\delta |y|}, \ \ |\upsilon|\leq 2,
\end{align*}
with $C,\delta>0$,
$\calp_k:=(\lbb_k,\alpha_k,\beta_k,\gamma_k,\theta_k) \in C([t^*,T_*]; \mathbb{Y})$
satisfies that for $t\in [t^*,T_*]$,
\begin{align} \label{aj-xj}
  \half\leq\frac{\lambda_{k}(t)}{w_{k}(T-t)}\leq2,\ \
  |\alpha_k(t)-x_k|\leq  \min\limits_{j\not =k}\{\frac {1}{12} |x_j-x_k|\} \wedge \frac 12,\ \
  |\beta_k(t)| + |\g_k(t)| \leq 1,
\end{align}
and
\begin{align} \label{T-M-0}
     (T-t^*) (1+ \max_{1\leq k\leq K} |x_k|) \leq  1.
\end{align}
Then,
there exist $C, \delta>0$ such that
for any $1\leq k\not =l\leq K$,
$m\in \mathbb{N}$ and
multi-index $\upsilon$ with $|\upsilon|\leq 2$,
\begin{align} \label{Gj-Gl-decoup}
  \int\limits_{\bbr^d} |x-\a_l|^n |\partial^\upsilon G_l(t)| |x-\a_k|^m |G_k(t)| dx
   \leq Ce^{-\frac{\delta}{T-t}}, \ \ t\in [t^*,T_*].
\end{align}
Moreover,
let $\Phi_k$ be defined in \eqref{phi-local} below. Then,
for any $h\in L^1$ or $L^2$,
$1\leq k\not = l\leq K$,
$m,n\in \mathbb{N}$ and multi-index $\upsilon$ with $|\upsilon |\leq 2$,
\begin{align}  \label{Gj-hl-decoup}
   \int\limits_{\bbr^d} |x-\a_l|^n |\partial^\upsilon  G_l(t)| |x-\a_k|^m |h| \Phi_kdx
   \leq Ce^{-\frac{\delta}{T-t}} \min\{\|h\|_{L^1}, \|h\|_{L^2}\}, \ \ t\in [t^*,T_*].
\end{align}
\end{lemma}

In the sequel,
we take $t$ and $T^*$ close to $T$ such that \eqref{aj-xj} and \eqref{T-M-0} hold.
In particular, $\lbb_k$ is comparable with $T-t$:
\begin{align} \label{lbb-t-comp}
   \frac 12 w_k (T-t) \leq  \lbb_k(t) \leq 2 w_k (T-t).
\end{align}
Hence, Lemma \ref{Lem-decoup-U} is applicable.  \\

\paragraph{\bf The regular profile $z$.}
The main estimates of  regular profile are contained in the following lemma.

\begin{lemma}  \label{Lem-z}
Let $z^*$ satisfy Hypothesis $(H2)$.
Let $z$ be the corresponding solution to equation \eqref{equa-z}.
For every $1\leq k\leq K$,
define the renormalized variables $\ve_{z,k}$ by
\begin{align} \label{ez-def}
   z(t,x) = \lbb_k^{-\frac d2} \ve_{z,k}(t,\frac{x-\a_k}{\lbb_k})e^{i \theta_k},
\end{align}
Then, for $\a^*=\a^*(T,m)$ sufficiently small, the following estimates hold:

$(i)$ (Smallness.) For $t$ close to $T$,
\begin{align}
   & \|z\|_{L^\9(t,T; H^{2m+d+2})} \leq C_{T,m} \a^*,  \label{z-LtHm-a*}  \\
   & \|\partial_t z\|_{C([t,T]; L^2)} \leq C_T \a^*.  \label{z-CtL2-a*}
\end{align}
In particular,
\begin{align} \label{vez-L2-H1}
  & \| \ve_{z,k}\|_{L^2} \leq C \a^*,\ \ \|\na \ve_{z,k}\|_{L^2} \leq C_T \a^* (T-t).
\end{align}

If in addition $xz\in H^1$,
then for $t$ close to $T$,
\begin{align} \label{xz-H1-small}
   \|xz\|_{L^\9(t,T;H^1)} \leq C_T \a^*.
\end{align}

$(ii)$ (Interaction between the profiles $U$ and $z$.)
If in addition $P(t) = \calo(T-t)$ and
$|\a_k-x_k|<\frac 12$ for $t$ close to $T$,
then for any $\delta>0$, there exists $C_{T,m,\delta}>0$ such that
\begin{align}    \label{eyvez-bdd}
  & \sum\limits_{|\upsilon |\leq 2} \|e^{-\delta|y|} \partial_y^\upsilon  \ve_{z,k} \|_{L^\9}
     \leq C_{T,m,\delta}  \a^* (T-t)^{m+1+\frac d2}.
\end{align}
\end{lemma}

\begin{remark}
Note that, by  the exponential decay \eqref{Q-decay} of ground state,
\begin{align*}
   \int U_k(t,x) \ol{z(t,x)} dx
   = \int Q_k(y)  \ol{ \ve_{z,k} (t,y)} dy
   \leq C \|e^{-\delta|y|} \ve_{z,k}\|_{L^\9}.
\end{align*}
Hence, estimate \eqref{eyvez-bdd} controls the interactions between the blow-up profile
and the regular profile.
As explained in Introduction,
this estimate in the NLS case follows from
Taylor's expansion of $z$ and differentiating equation \eqref{equa-NLS} enough times to
get high temporal orders.
For more general equation \eqref{equa-NLS-perturb},
including particularly the SNLS \eqref{equa-SNLS}
where the coefficients $a_1, a_0$ are only $C^{\frac 12-}_t$-regular in time,
we shall use a different inductive expansion of solutions
and the comparability between the spatial size $|x-\a_k|$ and the temporal size $T-t$,
due to the well localization of blow-up profile $U_k$.
\end{remark}

{\bf Proof of Lemma \ref{Lem-z}.}
$(i)$.
For simplicity, we set $p=2+\frac 4d$, $p'=\frac{2d+4}{d+4}$.
Applying the derivative $\<\na\>^{2m+d+2}$ to both sides of equation \eqref{equa-z} we have
\begin{align}
   i\partial_t \<\na\>^{2m+d+2} z + (\Delta  + a_1\cdot \na + a_0) \<\na\>^{2m+d+2} z
   + [\<\na\>^{2m+d+2}, a_1\cdot \na +a_0]z + \<\na\>^{2m+d+2}(|z|^{\frac 4d}z) =0,
\end{align}
with $\<\na\>^{2m+d+2} z(T) = \<\na\>^{2m+d+2} z^*$,
where $[\<\na\>^{2m+d+2}, a_1\cdot \na +a_0]$ is the commutator
$\<\na\>^{2m+d+2}(a_1\cdot \na +a_0) - (a_1\cdot \na +a_0) \<\na\>^{2m+d+2}$.
Then using the Strichartz and local smoothing estimates
(see \cite[Theorem 2.11]{Z17})
we have
\begin{align}
   \|z\|_{L^\9(t,T;H^{2m+d+2})}
   \leq& C_T \(\| z_*\|_{H^{2m+d+2}} + \| [\<\na\>^{2m+d+2}, a_1\cdot \na +a_0]z\|_{L^2(t,T;H^{-\frac 12}_{1})}
              + \|\<\na\>^{2m+d+2} (|z|^{\frac 4d}z)\||_{L^{p'}(t,T; L^{p'})} \) \nonumber \\
   \leq& C_T \(\a^* + \|z\|_{L^2(t,T;H^{2m+d+\frac 32}_{-1})}
                +  \| z\|_{L^\9(t,T;H^{2m+d+2})}^{1+\frac 4d} \).
\end{align}
Then, using   the interpolation (see \cite[Lemma 3.6]{Z17})
\begin{align}
  \|z\|_{H^{2m+d+\frac 32}_{-1}}
  \leq C \delta^\frac 12 \|z\|_{H^{2m+d+2}_{-1}} + C \delta^{-(2m+d+\frac 32)}\|z\|_{L^2}
   \leq C \delta^\frac 12  \|z\|_{H^{2m+d+2}} +  C \delta^{-(2m+d+\frac 32)}  \a^*,
\end{align}
where the last step is due to $\<x\>^{-1} \leq 1$
and the mass conservation $\|z\|_{L^2} = \|z^*\|_{L^2} \leq \a^*$,
we lead to
\begin{align}
   \|z\|_{L^\9(t,T;H^{2m+d+2})}
    \leq& C_T \((1+T^{\frac 12} \delta^{-(2m+d+\frac 32)} )\a^*
               + T^{\frac 12} \delta^\frac 12  \|z\|_{L^\9(t,T;H^{2m+d+2})}
                +   \| z\|_{L^\9(t,T;H^{2m+d+2})}^{1+\frac 4d}\).
\end{align}
Here and in the sequel, the constant $C_T$ may change from line to line.

Taking $\delta$ small enough such that $C_T T^\frac 12 \delta^\frac 12 <1/2$
we obtain
\begin{align}
   \|z\|_{L^\9(t,T;H^{2m+d+2})}
    \leq& C_{T,m} \(\a^*   +   \| z\|_{L^\9(t,T;H^{2m+d+2})}^{1+\frac 4d}\).
\end{align}
Hence,
taking $\a^*$ small enough we obtain  \eqref{z-LtHm-a*}.

Estimate \eqref{z-CtL2-a*} then follows from \eqref{z-LtHm-a*} and equation \eqref{equa-z},
and estimates in \eqref{vez-L2-H1} follow directly from the identities:
\begin{align} \label{vez-z}
   \ve_{z,k} (t,y)
   = \lbb_k^{\frac d2} z(t,\lbb_ky + \a_k) e^{-i\theta_k},
\end{align}
and
\begin{align} \label{navez-z}
   \na \ve_{z,k}(t,y)
   = \lbb_k^{\frac d2+1}   \na z(t,\lbb_ky + \a_k)   e^{-i\theta_k}.
\end{align}

It remains to prove \eqref{xz-H1-small}.
For this purpose,
we derive from \eqref{equa-z} that,
for every $1\leq j\leq d$,
\begin{align}
   i\partial_t (x_jz) + \Delta (x_jz) + a_1\cdot \na (x_jz) + a_0 (x_jz)
   - 2 \partial_j z - a_{1,j} z + x_j f(z) =0,
\end{align}
and $x_jz(T) =x_jz^*$,
$a_{1,j}$ is the $j$-th component of the vector $a_1$.
Then, applying Strichartz estimates
and using \eqref{z-LtHm-a*} we get
\begin{align}
   \|x_jz\|_{L^\9(t,T; L^2)}
   \leq& C_T \( \|x_jz^*\|_{L^2} +  \|2 \partial_j z + a_{1,j} z - x_j f(z)\|_{L^1(t,T; L^2)} \) \nonumber \\
   \leq& C_T \( \|x_jz^*\|_{L^2}
                + (T-t) \|z\|_{L^\9(t,T; H^1)}
                +  (T-t) \|z\|_{L^\9(t,T; L^\9)}^{\frac 4d} \|x_j z\|_{L^\9(t, T; L^2)}  \) \nonumber  \\
   \leq& C_T \( \a^* + (T-t)\a^* \|x_j z\|_{L^\9(t,T; L^2)} \),
\end{align}
which yields that for $t$ close to $T$,
\begin{align}  \label{xz-CtL2-a*}
    \|x_j z\|_{L^\9(t,T;L^2)} \leq C_T \a^*.
\end{align}

Moreover,
for every $1\leq l\leq d$, $x_j \partial_l z$ satisfies
\begin{align}
  i\partial_t (x_j \partial_l z) +  \Delta  (x_j \partial_l z)   + a_1 \cdot \na  (x_j \partial_l z)    + a_0  (x_j \partial_l z)
  + \mathcal{N}  =0,
\end{align}
where $ x_j \partial_l z(T)= x_j \partial_l z^*$, and
\begin{align}
   \mathcal{N} = -2 \partial_{jl}z - a_{1,j} \partial_l z + x_j (\partial_l a_1)\cdot \na z
      + x_j(\partial_l a_0)z + x_j \partial_l f(z) .
\end{align}
Hence,
by Strichartz estimates, \eqref{decay} and \eqref{z-LtHm-a*},
\begin{align*}
   \|x_j \partial_l z\|_{L^\9(t,T; L^2)}
   \leq& C_T \(\|x_j \partial_l z^*\|_{L^2} + \|\mathcal{N}\|_{L^1(t,T;L^2)}  \) \nonumber \\
   \leq& C_T \( \|x_j \partial_l z^*\|_{L^2} + (T-t)\|z\|_{L^\9(t,T;H^2)}
               + (T-t) \|z\|^{\frac 4d}_{L^\9(t,T; L^\9)} \|x_j\partial_l z\|_{L^\9(t,T; L^2)} \) \nonumber \\
   \leq& C_T \( \a^* + (T-t)\a^* \|x_j\partial_l z\|_{L^\9(t,T; L^2)} \),
\end{align*}
which yields that for $t$ close to $T$,
\begin{align}
   \|x_j \partial_l z\|_{L^\9(t,T; L^2)}  \leq C_T \a^*.
\end{align}
Thus, estimate \eqref{xz-H1-small} is proved.

$(ii)$.
Let us set $m^*:=2m+1$
and define the   operator $\mathcal{D}_t$ by
\begin{align*}
   \mathcal{D}_t := - i\( \Delta + a_1(t)\cdot \na + a_0(t) + |z(t)|^{\frac 4d} \).
\end{align*}
Then, by equation \eqref{equa-NLS-perturb} and the mean valued theorem,
\begin{align} \label{z-z*-t-expan}
   z(t) =  z^* + \int_t^T \mathcal{D}_{r} z(r) dr
        =  z^* + (T-t) \mathcal{D}_{t_1} z(t_1),
\end{align}
where  $t_1\in (t,T)$.
Further  expansion of  $z(t_1)$ by \eqref{equa-NLS-perturb} yields
\begin{align*}
   z(t) =&  z^* + (T-t) \cald_{t_1} \( z^* + \int_{t_1}^T \cald_{r} z(r) dr\) \nonumber \\
        =& z^* + (T-t) \cald_{t_1}  z^*  + (T-t)(T-t_1) \cald_{t_1} \circ \cald_{t_2} z(t_2),
\end{align*}
where $t_2 \in (t_1, T)$.
Then, further expansion by \eqref{z-z*-t-expan} and inductive arguments lead to
\begin{align} \label{z-calf-expan}
  z(t) = z^* + \sum\limits_{j=1}^{n} \prod \limits_{l=0}^{j-1} (T-t_l)  \cald_{t_1}   \circ \cdots \circ \cald_{t_j} z^*
            +  \prod \limits_{l=0}^{n} (T-t_l) \cald_{t_1} \circ \cdots \circ \cald_{t_{n+1}} z(t_{n+1}),
\end{align}
where $t_0:=t$, $t_{l}\in (t,T)$, $1\leq l\leq n$.
By   \eqref{z*-flat},
\begin{align} \label{cald-z*-Taylor}
   | \cald_{t_1}   \circ \cdots \circ \cald_{t_j} z^* (x)|
   \leq C_T |x-x_k|^{m^*-2j},
\end{align}
and by the Sobolev embedding $H^{2m+2+d} \hookrightarrow C^{2(m+1)}_b$,
for $2\leq n\leq m$,
\begin{align} \label{cald-z-L9}
   \|\cald_{t_1} \circ \cdots \circ \cald_{t_{n+1}} z(t_{n+1})\|_{L^\9}
   \leq C_{T,m} \|z\|_{L^\9(t,T;C_b^{2(m+1)})}
   \leq C_{T,m}  \|z\|_{L^\9(t,T; H^{2m+2+d})}.
\end{align}
Hence, we derive that for $2\leq n\leq m$,
\begin{align} \label{z-t-expan}
  |z(t)| \leq C_{T,m} \(|x-x_k|^{m^*} +  \sum\limits_{j=1}^{n}  (T-t)^j |x-x_k|^{m^*-2j} + (T-t)^{n+1}\),
  \ \ for\ |x-x_k|<1.
\end{align}

Note that,
\begin{align} \label{lbby-a-xk}
   |\lbb_k y + \a_k-x_k| \leq C P \<y\> \leq C(T-t)\<y\>.
\end{align}
Moreover, since $|\a_k -x_k|<\frac 12$,
$|y|>\frac{1}{2\lbb_k}$ in the regime
$\{y\in \bbr^d: |\lbb_k y +\a_k - x_k|\geq 1\}$,
and so
\begin{align*}
    \|e^{-\delta |y|} z(t,\lbb_ky + \a_k) I_{|\lbb_k y +\a_k - x_k|\geq 1}\|_{L^\9}
    \leq C_{T,m}e^{-\frac{\delta}{T-t}}.
\end{align*}
Taking into account \eqref{vez-z}, \eqref{z-t-expan} and \eqref{lbby-a-xk}
we obtain
\begin{align} \label{vezk-Tt-esti}
    \|e^{-\delta |y|} \ve_{z,k}(y)\|_{L^\9}
    \leq& \lbb_k^{\frac d2} \|e^{-\delta|y|} z(\lbb_k y+\a_k)
               (I_{|\lbb_k y +\a_k - x_k| < 1} + I_{|\lbb_k y +\a_k - x_k|\geq 1}) \|_{L^\9} \nonumber \\
    \leq& C_{T,m}    (T-t)^\frac d2 \((T-t)^{m^*-n} + (T-t)^{n+1}\)
          + C_{T,m} e^{-\frac{\delta}{T-t}}.
\end{align}
This yields that, for $n=m$ and $t$ close to $T$,
\begin{align}
    \|e^{-\delta |y|} \ve_{z,k}(y)\|_{L^\9} \leq C_{T,m,\delta} (T-t)^{m+1 + \frac d2}.
\end{align}

Similarly,
by \eqref{z-calf-expan},
for any multi-index $|\upsilon |\leq 2$,
\begin{align} \label{z-calf-expan-pnu}
  \partial_x^\upsilon  z(t)
    =& \partial_x^\upsilon  z^*
       + \sum\limits_{j=1}^{n} \prod \limits_{l=0}^{j-1} (T-t_l) \partial_x^\upsilon  \circ \cald_{t_1}   \circ \cdots \circ \cald_{t_j} z^*
            +  \prod \limits_{l=0}^{n} (T-t_l) \partial_x^\upsilon  \circ \cald_{t_1} \circ \cdots \circ \cald_{t_{n+1}} z(t_{n+1}) .
\end{align}
As in \eqref{cald-z*-Taylor} and \eqref{cald-z-L9}, we have
\begin{align} \label{pnu-cald-z*-Taylor}
   | \partial_x^\upsilon  \circ \cald_{t_1}   \circ \cdots \circ \cald_{t_j} z^* (x)|
   \leq C_T |x-x_k|^{m^*-2j-|\upsilon|},
\end{align}
and for $n\leq m-1$,
\begin{align} \label{pnu-cald-z-L9}
   \| \partial_x^\upsilon  \circ \cald_{t_1} \circ \cdots \circ \cald_{t_{n+1}} z(t_{n+1})\|_{L^\9}
     \leq C_{T,m}  \|z\|_{L^\9(t,T; H^{2m+2+d})}.
\end{align}
Thus,  for any multi-index $|\upsilon |\leq 2$ and $n\leq m-|\upsilon |$,
\begin{align}
    | \partial_x^\upsilon  z(t)|
    \leq& C_{T,m} \(|x-x_k|^{m^*-|\upsilon |}
             + \sum\limits_{j=1}^n (T-t)^j |x-x_k|^{m^*-2j-|\upsilon |}
             + (T-t)^{n+1} \), \ \ for\ |x-x_k|<1.
\end{align}
Taking into account
\begin{align} \label{pnu-vez-z}
     \partial_x^\upsilon  \ve_{z,k} (t,y)
   = \lbb_k^{\frac d2+|\upsilon |} \partial_y^\upsilon  z(t,\lbb_ky + \a_k) e^{-i\theta_k}
\end{align}
and arguing as in the proof of \eqref{vezk-Tt-esti}
we get
\begin{align}
   \|e^{-\delta|y|} \partial_y^\upsilon  \ve_{z,k}(t) \|_{L^\9}
   \leq& C_{T,m} (T-t)^{\frac d2 +|\upsilon |} \( (T-t)^{m^*-|\upsilon |} + \sum\limits_{j=1}^n (T-t)^{m^*-j-|\upsilon |} + (T-t)^{n+1}\)
         +  C_{T,m} e^{-\frac{\delta}{T-t}} \nonumber \\
   \leq& C_{T,m} \( (T-t)^{m^*-n+\frac d2} + (T-t)^{n+|\upsilon |+1+\frac d2}
         +   e^{-\frac{\delta}{T-t}}  \) \nonumber \\
   \leq& C_{T,m,\delta} (T-t)^{m+1+\frac d2},
\end{align}
where in the last step we chose $n=m-|\upsilon |$ for $1\leq |\upsilon |\leq 2$.

Therefore, the proof of Lemma \ref{Lem-z} is complete. \hfill $\square$

For the coefficients of  lower order perturbations,
we have the following estimates.

\begin{lemma}  \label{Lem-wtbwtc}
For any multi-index $\upsilon $, $|\upsilon |\leq 2$,
set $\wt{\partial^\upsilon  \phi_{l,k}}(y) := (\partial^\upsilon  \phi_l)(\lbb_k y + \a_k)$,
$1\leq l\leq N$.
Then,
\begin{align}   \label{phik-Taylor}
|\wt{ \partial^\upsilon   \phi_{l,k}} (y)|
 \leq C P^{\upsilon_*+1-|\upsilon |} \<y\>^{\upsilon_*+1}, \ \ 0\leq |\upsilon |\leq \upsilon_*,
\end{align}
where $\upsilon_*$ is the index of flatness in Hypothesis $(H1)$.
In particular,
for $\wt{a}_{1,k}(t,y):=a_1(t,\lbb_{k}y+\alpha_{k})$,
$\wt{a}_{0,k}(t,y):=a_0(t, \lbb_{k}y+\alpha_{k})$,
we have that for any multi-index $\upsilon $, $|\upsilon |\leq 2$,
there exists $C>0$ such that
\begin{align}
   & |\partial^\upsilon _y (\wt a_{1,k} (t,y))| \leq C \lbb_k^{|\upsilon |} P^{\upsilon_* -|\upsilon |} \<y\>^{\upsilon_*+1}, \label{a1-Pnu} \\
   & |\partial^\upsilon _y (\wt a_{0,k} (t,y))| \leq C \lbb_k^{|\upsilon |} P^{\upsilon_*-1 -|\upsilon |} \<y\>^{2\upsilon_*+2}. \label{a0-Pnu}
\end{align}
\end{lemma}

{\bf Proof.}
By Taylor's expansion and \eqref{degeneracy},
\begin{align*}
| \wt{\partial^\upsilon   \phi_{l,k}} (y)|
\leq C(\lambda_{k} y+\alpha_{k}-x_k)^{\upsilon_*+1-|\upsilon |}
\leq C P^{\upsilon_*+1-|\upsilon |} \<y\>^{\upsilon_*+1}, \ \ 0\leq |\upsilon |\leq \upsilon_*.
\end{align*}
This yields \eqref{phik-Taylor}.
Estimates \eqref{a1-Pnu} and \eqref{a0-Pnu} then follow from \eqref{phik-Taylor},
\eqref{a1-loworder} and \eqref{a0-loworder}.
\hfill $\square$  \\

\paragraph{\bf The remainder profile $R$.}
Lemma \ref{Lem-R-D-Lp} below permits to control the $H^1$ and $L^p$-norms of remainder.

\begin{lemma} (\cite[Lemma 2.7]{SZ20})  \label{Lem-R-D-Lp}
There exists $C>0$ such that
\begin{align}
   & \|R\|_{H^1} \leq C (T-t)^{-1} D ,\ \ \|R\|_{L^2} \|\na R\|_{L^2} \leq (T-t)^{-1} D^2 , \label{R-H1} \\
   & \|R\|_{L^p}^p \leq C (T-t)^{-d(\frac p2 -1)} D^p . \label{R-Lp-D}
\end{align}
\end{lemma}

In order to deal with the multi-bubble case,
it is useful to decompose the remainder $R$ into $K$ localized profiles  concentrating at
the singularities.
As in \cite{SZ20},
since equation \eqref{equa-NLS} is invariant under orthogonal transforms,
we may take an orthonormal basis $\{{\bf v_j}\}_{j=1}^d$ of $\bbr^d$,
such that $(x_j-x_l)\cdot {\bf v_1}\neq0$ for any $1\leq j\neq l\leq K$.
Hence, without loss of generality,
we assume that $x_1\cdot {\bf v_1}<x_2\cdot {\bf v_1}<\cdots<x_K\cdot {\bf v_1}$.
Then, set
\begin{align} \label{sep-xj-0}
\sigma :=\frac{1}{12}\min_{1\leq k \leq K-1}\{(x_{k+1}-x_k)\cdot {\bf v_1}\}> 0.
\end{align}
Let $\Phi(x)$ be a smooth function on $\R^d$ such that $0\leq \Phi(x)\leq 1$,
$|\na \Phi(x)| \leq C \sigma^{-1}$,
$\Phi(x)=1$ for $x\cdot {\bf v_1}\leq 4\sigma$
and $\Phi(x)=0$ for $x\cdot {\bf v_1} \geq 8\sigma$.
Define the localization functions $\{\Phi_k\}$ by
\be\ba \label{phi-local}
&\Phi_1(x) :=\Phi(x-x_1), \ \ \Phi_K(x) :=1-\Phi(x-x_{K-1}),  \\
&\Phi_k(x) :=\Phi(x-x_{k})-\Phi(x-x_{k-1}),\ \ 2\leq k\leq K-1.
\ea\ee
One has the partition of unity $1= \sum_{j=1}^K \Phi_k$.
Then,
\begin{align}  \label{R-Rj}
   R =\sum_{k=1}^{K}R_{k}, \ \ with\ \ R_k:= R \Phi_k.
\end{align}
The corresponding renormalized remainders $\ve_{k}$, $1\leq k\leq K$, are defined by
\begin{align} \label{Rj-ej}
R_{k}(t,x) =\lbb_{k}^{-\frac d2} \varepsilon_{k} (t,\frac{x-\a_{k}}{\lbb_{k}}) e^{i\theta_{k}}.
\end{align}

The following  almost orthogonality between profiles $\{R_k\}$
and $\{U_k\}$ is a consequence of the orthogonality \eqref{ortho-cond-Rn-wn} and
the decoupling Lemma \ref{Lem-decoup-U}.

\begin{lemma} (Almost orthogonality \cite[Lemma 4.4]{SZ20})  \label{Lem-almost-orth}
Let $t^*$ be as in Theorem \ref{Thm-Mod-bdd}.
Then, for $t^*$ and $T_*$ close to $T$,
there exists $\delta >0$ such that for every
$1\leq k\leq K$ and any  $t\in [t^*,T_*]$,
\be\ba \label{Orth-almost}
&|{\rm Re}\int (x-\a_{k}) U_{k} \ol{R_{k}}dx|
  + |{\rm Re} \int |x-\a_{k}|^2 U_{k} \ol{R_{k}}dx|\leq Ce^{-\frac{\delta}{T- t}}\|R\|_{L^2},  \\
& |{\rm Im}\int \nabla U_{k} \ol{R_{k}}dx|
  + |{\rm Im}\int \Lambda_k U_{k} \ol{R_{k}}dx|
  + |{\rm Im}\int  \varrho_{k} \ol{R_{k}}dx|\leq Ce^{-\frac{\delta}{T- t}}\|R\|_{L^2}.
\ea\ee
\end{lemma}

Furthermore, by \eqref{equa-NLS-perturb}, \eqref{equa-z} and \eqref{u-dec},
the remainder $R$ satisfies the equation
\begin{align} \label{equa-R}
   i\partial_t R +\Delta R + a_1 \cdot \nabla R+ a_0 R +(f(v)-f(U+z))=-\eta,
\end{align}
where
\begin{align} \label{etan-Rn}
    \eta =& i\partial_t U + \Delta U + a_1 \cdot \nabla U+ a_0 U +f(U+z) - f(z),
\end{align}
and $f(z):= |z|^\frac 4d z$, $f(U+z)$ is defined similarly.

The estimates of $\eta$ are contained in Lemma \ref{Lem-eta} below.
\begin{lemma}  \label{Lem-eta}
Suppose that $P =\calo(T-t)$
and $|\a_k - x_k| <\frac 12$ for any $t\in [t^*,T_*]$.
Then,
\begin{align}
   & |\eta(t,x)| \leq C (T-t)^{-\frac d2 -2}
     \sum\limits_{k=1}^K  \(Mod +  |\ve_{z,k}(y)|+(T-t)^{\upsilon_*+1} \) e^{-\delta |y|} \bigg|_{y=\frac{x-\a_k}{\lbb_k}}
       + C \wt{\eta}, \label{eta-tx}
\end{align}
where $\wt \eta$ satisfies
$
     \|\wt \eta(t) \|_{L^2}
     \leq C e^{-\frac{\delta}{T-t}},
$
and   for any multi-index $\upsilon$ with $|\upsilon|\leq 2$,
\begin{align}
   & \|\partial_x^\upsilon \eta(t)\|_{L^2}
      \leq C(T-t)^{-2-|\upsilon|}\(Mod+\a^* (T-t)^{m+1+\frac d2}+(T-t)^{\upsilon_*+1} \).  \label{eta-L2}
\end{align}
\end{lemma}

{\bf Proof.}
Let $\Psi_k$ be as in \eqref{Psik-def}.
We decompose $\eta$ into four parts:
\begin{align} \label{eta1-eta4}
   \eta   =& \eta_1 + \eta_2 + \eta_3+\eta_4,
\end{align}
where
\begin{align}
     & \eta_1= \sum\limits_{k=1}^K \frac{ e^{i \theta_k(t)}}{\lbb_k(t)^{2+\frac d2}} \Psi_k(t, \frac{x-\a_k(t)}{\lbb_k(t)}), \label{eta1-def} \\
     & \eta_2 = f\(U + z\) - f\(U\) - f(z), \label{eta2-def}\\
     & \eta_3 = f(U) - \sum\limits_{k=1}^K f(U_k),    \label{eta3-def} \\
      & \eta_4=a_1 \cdot \nabla U+ a_0 U. \label{eta4-def}
\end{align}

By the exponential decay \eqref{Q-decay} of ground state
and $\|\ve_{z,k}\|_{L^\9} \leq C$,
\begin{align} \label{eta124-pt-esti}
   |\eta_1+\eta_2 |
   \leq C (T-t)^{-\frac d2 -2} \sum\limits_{k=1}^K \(Mod_k +  |\ve_{z,k}(t,y)|\)
   e^{- \delta |y|} \bigg|_{y=\frac{x-\a_k}{\lbb_k}},
\end{align}
and $\wt \eta := |\eta_3|$ contains different blow-up profiles,
and thus, by Lemma \ref{Lem-decoup-U},
\begin{align}
   \|\wt \eta(t) \|_{L^2} \leq C e^{-\frac{\delta}{T-t}}.
\end{align}
Moreover, since
\begin{align}
   \eta_4(t,x)
   = \sum\limits_{k=1}^K \lbb_k^{-\frac d2-1}(t) \wt a_{1,k}(t,y) \na Q_k(t,y) e^{i\theta_k}
     + \lbb_k^{-\frac d2} (t)  \wt a_{0,k} (t,y) Q_k(t,y)  e^{i\theta_k} \big|_{y=\frac{x-\a_k}{\lbb_k}},
\end{align}
where $\wt a_{1,k}, \wt a_{0,k}$ are as in Lemma \ref{Lem-wtbwtc},
using Lemma \ref{Lem-wtbwtc}, \eqref{Q-decay} and $P\leq C(T-t)$ we get
\begin{align} \label{eta3-pt-esti}
    |\eta_4(t,x)|
    \leq& C\sum\limits_{k=1}^K \( (T-t)^{-\frac d2-1}P^{\upsilon_*} \<y\>^{\upsilon_*+1} e^{-\delta|y|}
          +  (T-t)^{-\frac d2} P^{\upsilon_*-1} \<y\>^{2\upsilon_*}e^{-\delta|y|} \) \big|_{y=\frac{x-\a_k}{\lbb_k}}  \nonumber \\
    \leq& C \sum\limits_{k=1}^K (T-t)^{\upsilon_*-\frac d2-1}  e^{-\frac \delta 2|y|}  \big|_{y=\frac{x-\a_k}{\lbb_k}}.
\end{align}
Hence, \eqref{eta124-pt-esti} and \eqref{eta3-pt-esti} together yield \eqref{eta-tx}.

Concerning \eqref{eta-L2}, by \eqref{Psik-def}, it is clear that
\begin{align}
    \|\partial^\upsilon_x \eta_1\|_{L^2} \leq C (T-t)^{-2-|\upsilon|} Mod.
\end{align}
Moreover, expanding $f$ and then using the exponential decay \eqref{Q-decay} of ground state
we have
\begin{align}
    \|\partial^\upsilon_x \eta_2\|_{L^2}
     \leq C  (T-t)^{-2-|\upsilon|} \sum\limits_{|\upsilon|\leq 2} \|e^{-\delta |y|} \partial_y^\upsilon \ve_{z,k}\|_{L^\9}
          + C e^{-\frac{\delta}{T-t}}.
\end{align}
Since $\eta_3$ contains the interactions between different blow-up profiles,
by Lemma \ref{Lem-decoup-U},
\begin{align}
    \|\partial^\upsilon_x \eta_3\|_{L^2}  \leq C e^{-\frac{\delta}{T-t}}.
\end{align}
At last, applying Lemma \ref{Lem-wtbwtc} we also infer that
\begin{align}
  \|\partial^\upsilon_x \eta_4\|_{L^2}
  \leq \sum\limits_{k=1}^K \lbb_k^{-|\upsilon|+\frac d2}
        \|\partial_y^\upsilon(\wt a_1 \lbb_k^{-\frac d2-1} \na Q_k + \wt a_0 \lbb_k^{-\frac d2} Q_k) \|_{L^2}
  \leq  C (T-t)^{\upsilon_*-|\upsilon|-1}.
\end{align}

Therefore, putting the above estimates altogether we obtain \eqref{eta-L2}.
\hfill $\square$

\section{Localized mass and (generalized) energy} \label{Sec-Mass-Energy}

This section is devoted to the key estimates of localized mass, energy and the generalized energy.

\subsection{Control of localized mass}     \label{Subsec-Loc-Mass}

Recall that the localized mass is defined by
\begin{align} \label{Mj-def}
 M_k : =2{\rm Re} \< R_{k}, {U_{k}}\>+\int |R|^2\Phi_kdx,
\end{align}
where $R_k=R\Phi_k$
and $\{\Phi_k\}$ are the localization functions given by \eqref{phi-local}.

The main estimate
is contained in Theorem \ref{Thm-mass-local} below.

\begin{theorem} (Control of localized mass)    \label{Thm-mass-local}
Suppose that $P= \calo(T-t)$.
Then, there exists $C>0$
such that for every $1\leq k\leq K$,
\begin{align} \label{mass-local-esti}
   | M_k(t) |  \leq C \int_t^{T_*} \( \a^*D + \frac{D^2}{T-s} \) ds
        + C\(\a^*D +  \a^* (T-t)^{m+1+\frac d2} \),\ \ t\in [t^*, T_*].
\end{align}
\end{theorem}

{\bf Proof.}
On one hand, the geometrical decomposition \eqref{u-dec}  yields the expansion:
\begin{align} \label{u2Phij-expan}
\int |v (t)|^2\Phi_kdx
&= \int (|U|^2+|z|^2+|R|^2)\Phi_kdx  +2{\rm Re}\int ( R\overline{U}+z \ol{R}+z \overline{U})\Phi_k dx  \nonumber \\
&= \int |U|\Phi_k^2dx+\int |z |^2\Phi_kdx +\int |R |^2\Phi_kdx \nonumber  \\
&\quad+2{\rm Re}\int R_k  \overline{U_k}  dx+2{\rm Re}\int z \overline{U_k} dx +2{\rm Re}\int z \overline{R}_k  dx
      + \calo(e^{-\frac{\delta}{T-t}}\|R\|_{L^2}),
\end{align}
where the last step is due to Lemma \ref{Lem-decoup-U}
and $\delta>0$.

On the other hand,
since $v(T_*) = S(T_*) + z(T_*)$,
we have
\begin{align} \label{uT2Phij-expan}
\int |v(T_*)|^2\Phi_kdx
= \int (|S(T_*)|^2+|z(T_*)|^2)\Phi_kdx  +2{\rm Re}\int (z\overline{S} )(T_*)\Phi_k dx  .
\end{align}

Note that,
the integrations $\int |z|^2 \Phi_k dx$ and  $\int |z(T_*)|^2 \Phi_k dx$
only contribute a small constant $(\a^*)^2$,
which, however, is  insufficient to close the bootstrap arguments later.
The key point is that one more factor $D$ can be explored by
subtracting  \eqref{uT2Phij-expan} from \eqref{u2Phij-expan}
and then using both the dynamics generated by equations \eqref{equa-NLS-perturb} and \eqref{equa-z}.

To be precise, we derive from \eqref{u2Phij-expan} and \eqref{uT2Phij-expan} that
\begin{align} \label{uL2-phi}
   | M_k (t) |
  \leq& \bigg|\int (|v(t)|^2-|v(T_*)|^2)\Phi_k dx -\int  (|z(t)|^2- |z(T_*)|^2)\Phi_kdx\bigg|  \nonumber \\
   & +  \bigg|\int (|U(t)|^2-|S(T_*)|^2)\Phi_kdx \bigg|   \nonumber \\
  &  +2\(\bigg|\int (z \overline{U}_k )(t) dx\bigg|+ \bigg|\int (z \overline{R}_k )(t) dx\bigg|
     + \bigg|\int (z \overline{S} )(T_*)\Phi_kdx\bigg|\)
     + C e^{-\frac{\delta}{T-t}} \|R\|_{L^2}   \nonumber \\
   =&: K_1 + K_2 +K_3 + C e^{-\frac{\delta}{T-t}} \|R\|_{L^2} .
\end{align}

Let us first treat the easier two terms $K_2$, $K_3$.
Actually, it holds that (see \cite[(5.22),(5.23)]{SZ20})
\begin{align}
  & \int |U (t)|^2\Phi_kdx
   =\|Q\|_{L^2}^2 + \calo(e^{-\frac{\delta}{T-t}}),  \label{U-Q-esti} \\
 & \int |S(T_*)|^2\Phi_k dx
  =\|Q\|_2^2+ \calo(e^{-\frac{\delta}{T-T_*}})
  =\|Q\|_2^2+ \calo(e^{-\frac{\delta}{T-t}}),   \label{utn-Q-esti}
\end{align}
which yields that
\begin{align} \label{uT-Ut}
   K_2 (t) \leq C e^{-\frac{\delta}{T-t}}.
\end{align}
Moreover, by \eqref{z-LtHm-a*} and \eqref{eyvez-bdd},
\begin{align}\label{r-rz}
K_3(t)
 \leq& C \(\|z\|_{L^2} \|R(t)\|_{L^2} +  \|e^{-\delta|y|} (|\ve_{z,k}(t)|+|\ve_{z,k}(T_*)|)\|_{L^\9} + e^{-\frac{\delta}{T-T_*}} \) \nonumber \\
 \leq& C \(\al^* D + \a^* (T-t)^{m+1+\frac d2} \)
\end{align}

Hence, it remains to treat the first term $K_1$ on the R.H.S. of \eqref{uL2-phi}.

For this purpose,
we derive from equation \eqref{equa-NLS-perturb} that
\begin{align} \label{du2-bc}
 \frac{d}{dt}\int |v|^2\Phi_kdx
 =& {\rm Im}\int (2\ol{v}\nabla v+ a_1 |v|^2)\cdot  \nabla\Phi_k dx.
\end{align}
Similarly, by equation \eqref{equa-z},
\begin{align}  \label{dz2-bc}
 \frac{d}{dt}\int |z|^2\Phi_kdx
 =& {\rm Im}\int (2\ol{z}\nabla z+ a_1 |z|^2)\cdot  \nabla\Phi_k dx.
\end{align}
Thus,
\begin{align} \label{du2-dz2}
  \frac{d}{dt}\int |v|^2\Phi_kdx -\frac{d}{dt}\int |z|^2\Phi_kdx
 =  {\rm Im} \int 2(\ol{v}\nabla v-\ol{z}\nabla z) \cdot \na \Phi_k + (|v|^2-|z|^2) a_1 \cdot \na \Phi_k dx.
\end{align}

Note that, by \eqref{u-dec},
\begin{align} \label{unau-znaz}
   \ol{v}  \na v  -  \ol{z}  \na z
   = \ol{U}  \na (U+R+z) + (\ol R+\ol z) \na {U}
      + \ol{z} \na R + \ol{R}   \na z + \ol{R} \na R .
\end{align}
Since $P=\calo(T-t)$, $|x_k - \a_k(t)|\leq \sigma$, $1\leq k\leq K$,
${\rm supp} \na \Phi_k \subseteq \cap_{k=1}^K \{x: |x-\a_k|\geq 3 \sigma\}$.
By  \eqref{Q-decay},
\begin{align} \label{UnaURz}
   \bigg| \int \( \ol{U}\cdot \na (U+R+z) + (\ol R+\ol z) \na {U} \) \cdot \na \Phi_k dx \bigg|
  \leq C e^{-\frac{\delta}{T-t}}.
\end{align}
Moreover,   the integration by parts formula yields
\begin{align}  \label{znaR}
   \bigg| \int \ol{z} \na R \cdot \na \Phi_k dx \bigg|
   =  \bigg| \int R  \na \ol{z} \cdot \na \Phi_k + R \ol{z} \Delta \Phi_k  dx \bigg|
  \leq C \|z\|_{H^1} \|R\|_{L^2}.
\end{align}
Thus, it follows from \eqref{unau-znaz}-\eqref{znaR} that
\begin{align}  \label{du2-dz2.1}
   \bigg| {\rm Im} \int (\ol{v} \na v - \ol{z} \na z ) \cdot \na \Phi_k dx \bigg|
   \leq C \(  \|z\|_{H^1} \|R\|_{L^2} +  \|R\|_{L^2} \|\na R\|_{L^2}  + e^{-\frac{\delta}{T-t}} \).
\end{align}

Similarly, we have
\begin{align}   \label{du2-dz2.2}
    \bigg| {\rm Im} \int (|v|^2 - |z|^2) a_1 \cdot \na \Phi_k dx  \bigg|
    \leq& C \int |\na \Phi_k|\  \big| |U|^2 + |R|^2 + 2 {\rm Re} (R \ol{U} + z\ol{U}  +z \ol{R}) \big| dx \nonumber \\
    \leq& C \( \|z\|_{L^2} \|R\|_{L^2} + \|R\|_{L^2}^2 + e^{-\frac{\delta}{T-t}}\).
\end{align}

Hence,
we conclude from \eqref{z-LtHm-a*}, \eqref{du2-dz2}, \eqref{du2-dz2.1} and \eqref{du2-dz2.2} that
\begin{align*}
      \bigg|\frac{d}{dt}\int |v|^2\Phi_kdx -\frac{d}{dt}\int |z|^2\Phi_kdx \bigg|
 \leq&  C \( \|z\|_{H^1} \|R\|_{L^2} + \|R\|_{L^2}^2 + \|R\|_{L^2}\|\na R\|_{L^2}  + e^{-\frac{\delta}{T-t}} \) \nonumber \\
 \leq& C \(\a^* D + \frac{D^2}{T-t} + e^{-\frac{\delta}{T-t}} \).
\end{align*}
where $\delta >0$.
Integrating both sides
we then obtain
\begin{align}  \label{u-utn}
  K_1 \leq C\int_{t}^{T_*} \a^* D  + \frac{D^2}{T-s} \ ds
+C e^{-\frac{\delta}{T-t}}.
\end{align}

Therefore,
plugging \eqref{uT-Ut}, \eqref{r-rz} and \eqref{u-utn} into \eqref{uL2-phi}
we obtain \eqref{mass-local-esti}.
The proof is complete.
\hfill $\square$

\subsection{Refined estimate of $\beta$} \label{Subsec-Energy}

In this subsection we shall derive the refined estimate of
parameter $\beta=(\beta_k)$
from the energy functional, defined by
\begin{align}   \label{energy-conserv}
  E(v)  := \half\int_{\R^d}|\nabla v|^2dx-\frac{d}{2d+4}\int_{\R^d}|v|^{2+\frac{4}{d}}dx.
\end{align}

Unlike in the NLS case,
the energy of solutions to \eqref{equa-NLS-perturb} is no longer conserved,
it is thus important to first control the variation of energy.
This is the content of Lemma \ref{Lem-energy} below.

\begin{lemma} (Variation of the energy) \label{Lem-energy}
Suppose $P = \calo(T-t)$.
Then, there exists $C>0$
such that
\begin{align} \label{esti-Eut-Eutn}
     \bigg| \frac{d}{dt} E(v) - \frac{d}{dt} E(z) \bigg|
     \leq C \( \a^*D + \frac{D^2}{(T-t)^2}  + (T-t)^{\upsilon_*-3} \),\ \ \forall t\in[t^*,T_*].
\end{align}
\end{lemma}

\begin{remark}
As in the proof of Theorem \ref{Thm-mass-local},
the key point here is that
one more factor $D$ can be gained from the difference between the energies of $v$ and $z$.
\end{remark}

{\bf Proof.}
As  in the previous case of localized mass,
we consider the difference between
two energies of $u$ and $z$.
Using \eqref{equa-NLS-perturb} and \eqref{equa-z} we compute,
as in \cite[(5.26)]{SZ20},
\begin{align} \label{dtE}
\frac{d}{dt}E(v)
=&-2\sum\limits_{l=1}^N h_l {\rm Re}\int \nabla^2 \phi_l(\nabla v,\nabla \ol{v})dx
+\frac{1}{2}\sum\limits_{l=1}^N  h_l \int \Delta^2 \phi_l|v|^2dx  \nonumber \\
& +\frac{2}{d+2}\sum\limits_{l=1}^N   h_l \int\Delta \phi_l|v|^{2+\frac{4}{d}}dx
 -\sum\limits_{j=1}^d {\rm Im}\int \na\(\sum\limits_{l=1}^N \partial_j \phi_l  h_l\)^2 \cdot \na v \ol{v} dx,
\end{align}
and
\begin{align} \label{dtz}
\frac{d}{dt}E(z)
=&-2\sum\limits_{l=1}^N  h_l {\rm Re}\int \nabla^2 \phi_l(\nabla z,\nabla \ol{z})dx
+\frac{1}{2}\sum\limits_{l=1}^N  h_l \int \Delta^2 \phi_l|z|^2dx  \nonumber \\
& +\frac{2}{d+2}\sum\limits_{l=1}^N   h_l \int\Delta \phi_l|z|^{2+\frac{4}{d}}dx
 -\sum\limits_{j=1}^d {\rm Im}\int \na\(\sum\limits_{l=1}^N \partial_j \phi_l  h_l\)^2 \cdot \na z \ol{z} dx.
\end{align}

In order to control the difference $\frac{d}{dt}E(v) - \frac{d}{dt}E(z)$,
we first see that,
by \eqref{u-dec}, integration by parts formula and \eqref{eyvez-bdd},
\begin{align} \label{na2phi-nau-naz}
   & \bigg| {\rm Re} \int \na^2 \phi_l(\na v, \na \ol{v}) dx -   {\rm Re} \int \na^2 \phi_l(\na z, \na \ol{z}) dx \bigg|   \nonumber \\
  =& \bigg|  \sum\limits_{i,j=1}^d {\rm Re}
  \int \partial_{ij} \phi_l \partial_i (U+ R) \partial_j (\ol{U+R})
       -  \partial_{ijj} \phi_l \partial_j z  (\ol{U+R})
       - \partial_{iij} \phi_l ({U+R})  \partial_j \ol{z} \nonumber \\
   & \qquad  \qquad - \partial_{ij} \phi_l  \partial_{ij}  z  (\ol{U+R})
     - \partial_{ij} \phi_l (U+R) \partial_{ij} \ol{z} dx \bigg|   \nonumber \\
  \leq& C \sum\limits_{i,j=1}^d  \bigg( \|R\|_{H^1}^2
          + \int (|\partial_j z| + |\partial_{ij} z|)  |R|dx
          + \bigg| \int \partial_{ij} \phi_l (\partial_i U \partial_j \ol{U} + \partial_i U \partial_j \ol{R} + \partial_i R \partial_j \ol{U}) dx \bigg| \nonumber \\
      &  \qquad + \bigg| \int z \partial_i (\partial_{ijj} \phi_l \ol{U}) + \ol{z} \partial_j (\partial_{iij} \phi_l {U}) dx  \bigg|
        + \bigg|\int z\partial_{ij} (\partial_{ij} \phi_l \ol{U}) dx \bigg|   \bigg) \nonumber \\
  \leq& C \(\|z\|_{H^2} \|R\|_{L^2}  + \|R\|_{H^1}^2
            + (T-t)^{\upsilon_*-3} (1+ \sum_{k=1}^K \|\na \ve_k\|_{L^2})
            + (T-t)^{\upsilon_*-3}\sum_{k=1}^K \| e^{-\delta |y|}   \ve_{z,k}\|_{L^\9}
              + e^{-\frac{\delta}{T-t}}\).
\end{align}

Similarly,
\begin{align}  \label{D2phi-u2-z2}
    \bigg| \int \Delta^2 \phi_l (|v|^2 -|z|^2) dx  \bigg|
    \leq& C\( \|R\|_{L^2}^2 + \int |z| |R| dx + \int |\Delta^2 \phi_l| (|U|^2 + |R U| + |z U|) \) dx \nonumber \\
    \leq& C \bigg(\|R\|_{L^2}^2 + \|z\|_{L^2}  \|R\|_{L^2} + (T-t)^{\upsilon_*-3}
          + (T-t)^{\upsilon_* -3} \sum\limits_{k=1}^K  \| e^{-\delta |y|} \ve_{z,k}\|_{L^\9}  \bigg),
\end{align}
and
\begin{align} \label{Dphifv-fU}
   \bigg| \int \Delta \phi_l (|v|^{2+\frac 4d} - |z|^{2+\frac 4d}) dx \bigg|
    \leq& C \int |\Delta \phi_l| \(|U+R|^{2+\frac 4d}+ |z|^{1+\frac 4d} |U+R| \) dx \nonumber \\
    \leq& C\( \|R\|_{H^1}^{2+\frac 4d}  +\|z\|_{H^1}^{1+\frac 4d}  \|R\|_{L^2}
          + (T-t)^{\upsilon_*-3}
          +  (T-t)^{\upsilon_*-3} \sum\limits_{k=1}^K  \| e^{-\delta |y|}\ve_{z,k} \|^{1+\frac 4d}_{L^\9} \).
\end{align}

Moreover, by \eqref{unau-znaz},
\begin{align} \label{nauu-nazz}
    & {\rm Im} \int \na \(\sum\limits_{j=1}^d \partial_j \phi_l\)^2 \cdot (\na v\ol{v} - \na z \ol{z}) dx  \nonumber \\
  =& {\rm Im} \int \na \(\sum\limits_{j=1}^d  \partial_j \phi_l\)^2 \cdot
   \bigg( \na U (\ol{U} + \ol{R} + \ol{z})
           + (\na R + \na z) \ol U
         + \na R \ol{R} +  \na R \ol{z} + \na z \ol{R} \bigg) dx.
\end{align}
Note that,
by the integration by parts formula,
\eqref{eyvez-bdd} and \eqref{phik-Taylor},
\begin{align}
     &  \bigg|  \int \na \(\sum\limits_{j=1}^d \partial_j \phi_l\)^2 (\na z \ol{U} + \na R \ol{z}) dx \bigg|    \nonumber \\
   =& \bigg| \int \Delta  \(\sum_{j=1}^d \partial_j \phi_l\)^2 (z \ol{U} + R \ol{z}) dx
       + \int \na \(\sum_{j=1}^d \partial_j \phi_l\)^2 \cdot (z \na \ol{U} + R \na \ol{z}) dx \bigg|  \nonumber \\
   \leq& C\( \|z\|_{H^1} \|R\|_{L^2}
          + (T-t)^{2\upsilon_*-2} \sum_{k=1}^K \| e^{-\delta |y|} \ve_{z,k}\|_{L^\9}\).
\end{align}
Using \eqref{phik-Taylor} again we also have
\begin{align}
  & \bigg| {\rm Im} \int \na \(\sum\limits_{j=1}^d \partial_j \phi_l\)^2
         \cdot \( \na U(\ol{U} + \ol{R} + \ol{z}) + \na R \ol{U} \) dx \bigg| \nonumber \\
  \leq& C \( (T-t)^{2\upsilon_*-2}
          \(1+ \sum\limits_{k=1}^K \|e^{-\delta|y|} \ve_{z,k}\|_{L^\9} + D\)
          + e^{-\frac{\delta}{T-t}} \) \nonumber \\
  \leq& C (T-t)^{2\upsilon_*-2}.
\end{align}
Plugging these into \eqref{nauu-nazz} we get
\begin{align} \label{napphij2}
     & \bigg| {\rm Im} \int \na \(\sum \partial_j \phi_h\)^2 \cdot (\na v \ol{v} - \na z \ol{z}) dx \bigg|  \nonumber \\
  \leq&  C\( \|z\|_{H^1} \|R\|_{L^2} + \|R\|_{L^2} \|\na R\|_{L^2}
           +  (T-t)^{2\upsilon_*-2}
        +  (T-t)^{2\upsilon_*-2} \sum\limits_{k=1}^K \| e^{-\delta |y|} \ve_{z,k}\|_{L^\9} \).
\end{align}

Therefore,
we conclude from the estimates \eqref{na2phi-nau-naz},
\eqref{D2phi-u2-z2}, \eqref{Dphifv-fU} and \eqref{napphij2} that
\begin{align*}
  \bigg|\frac{d}{dt}E(v)-\frac{d}{dt}E(z)\bigg|
 \leq & C\(   \(\|z\|_{H^2} + \|z\|_{H^1}^{1+\frac 4d}\) \|R\|_{L^2}
             + \|R\|_{H^1}^2  + \|R\|_{H^1}^{2+\frac 4d}
        + (T-t)^{\upsilon_*-3}   \) \nonumber \\
 \leq& C\(\a^* D + \frac{D^2}{(T-t)^2} + (T-t)^{\upsilon_*-3} \),
\end{align*}
which yields \eqref{esti-Eut-Eutn},
thereby finishing the proof.
\hfill $\square$

We are now in position to derive the refined estimate of $\beta=(\beta_k)$,
which is essentially a consequence of
the coercivity of energy around the ground state.

\begin{theorem} (Improved estimate of $\beta$) \label{Thm-Energy}
Suppose that
$P =\calo (T-t)$ and $D=o(1)$.
Then, for any  $t\in[t^*,T_*]$,
\begin{align}  \label{energy-esti}
     \sum_{k=1}^K  \frac{|\beta_{k}|^2}{2\lambda_{k}^2} \|Q\|^2_{L^2}
\leq  \frac{\|yQ\|_{L^2}^2}{8} \sum_{k=1}^K  \(w_k^2 - \frac{\g_k^2}{\lbb_k^2}\)
       + \calo(Er),
\end{align}
where the error term
\begin{align} \label{Er}
    Er :=  \int_t^T \(\a^*D + \frac{D^2}{(T-s)^2}\) ds
           + \a^* D + \sum\limits_{k=1}^K \frac{|M_k|}{(T-t)^2}
              + (T-t)^{\upsilon_*-2} + \a^*(T-t)^{m-1+\frac d2}.
\end{align}
\end{theorem}

{\bf Proof.}
Let $F(v):= \frac{d}{2d+4}|v|^{2+\frac 4d}$,
$F(U+z)$ and $F(z)$ are defined similarly.
Set $f(v):= |v|^{\frac 4d} v$.
Rewrite
\begin{align}
   E(v) = E(v) + \sum\limits_{k=1}^K
          \frac{1}{\lbb_k^2} {\rm Re} \int \ol{U_k} R_k
          + \frac 12 |R|^2 \Phi_k dx
          - \sum\limits_{k=1}^K \frac{1}{2\lbb_k^2} M_k.
\end{align}
Using \eqref{u-dec} and \eqref{g-gzz-expan} we expand
\begin{align*}
   F(v) = F(U+z) + F'(U+z)\cdot R + F''(U+z,R)\cdot R^2.
\end{align*}
Applying \eqref{g-gzz-expan} again to $F(U+z)$ we get
\begin{align}
  F(v) = F(U) + F'(U)\cdot z + F''(U,z)\cdot z^2 + F'(U+z)\cdot R + F''(U+z,R)\cdot R^2.
\end{align}
Then, taking into account
$F'(U)\cdot z = {\rm Re} (f(U) \ol{z})$,
$F'(U+z)\cdot z = {\rm Re} (f(U+z) \ol{z})$
and the expansion
\begin{align*}
   \frac 12 \|\na v\|_{L^2}^2
   =  \frac 12 \|\na U\|_{L^2}^2  +  \frac 12 \|\na z\|_{L^2}^2  +  \frac 12 \|\na R\|_{L^2}^2
     - {\rm Re} \<\Delta U+ \Delta z, R\>
     - {\rm Re} \<\Delta U, z\>,
\end{align*}
we obtain
\begin{align} \label{Eu-expan}
E(v) =& E(U)+E(z) -\sum\limits_{k=1}^K \frac{1}{2\lbb_k^2} M_k
   -  \( {\rm Re}\int (\Delta U+\Delta z+|U+z|^{\frac 4d}(U+z)) \overline{R} dx
          - \sum\limits_{k=1}^K \frac{1}{\lbb_k^2} {\rm Re} \int \ol{U_k} R_k dx \)
    \nonumber \\
    &+ \( \int \frac{1}{2}|\nabla R|^2dx
           + \sum\limits_{k=1}^K \frac{1}{2\lbb_k^2} \int |R|^2 \Phi_k dx
           -{\rm Re}\int  F''(U+z, R)\cdot R^2 dx \) \nonumber \\
   &-  {\rm Re}\int ({\Delta U+|U|^{\frac 4d}U}) \overline{z} dx
      -{\rm Re}\int (F''(U, z)\cdot z^2 - F(z))  dx \nonumber \\
   =&: E(U) + E(z) -\sum\limits_{k=1}^K \frac{1}{2\lbb_k^2} M_k  + \sum\limits_{l=1}^4 E_l.
\end{align}
Note that,
$E_1$ and $E_2$ are ordered by the homogeneity of $R$,
and $E_3$ and $E_4$ contain the perturbations with the regular profile $z$.

Next we estimate $E_l$ separately, $1\leq l\leq 4$.

For  the linear term $E_1$,
by Lemma \ref{Lem-decoup-U},
\begin{align}
E_1=& - \sum_{k=1}^{K} {\rm Re} \int ({\Delta U_{k}- \lbb_{k}^{-2} U_{k}+|U_{k}|^{\frac 4d}U_{k}}) \overline{R_{k}} dx  \nonumber\\
& - {\rm Re}\int (\Delta z+|U+z|^{\frac 4d}(U+z)-|U|^{\frac 4d}U) \overline{R} dx
   + \calo(e^{-\frac{\delta}{T-t}}\|R\|_{L^2}) \nonumber\\
=&:E_{11}+E_{12}+ \calo(e^{-\frac{\delta}{T-t}}\|R\|_{L^2}).
\end{align}
Using the identity \eqref{equa-Qk}
and the almost orthogonality \eqref{Orth-almost}
we have (see \cite[(3.38)]{CSZ21}),
\begin{align} \label{E11}
E_{11}= & - \sum\limits_{k=1}^K  \frac{1}{\lambda_{k}^2} {\rm Im} \int(\gamma_{k}\Lambda Q_{k} -2\beta_{k}\cdot\nabla Q_{k}) \ol{\varepsilon_k} dx
          - \sum\limits_{k=1}^K \frac{1}{\lbb_{k}^2} {\rm Re} \int |\beta_{k} - \frac {\gamma_{k}}{ 2} y|^2 Q_{k} \ol{\varepsilon_k} dx \nonumber \\
      =&  - \sum_{k=1}^{K}\frac{|\beta_k|^2}{\lbb_k^2}{\rm Re} \int  U_{k} \ol{R}_k dx
           + \calo(e^{-\frac{\delta}{T-t}}\|R\|_{L^2}) \nonumber \\
      =&  - \sum_{k=1}^{K}\frac{|\beta_k|^2}{2\lbb_k^2} M_k + \calo\(\|R\|_{L^2}^2 + e^{-\frac{\delta}{T-t}}\|R\|_{L^2}\).
\end{align}
Moreover, by   \eqref{z-LtHm-a*} and \eqref{eyvez-bdd},
\begin{align}  \label{E12}
|E_{12}|
&\leq C \int (|\Delta z|+|U|^{\frac{4}{d}}|z|+|z|^{1+\frac{4}{d}})|R|dx  \nonumber \\
&\leq C \(\|\Delta z\|_{L^2}+  \sum\limits_{k=1}^K \lbb_k^{-2}\|e^{-\delta |y|} \ve_{z,k}\|_{L^\9}
 +\|z\|^{1+\frac4d}_{H^1} + e^{-\frac{\delta}{T-t}}\)\|R\|_{L^2} \nonumber \\
 &\leq C \( \a^*\|R\|_{L^2} +   e^{-\frac{\delta}{T-t}}\).
\end{align}
Thus, we obtain
\begin{align}\label{E1}
  E_1 = - \sum\limits_{k=1}^K \frac{|\beta_k|^2}{2\lbb_k^2} M_k
        + \calo(\a^* \|R\|_{L^2}   + e^{-\frac{\delta}{T-t}})
        + o\(\frac{D^2}{(T-t)^2}\).
\end{align}

Concerning the second term $E_2$,
set
\begin{align}
   \wt E_2 := - {\rm Re} \int F''(U+z, R)\cdot R^2 - F''(U, R)\cdot R^2 dx.
\end{align}
We estimate
\begin{align} \label{E2-expan}
 E_2 =& \frac 12 \int |\na R|^2 dx
         + \sum\limits_{k=1}^K \frac{1}{2\lbb_k^2} \int |R|^2 \Phi_k dx
         - {\rm Re} \int F''(U,R)\cdot R^2 dx
         + \wt E_2 \nonumber \\
      =& \frac 12 \int |\na R|^2 dx
         + \sum\limits_{k=1}^K \frac{1}{\lbb_k^2}  |R|^2 \Phi_k
         - (1+\frac 2d) |U|^\frac 4d |R|^2
         - \frac 2d |U|^{\frac 4d-2} U^2 \ol{R}^2 dx
         + \wt E_2
         + \calo\(\frac{D^3}{(T-t)^2}\) \nonumber \\
      \geq& \wt C \frac{D^2}{(T-t)^2}  + \wt E_2
            + \calo\(\sum\limits_{k=1}^K \frac{M_k^2}{(T-t)^2}   + e^{-\frac{\delta}{T-t}}\),
\end{align}
where $\wt C>0$, the error in the second step is caused by the remainders of orders higher than two
(see \cite[(3.34)]{CSZ21}),
and the last step is mainly due to the local coercivity of linearized operator in
Lemma \ref{Lem-coer-f-local},
see the proof of \cite[(3.39)]{CSZ21},
and $\frac{D^3}{(T-t)^2} = o\(\frac{D^2}{(T-t)^2}\)$.
The error $\wt E_2$ can be bounded by
\begin{align}
   |\wt E_2|
   \leq& C \int (|U|^{\frac 4d-1} + |R|^{\frac 4d-1} + |z|^{\frac 4d-1}) |z||R|^2 dx \nonumber \\
   \leq& C \((T-t)^{-2} \sum\limits_{k=1}^K \|e^{-\delta |y|}\ve_{z,k}\|_{L^\9} \|R\|_{L^2}^2
           + \|z\|_{L^\9} \|R\|_{L^{1+\frac 4d}}^{1+\frac4d}
           + \|z\|_{L^\9}^{\frac 4d} \|R\|_{L^2}^2 + e^{-\frac{\delta}{T-t}}\) \nonumber \\
   \leq& C \((T-t)^{-2} \sum\limits_{k=1}^K \|e^{-\delta |y|}\ve_{z,k}\|_{L^\9}+ \a^* (T-t)^{-2+\frac d2} + \a^*\) D^2
         + C e^{-\frac{\delta}{T-t}} ,
\end{align}
where we also used \eqref{R-Lp-D}
and $D=\calo(1)$ in the last step.

Thus, for $t$ close to $T$ such that
$C \(  \sum_{k=1}^K \|e^{-\delta |y|}\ve_{z,k}\|_{L^\9} + 2(T-t)^{\frac d2} \) \leq \frac {1}{2} \wt C$,
it follows that
\begin{align} \label{E2}
   E_2  \geq& \frac {\wt C}{2} \frac{D^2}{(T-t)^2}
            + \calo\(\sum\limits_{k=1}^K \frac{M_k^2}{(T-t)^2}  + e^{-\frac{\delta}{T-t}}\),
\end{align}

The last two terms $E_3$ and $E_4$ can be estimated easily by using
\eqref{Q-decay} and \eqref{eyvez-bdd}:
\begin{align}
  & |E_3|   \leq C (T-t)^{-2} \sum\limits_{k=1}^K \| e^{-\delta|y|}\ve_{z,k}\|_{L^\9} + C e^{-\frac{\delta}{T-t}}
            \leq C \a^*(T-t)^{m-1+\frac d2},  \label{E3} \\
  & |E_4| \leq C \sum\limits_{j=2}^{1+\frac 4d} \int |U|^{2+\frac 4d-j}|z|^j dx
\leq C  (T-t)^{-2} \sum\limits_{k=1}^K  \| e^{-\delta|y|}\ve_{z,k}\|_{L^\9}
     + C e^{-\frac{\delta}{T-t}}
\leq C \a^*(T-t)^{m-1+\frac d2}. \label{E4}
\end{align}

Thus,
combining \eqref{Eu-expan}, \eqref{E1}, \eqref{E2}, \eqref{E3} and \eqref{E4}
we conclude that for some $C>0$,
\begin{align}  \label{enr-exp1}
E(v) \geq & E(U)+E(z) + C \frac{D^2}{(T-t)^2}
          - \sum\limits_{k=1}^K \frac{1+|\beta_k|^2}{2\lbb_k^2} M_k \nonumber \\
          &+ \calo\(\sum\limits_{k=1}^K \frac{M_k^2}{(T-t)^2} + \a^*D
           +   \a^*(T-t)^{m-1+\frac d2} \).
\end{align}

Furthermore,
since $v(T_*) = S(T_*)+z(T_*)$ we derive that
\begin{align} \label{EuT-expan}
E(v(T_*))=& E(S(T_*))+E(z(T_*))+{\rm Re}\int \nabla S(T_*)\nabla\bar{z}(T_*)dx \nonumber \\
          & -\int F(v(T_*))-F(S(T_*))-F(z(T_*))dx.
\end{align}
As in \eqref{E3} and \eqref{E4},
by \eqref{eyvez-bdd} and $T-T_* \leq T-t$,
\begin{align} \label{STnazT-esti}
  \bigg|{\rm Re}\int \nabla S(T_*)\nabla\bar{z}(T_*)dx \bigg|
  \leq& C(T-t)^{-2} \sum\limits_{k=1}^K  \| e^{-\delta|y|} \ve_{z,k}(T_*)\|_{L^\9}
        + C e^{-\frac{\delta}{T-T_*}}  \nonumber \\
  \leq& C \a^* (T-t)^{m-1+\frac d2},
\end{align}
and
\begin{align} \label{Fu-Fz-esti}
 \bigg|\int F(v(T_*))-F(S(T_*))-F(z(T_*))dx\bigg|
  \leq& C (T-t)^{-2} \sum\limits_{k=1}^K  \| e^{-\delta|y|} \ve_{z,k}(T_*)\|_{L^\9}  + Ce^{-\frac{\delta}{T-T_*}}  \nonumber \\
  \leq& C \a^* (T-t)^{m-1+\frac d2}.
\end{align}
Thus, it follows from \eqref{EuT-expan}, \eqref{STnazT-esti} and \eqref{Fu-Fz-esti} that
\begin{align}\label{enr-exp2}
E(v(T_*))=E(S(T_*))+E(z(T_*)) +  \calo\( \a^*(T-t)^{m-1+\frac d2} \).
\end{align}

Now, plugging \eqref{enr-exp2} into \eqref{enr-exp1} we derive
\begin{align}
E(U(t)) + C \frac{D^2}{(T-t)^2}
  \leq& \(E(v(t)) -E(v(T_*))\)- \(E(z(t))-E(z(T_*))\)+E(S(T_*))\nonumber\\
      &+ \sum\limits_{k=1}^K \frac{1+|\beta_k|^2}{2\lbb_k^2} M_k
          + \calo\(\sum\limits_{k=1}^K \frac{M_k^2}{(T-t)^2}
           + \a^*D +  \a^*(T-t)^{m-1+\frac d2} \).
\end{align}

Thus, by the variation control \eqref{esti-Eut-Eutn},
\begin{align} \label{EUt-EST}
E(U(t)) +  C \frac{D^2}{(T-t)^2}
    \leq&  E(S(T_*))
            + \calo(Er),
\end{align}
where $Er$ is as in \eqref{Er}.
Moreover,  \eqref{Uj-Qj-Q}
and Lemma \ref{Lem-decoup-U} yield
\begin{align}
  & E(U(t))
 = \sum_{k=1}^K\( \frac{|\beta_{k}|^2}{2\lambda_{k}^2} \|Q\|^2_{L^2}
+ \frac{\gamma_{k}^2}{8\lambda_{k}^2} \|yQ\|_{L^2}^2\)
+  \calo(e^{-\frac{\delta}{T-t}}) ,  \label{EUt} \\
  & E(S(T_*)) =\sum_{k=1}^K\frac{w_k^2}{8}\|yQ\|_{L^2}^2
  + \calo(e^{-\frac{\delta}{T-T_*}}), \label{EST}
\end{align}

Therefore, plugging \eqref{EUt} and \eqref{EST} into \eqref{EUt-EST}
we obtain \eqref{energy-esti}.
The proof is complete.
\hfill $\square$

\subsection{Monotonicity of generalized energy}  \label{Subsec-General-Energy}
This subsection is mainly devoted to the crucial
monotonicity property of generalized energy.

Let $\chi(x)=\psi(|x|)$ be a smooth radial function on $\R^d$,
where $\psi$ satisfies
$\psi'(r) = r$ if $r\leq 1$,
$\psi'(r) = 2- e^{-r}$ if $r\geq2$,
and
\be\label{chi}
 \bigg|\frac{\psi^{'''}(r)}{\psi^{''}(r)} \bigg|\leq C,
\ \ \frac{\psi'(r)}{r}-\psi^{''}(r) \geq0.
\ee
Let $\chi_A(x) :=A^2\chi(\frac{x}{A})$, $A\geq 1$,
$f(v):= |v|^{\frac 4d} v$
and $F(v):= \frac{d}{2d+4} |v|^{2+\frac 4d}$.

The generalized energy, adapted to the multi-bubble case, is defined by
\begin{align} \label{I-def}
   \mathscr{I}(t)
   := &\frac{1}{2}\int |\nabla R|^2dx+\frac{1}{2}\sum_{k=1}^K\int\frac{1}{\lambda_{k}^2} |R|^2 \Phi_kdx
           -{\rm Re}\int F(v)-F(U+z)-f(U+z)\ol{R}dx \nonumber \\
&+\sum_{k=1}^K\frac{\gamma_{k}}{2\lambda_{k}}{\rm Im} \int (\nabla\chi_A) \(\frac{x-\alpha_{k}}{\lambda_{k}}\)\cdot\nabla R \ol{R}\Phi_kdx
   =: \mathscr{I}^{(1)}+ \mathscr{I}^{(2)}.
\end{align}
where $\mathscr{I}^{(1)}$ mainly contains the quadratic terms of remainder (up to acceptable errors)
and  $\mathscr{I}^{(2)}$ is a Morawetz type functional.
The key monotonicity property
is formulated below.

\begin{theorem} (Monotonicity of generalized energy) \label{Thm-I-mono}
Suppose that $P =\calo(T-t)$,
$|\beta_k|+ D(t) = \calo((T-t)^2)$.
Then, there exist $C_1, C_2>0$
such that for $A$ large enough and  $t\in[t^*,T_*]$,
\begin{align} \label{dIt-mono-case2}
\frac{d \mathscr{I}}{dt}
\geq C_1
     \sum_{k=1}^{K}\frac{\g_k}{\lambda^2_{k}}
     \int (|\nabla R_{k}|^2+\frac{1}{\lbb_{k}^2} |R_{k}|^2 )
     e^{-\frac{|x-\alpha_{k}|}{A\lbb_k}}dx
    -C_2 A \mathcal{E}_r.
\end{align}
where
\begin{align} \label{calEr-def}
   \mathcal{E}_r
   =& \sum_{k=1}^{K} \frac{\lbb_k \dot{\lbb}_k + \g_k}{\lbb_k^{4} } M_k
      + \( \frac{Mod}{(T-t)^3} + \a^* (T-t)^{m-3+\frac d2} + (T-t)^{\upsilon_*-3}\)  D \nonumber \\
    &   + \frac{D^2}{(T-t)^2}
        +  \ve  \frac{D^2}{(T-t)^3}
          + \frac{M_k^2}{(T-t)^3}
      + e^{-\frac{\delta}{T-t}}.
\end{align}
\end{theorem}

The functionals $\scri^{(1)}$ and $\scri^{(2)}$ will be treated in
Lemmas \ref{Lem-I1t} and \ref{Lem-I2t} below, respectively.

\begin{lemma} (Control of $\mathscr{I}^{(1)}$)  \label{Lem-I1t}
Consider the situations as in Theorem \ref{Thm-I-mono}.
Then, there exists $C>0$ such that for any $t\in[t^*,T_*]$,
\begin{align} \label{I1t-case1}
 \frac{d \mathscr{I}^{(1)}}{dt}
\geq &\sum_{k=1}^{K}\frac{\gamma_k}{\lambda_k^4} \|\ve_k\|_{L^2}^2
     -\sum_{k=1}^{K}\frac{\gamma_k}{\lambda_k^4} {\rm Re}\int (1+\frac 2d)|Q_k|^\frac 4d|\ve_k|^2
      +\frac 2 d   |Q_k|^{\frac 4d -2} \ol{Q_k}^2\ve_k^2  dy   \nonumber  \\
&-\sum_{k=1}^{K}\frac{\gamma_k}{\lambda^4_k} {\rm Re}
      \int y \cdot\nabla \ol{Q_k}
       \ (f''(Q_k)\cdot \ve_k^2 ) dy
       -C \mathcal{E}_r,
\end{align}
where $\cale_r$ is the error as in \eqref{calEr-def}
but without the term $\frac{M_k^2}{(T-t)^3}$.
\end{lemma}

{\bf Proof.}
Let $\eta$ be as in \eqref{etan-Rn}.
Using equation \eqref{equa-R} and \eqref{g-gzz-expan}
we compute as in \cite[(5.32)]{SZ20} that
\begin{align} \label{equa-I1t}
 \frac{d \scri^{(1)}}{dt}
=& -\sum_{k=1}^K \dot{\lbb}_{k}\lbb_{k}^{-3} \int|R|^2\Phi_kdx
  -\sum_{k=1}^K \lbb_{k}^{-2} {\rm Im} \<f^\prime(U+z)\cdot R, R_k\>  \nonumber\\
  &  -{\rm Re} \< f''(U+z,R)\cdot R^2,  \partial_t (U+z) \>
    - \sum_{k=1}^K \lbb_{k}^{-2}{\rm Im} \< R \nabla\Phi_k, \na R\> \nonumber\\
  &-\sum_{k=1}^K \lbb_{k}^{-2} {\rm Im} \< f''(U+z,R)\cdot R^2, R_k\> \nonumber \\
& - {\rm Im} \<\Delta R -\sum\limits_{k=1}^K \lbb_k^{-2} R_k +f(v)-f(U+z), a_1\cdot \nabla R+a_0 R  \> \nonumber \\
 & -{\rm Im} \<\Delta R - \sum_{k=1}^K \lbb_{k}^{-2} R_k +f(v)-f(U+z), \eta\>  \nonumber\\
=& : \sum\limits_{l=1}^7 \scri^{(1)}_{t,l}.
\end{align}
The terms $\{\scri^{(1)}_{t,l}\}$ are estimated as follows:

{\it $(i)$ Estimate of $\scri^{(1)}_{t,1}$.}
Since $P=\calo(T-t)$, $D=\calo((T-t)^2)$,
by Theorem \ref{Thm-Mod-bdd},
\begin{align}   \label{Mod-t2}
   Mod = \calo((T-t)^2).
\end{align}
Hence, we compute
\begin{align*}
     - \frac{\dot{\lbb}_k}{\lbb_k^3} = \frac{\g_k}{\lbb_k^4}
   - \frac{\lbb_k\dot \lbb_k+\g_k}{\lbb_k^4}
   =  \frac{\g_k}{\lbb_k^4}  + \calo\(\frac{Mod}{\lbb_k^4}\)
   =  \frac{\g_k}{\lbb_k^4}  + \calo\(\frac{D^2}{(T-t)^2}\),
\end{align*}
which yields that
\begin{align} \label{I1t1-esti}
   \scri^{(1)}_{t,1}
  =&  \sum_{k=1}^{K} \frac{\g_k}{\lbb_k^4}\int |R|^2\Phi_kdx
     +\calo\(\frac{D^2}{(T-t)^2}\).
\end{align}

{\it $(ii)$ Estimates of $\scri^{(1)}_{t,2}$ and  $\scri^{(1)}_{t,3}$.}
Rewrite
\begin{align*}
    \scri^{(1)}_{t,2} +  \scri^{(1)}_{t,3}
    =& - \sum\limits_{k=1}^K \lbb_k^{-2}
      {\rm Im}  \< f'(U)\cdot R, R_k\>
      - {\rm Re} \<f''(U, R)\cdot R^2, \partial_t U\>  + er  ,
\end{align*}
where $er$ denotes the difference
\begin{align} \label{er-J1-J2}
  er:=& -\sum\limits_{k=1}^K \lbb_k^{-2} {\rm Im} \<f'(U+z)\cdot R - f'(U)\cdot R, R_k\> \nonumber \\
       & - \( {\rm Re} \<f''(U+z, R)\cdot R^2, \partial_t (U+z)\>
         - {\rm Re} \<f''(U, R)\cdot R^2, \partial_t U\>\)   \nonumber \\
      =:& er_1 + er_2.
\end{align}

We claim that
\begin{align} \label{er-Tt}
   er= \calo\( \a^* (T-t)^{-2}{D^2} \).
\end{align}
To this end,
by \eqref{f'v1-f'v2}, the renormalized variable $\ve_{z,k}$ in \eqref{ez-def},
\begin{align} \label{er-J1}
  |er_1| \leq& C (T-t)^{-2} \int (|U|^{\frac 4d-1} + |z|^{\frac 4d-1}) |z| |R|^2 dx \nonumber \\
       \leq& C \( (T-t)^{-4} \sum\limits_{k=1}^K \|e^{-\delta |y|}\ve_{z,k}\|_{L^\9} \|R\|_{L^2}^2
             +  (T-t)^{-2}  \|z \|_{L^\9}^{\frac 4d} \|R\|_{L^2}^2
             + e^{-\frac{\delta}{T-t}}  \|R\|_{L^2}^2 \) \nonumber \\
       \leq& C \(\a^* (T-t)^{m-3+\frac d2} + \a^* (T-t)^{-2} \) D^2.
\end{align}
Regarding the second term $er_2$, note that
\begin{align*}
   er_2
       =&   {\rm Re} \<f''(U+z, R)\cdot R^2, \partial_t z\>
          +  {\rm Re} \<f''(U+z, R)\cdot R^2 - f''(U, R)\cdot R^2, \partial_t U\> \nonumber \\
       =& : er_{21} + er_{22}.
\end{align*}
By \eqref{f''-bdd}
and estimates  \eqref{U-Lp}, \eqref{z-LtHm-a*}, \eqref{z-CtL2-a*}
and \eqref{R-Lp-D},
\begin{align}\label{er-J21}
  |er_{21}|
  \leq& C \int (|U|^{\frac 4d-1} +|z|^{\frac 4d-1} + |R|^{\frac 4d-1}) |R|^2 |\partial_tz| dx \nonumber \\
  \leq& C \|\partial_t z\|_{L^2}
          \(\|U\|_{L^{4(\frac 4d-1)}}^{\frac 4d-1}\|R\|_{L^8}^2
            + \|z\|_{L^\9}^{\frac 4d-1} \|R\|_{H^1}^{2}
            + \|R\|_{L^{\frac 8d+2}}^{\frac 4d +1}\) \nonumber \\
  \leq&  C \a^*  (T-t)^{-2}D^2.
\end{align}
Moreover,
since
by \eqref{equa-Ut} and $Mod =\calo(1)$,
\begin{align}  \label{pt-Uj-tx}
    \|\partial_t U_k(t)\|_{L^\9} \leq C (T-t)^{-\frac d2-2}.
\end{align}
Then, using \eqref{f''v1-f''v2}, \eqref{eyvez-bdd} and \eqref{R-Lp-D} we get
\begin{align} \label{er-J22}
   | er_{22} |
   \leq& C \int(|U|^{\frac 4d-2} + |R|^{\frac4d-2}+ |z|^{\frac 4d-2}) |z| |R|^2 |\partial_t U| dx \nonumber \\
   \leq& C (T-t)^{-4} \sum\limits_{k=1}^K \|e^{-\delta|y|}\ve_{z,k}\|_{L^\9} \|R\|_{L^2}^2
         + C (T-t)^{-d-2}  \sum\limits_{k=1}^K \|e^{-\delta|y|}\ve_{z,k}\|_{L^\9} \|R\|_{L^{\frac 4d}}^{\frac 4d}
         + C e^{-\frac{\delta}{T-t}} \|R\|_{L^2}^2  \nonumber \\
   \leq& C (T-t)^{-4} \sum\limits_{k=1}^K \|e^{-\delta|y|}\ve_{z,k}\|_{L^\9} D^2
         +   C e^{-\frac{\delta}{T-t}}D^2  \nonumber \\
   \leq& C \a^* (T-t)^{m-3+\frac d2} D^2.
\end{align}
Hence, plugging \eqref{er-J1}, \eqref{er-J21} and \eqref{er-J22} into \eqref{er-J1-J2}
we prove \eqref{er-Tt}, as claimed.

Thus,  since
$$|\beta_k|+ D(t) + Mod(t) =\calo((T-t)^2),$$
computing as in the proof of \cite[(4.18),(4.20)]{CSZ21} we obtain
\begin{align} \label{I1t2t3-esti}
   \scri^{(1)}_{t,2} +  \scri^{(1)}_{t,3}
    =&- \sum\limits_{k=1}^K \frac{\gamma_k}{\lambda_k^2} {\rm Re}
      \int (1+\frac 2d)|U_k|^{\frac {4}{d}}|R_k|^2
     + \frac 2 d |U_k|^{\frac 4d-2} \ol{U_k}^2 R_k^2 dx \nonumber \\
 &- \sum\limits_{k=1}^K\frac{\gamma_k}{\lambda_k} {\rm Re}
    \int (\frac{x-\alpha_k}{\lambda_k})\cdot\nabla \ol{U_k}
     \ f''(U_k)\cdot R_k^2 dx
        +\calo ( (T-t)^{-2} D^2(t)).
\end{align}

{\it $(iii)$ Estimate of $\scri^{(1)}_{t,4}$.}
We use the smallness in  {\rm Case (I)} or {\rm Case (II)}
to control this term
and have that (see  \cite[(4.21),(4.22)]{CSZ21})
\begin{align} \label{I1t4-esti-case2}
     |\scri^{(1)}_{t,4}|
     \leq C \ve (T-t)^{-3} {D^2}.
\end{align}

{\it $(iv)$ Estimate of $\scri^{(1)}_{t,5}$.}
Using \eqref{f''-bdd},  \eqref{U-Lp},
\eqref{R-Lp-D} and  $D\leq C(T-t)^2$
we get
\begin{align} \label{I155.0}
   | \scri^{(1)}_{t,5}|
   \leq& C (T-t)^{-2} \( \|U\|_{L^{2(\frac 4d-1)}}^{\frac 4d-1} \|R\|_{L^6}^3
            + \|z\|_{L^\9}^{\frac 4d-1} \|R\|_{L^3}^3
            + \|R\|_{L^{\frac 4d +2}}^{\frac 4d +2}\)  \nonumber  \\
   \leq& C (T-t)^{-2} \((T-t)^{-2} D^3 + \a^* (T-t)^{-\frac d2} D^3  + (T-t)^{-2} D^{\frac 4d+2}\) \nonumber \\
   \leq& C (T-t)^{-2}  D^2.
\end{align}

{\it $(v)$ Estimate of $\scri^{(1)}_{t,6}$.}
Using the integration by parts formula and
\eqref{fv1-fv2} we get
(see also \cite[(5.42)]{SZ20})
\begin{align} \label{DR-bnaR.1}
       & |{\rm Im}\<\Delta R - \sum\limits_{k=1}^K \lbb_k^{-2} R_k + f(v) - f(U+z), a_1 \cdot \na R\>|  \nonumber \\
    \leq & C\(\|\na R\|_{L^2}^2 + (T-t)^{-2} \|R\|_{L^2}^2\)
         + C   \int (|U|^{\frac 4d} + |z|^{\frac 4d} + |R|^{\frac 4d})  |R| |a_1\na R|  dx.
\end{align}
Then, by \eqref{a1-loworder}, the change of variables and
\eqref{phik-Taylor},
\begin{align}  \label{I1t6-b.3}
       & \int (|U|^{\frac 4d}+ |z|^{\frac 4d} + |R|^{\frac 4d})  |R| |a_1 \cdot \na R|  dx    \nonumber \\
  \leq& C (T-t)^{-2} \sum\limits_{l=1}^N \sum\limits_{k=1}^K
            \|e^{-\delta|y|}  \nabla\phi_l (\lbb_ky + \a_k) \|_{L^\9} \|R\|_{L^2} \|\na R\|_{L^2}     \nonumber \\
      &  + C \(\|z\|_{L^\9}^\frac 4d \|R\|_{L^2} \|\na R\|_{L^2}
              + \|R\|_{L^{\frac 8d +2}}^{\frac 4d+1} \|\na R\|_{L^2} \)
        + C e^{-\frac{\delta}{T-t}} \|R\|_{L^2}^2   \nonumber \\
  \leq& C \( (T-t)^{\upsilon_*-3} + \a^*(T-t)^{-1}\)  D^2
        + C (T-t)^{-3}D^{2+\frac 4d}
        + C e^{-\frac{\delta}{T-t}} D^2 .
\end{align}
Since $D \leq C (T-t)$, we come to
\begin{align}\label{I1t6-b}
   |{\rm Im} \<\Delta R -   \sum\limits_{k=1}^K \lbb_k^{-2} R_k + f(u) - f(U+z), a_1 \cdot \na R\>|
    \leq& C  (T-t)^{-2} D^2 .
\end{align}

Similarly, we have
\begin{align} \label{I1t6-c}
        |{\rm Im} \<\Delta R-  \sum\limits_{k=1}^K \lbb_k^{-2} R_k + f(v) - f(U+z),a_0 R\>|
   \leq& C\( \|R\|^2_{H^1} + \lbb_k^{-2}\|R\|_{L^2}^2
         + \int (|U|^{\frac 4d}+|z|^{\frac 4d} +|R|^{\frac 4d} )  |R|^2 dx \)  \nonumber  \\
   \leq& C \( (T-t)^{-2} D^2 + \a^* D^2 + (T-t)^{-2}D^{2+\frac 4d} \) \nonumber  \\
   \leq&  C (T-t)^{-2} D^2.
\end{align}

Thus, we conclude from \eqref{I1t6-b}  and \eqref{I1t6-c}  that
\begin{align} \label{I1t6-esti}
    |\scri^{(1)}_{t,6}|
  \leq  C   (T-t)^{-2} D^2.
\end{align}

{\it $(vi)$ Estimate of $\scri^{(1)}_{t,7}$.}
It remains to treat the  delicate inner product $\scri^{(1)}_{t,7}$
involving the $\eta$ term.
First, we claim that
\begin{align} \label{I1t7.0}
  \scri^{(1)}_{t,7}
  =& - \sum\limits_{k=1}^K {\rm Im}
  \< \Delta R_k- \lbb_k^{-2} R_k + f'(U_k)\cdot R_k, \eta\> \nonumber \\
  & + \calo\( \a^* (T-t)^{m-1+\frac d2}  D + (T-t)^{-2} D^2 + e^{-\frac{\delta}{T-t}}\).
\end{align}
This means that, the inner products involving $R_k$ of orders higher than one
are acceptable errors.

To this end,
we use the expansion \eqref{g-gzz-expan} to get
\begin{align} \label{fu-fUz}
  f(v)-f(U+z)
  =& f'(U+z)\cdot R + f''(U+z, R)\cdot R^2 \nonumber \\
  =& f'(U) \cdot R + \(f'(U+z)\cdot R - f'(U)\cdot R\) + f''(U+z, R)\cdot R^2.
\end{align}
Define the renormalized variable $\ve_{R,k}$   by
\begin{align*}
   R(t,x) = \lbb_k^{-\frac d2} \ve_{R,k} (t,\frac{x-\a_k}{\lbb_k})e^{i\theta_k}.
\end{align*}
Then,  by \eqref{f'v1-f'v2} and \eqref{eta-tx},
\begin{align} \label{I1t7-f'Uz}
       & {\rm Im}\<f'(U+z) \cdot R -f'(U)\cdot R, \eta \> \nonumber \\
   \leq&  C \int (|U|^{\frac 4d-1} + |z|^{\frac 4d-1}) |z| |R| |\eta| dx \nonumber \\
   \leq& C (T-t)^{-4}  \sum\limits_{k=1}^K
           \int \(e^{-\delta|y|}+|\ve_{z,k}|^{\frac 4d-1}\) |\ve_{z,k}|
           |\ve_{R,k}| \(Mod + |\ve_{z,k}| + (T-t)^{\upsilon_*+1}\) e^{-\delta |y|} dy
           + Ce^{-\frac{\delta}{T-t}} \nonumber \\
   \leq& C (T-t)^{-4}  \sum\limits_{k=1}^K  \|e^{-\delta|y|} \ve_{z,k} \|_{L^\9} \|R\|_{L^2}
               \(Mod + \|e^{-\delta|y|}\ve_{z,k}\|_{L^\9} +(T-t)^{\upsilon_*+1}\)
         + Ce^{-\frac{\delta}{T-t}}  \nonumber  \\
   \leq&C  \a^* (T-t)^{m-1+\frac d2} D  + Ce^{-\frac{\delta}{T-t}},
\end{align}
Moreover,
by \eqref{U-Lp}, \eqref{R-Lp-D}, \eqref{eta-L2}, \eqref{Mod-t2}
and $Mod\leq C(T-t)^2$,
\begin{align} \label{I1t7-f''Uz}
  {\rm Im}\<f''(U+z,R)\cdot R^2, \eta\>
  \leq& C  \|U+z\|_{L^{4(\frac 4d -1)}}^{\frac 4d-1} \|R\|_{L^8}^2 \|\eta\|_{L^2}
       +C \|R\|_{H^1}^{\frac 4d+1} \|\eta\|_{L^2} \nonumber\\
  \leq& C  (T-t)^{-2}\|\eta\|_{L^2} D^2  + C \|\eta\|_{L^2} \|R\|_{H^1}^2 \nonumber \\
  \leq& C (T-t)^{-2} D^2.
\end{align}
Thus, combining \eqref{fu-fUz}, \eqref{I1t7-f'Uz} and \eqref{I1t7-f''Uz}
and using Lemma \ref{Lem-decoup-U} we obtain \eqref{I1t7.0}, as claimed.

Next,
in order to treat the remaining linear terms on the R.H.S. of \eqref{I1t7.0},
we decompose $\eta$ into four parts $\eta = \sum_{l=1}^4 \eta_l$
as in \eqref{eta1-eta4}.
Note that, by Lemma \ref{Lem-decoup-U},  \eqref{eta2-def} and \eqref{eta3-def},
\begin{align}  \label{eta2-eta4-I1-esti}
      & | {\rm Im} \<\Delta R_k - \lbb_k^{-2}R_k + f'(U_k)\cdot R_k, \eta_2+\eta_3\> |    \nonumber \\
  \leq&   C \lbb_k^{-4} \sum\limits_{j=1}^{\frac 4d}
             \int |\na \ve_k| (|\ve_{z,k}|^j + |\na \ve_{z,k}||\ve_{z,k}|^{j-1}) e^{-\delta|y|} dy
        + C \sum\limits_{j=1}^{\frac 4d} \int |\ve_k| |\ve_{z,k}^j| e^{-\delta |y|} dy \nonumber \\
      &    + C \lbb_k^{-4} \sum\limits_{j=1}^{\frac 4d}
          \int |f'(Q_k)\ve_k| |Q_k|^{1+\frac 4d-j} |\ve_{z,k}|^{j} dy
           + C e^{-\frac{\delta}{T-t}}  \nonumber \\
  \leq& C \lbb_k^{-4} \|e^{-\delta |y|} (|\ve_{z,k}|+|\na \ve_{z,k}|) \|_{L^\9} D
         + C e^{-\frac{\delta}{T-t}}  \nonumber  \\
  \leq& C \a^* (T-t)^{m-3+\frac d2} D  + C e^{-\frac{\delta}{T-t}},
\end{align}
where the last step is due to \eqref{eyvez-bdd}.
Moreover, by \eqref{eta4-def} and Lemma \ref{Lem-wtbwtc},
\begin{align}
    & | {\rm Im} \<\Delta R_k - \lbb_k^{-2}R_k + f'(U_k)\cdot R_k, \eta_4\> |    \nonumber \\
    \leq&  C \lbb_k^{-2} |{\rm Im} \<\na \ve_k, \na(\lbb_k^{-1} \wt a_{1,k}\cdot \na Q_k + \wt a_{0,k} Q_k)\>|    \nonumber \\
         & + C \lbb_k^{-2} \int (|\ve_k|+|f'(Q_k)\ve_k|) | \lbb_k^{-1}  \wt a_{1,k}\cdot \na Q_k + \wt a_{0,k} Q_k|dy     \nonumber \\
    \leq& C \lbb_k^{\upsilon_*-3} \( \|\ve_k\|_{L^2} + \|\na \ve_k\|_{L^2} \) \nonumber \\
    \leq& C (T-t)^{\upsilon_*-3} D.
\end{align}
Hence, taking into account \eqref{eta1-def}
and using Lemma \ref{Lem-decoup-U} again we obtain
\begin{align} \label{Rlinear-eta}
      {\rm Im} \<\Delta R_k - \lbb_k^{-2} R_k + f'(U_k)\cdot R_k, \eta\>
  =& \lbb_k^{-4} {\rm Im} \<\Delta \ve_k - \ve_k + f'(Q_k)\cdot \ve_k, \Psi_k\> \nonumber \\
   &  + \calo\( (\a^*  (T-t)^{m-3+\frac d2}
     +  (T-t)^{\upsilon_*-3} ) D +e^{-\frac{\delta}{T-t}} \),
\end{align}
where $\Psi_k$ is given by \eqref{Psik-def}.

The analysis is now reduced to that of the inner product involving  $\Psi_k$.

By the proximity
\begin{align} \label{Qj-Q}
   Q_k = Q + \calo(P e^{-\delta |y|})
\end{align}
and the definition of linearized operators in \eqref{L+-L-},
\begin{align} \label{Rlinear-L-Psik}
   & {\rm Im} \<\Delta \ve_k -  \ve_k + f'(Q_k)\cdot \ve_k, \Psi_k\>   \nonumber \\
  =&   {\rm Im} \<\Delta \ve_k -  \ve_k + f'(Q)\cdot \ve_k, \wt \Psi_k\>
    + \calo(  P  Mod \|R\|_{L^2}) \nonumber \\
  =&  \<L_+\ve_{k,1}, \wt \Psi_{k,2}\> - \<L_-\ve_{k,2}, \wt \Psi_{k,1}\>
       + \calo( P Mod D).
\end{align}
where $\wt \Psi_k$ is defined as in \eqref{Psik-def}
with $Q$ replacing $Q_k$,
$\ve_{k,1} = {\rm Re} \ve_k$,
$\ve_{k,2} = {\rm Im} \ve_{k,2}$,
and $\wt \Psi_{k,1}$, $\wt \Psi_{k,2}$ are defined similarly.
Then, by \eqref{Psik-def} and the algebraic identities in \eqref{Q-kernel},
\begin{align}
   \<L_+\ve_{k,1}, \wt \Psi_{k,2}  \>
   =&-(\lbb_k \dot{\lbb}_k + \g_k)\<\ve_{k,1}, L_+\Lambda Q\> \nonumber \\
   =& 2(\lbb_k \dot{\lbb}_k + \g_k) {\rm Re}\<U_k, R_k\>
      + \calo( P Mod \|R\|_{L^2}),  \label{L+ve1-Psik2}
\end{align}
and
\begin{align}
    \<L_-\ve_{k,2}, \wt \Psi_{k,1}\>
   = & -(\lbb_k^2 \beta_k + \g_k \beta_k) \<\ve_{k,2}, L_-Q\>
        + \frac 14 (\lbb_k^2 \dot{\g}_k + \g_k^2) \<\ve_{k,2}, L_-|y|^2 Q\> \nonumber \\
   = & 2(\lbb_k^2 \beta_k + \g_k \beta_k) \<\ve_{k,2}, \na Q\>
        -(\lbb_k^2 \dot{\g}_k + \g_k^2) \<\ve_{k,2}, \Lambda Q\> \nonumber \\
   = &\calo\( P Mod \|R\|_{L^2} + e^{-\frac{\delta}{T-t}} \|R\|_{L^2} \), \label{L-ve2-Psik1}
\end{align}
where in the last step we also used the almost orthogonality \eqref{Orth-almost}.
Hence, plugging \eqref{Rlinear-L-Psik}, \eqref{L+ve1-Psik2} and \eqref{L-ve2-Psik1} into \eqref{Rlinear-eta}
we obtain
\begin{align} \label{Rlinear-eta-esti}
     & {\rm Im} \<\Delta R_k - \lbb_k^{-2} R_k + f'(U_k)\cdot R_k, \eta \> \nonumber \\
   =&  2 \lbb_k^{-4} (\lbb_k \dot{\lbb}_k + \g_k)
         {\rm Re}\<U_k, R_k\>    \nonumber \\
    & + \calo\( \((T-t)^{-3} Mod + \a^* (T-t)^{m-3+\frac d2} + (T-t)^{\upsilon_*-3}\) D + e^{-\frac{\delta}{T-t}}\).
\end{align}

Thus, combining \eqref{I1t7.0} and \eqref{Rlinear-eta-esti} together we arrive at
\begin{align} \label{I1t7-esti}
   \scri^{(1)}_{t,7}
   =& - 2 \sum_{k=1}^{K}\lbb_k^{-4} (\lbb_k \dot{\lbb}_k + \g_k)
        {\rm Re} \<U_k, R_k\>  \nonumber \\
    & + \calo\( \( (T-t)^{-3} Mod + \a^*  (T-t)^{m-3+\frac d2} +  (T-t)^{\upsilon_*-3}\) D  +  (T-t)^{-2}  D^2 +  e^{-\frac{\delta}{T-t}} \).
\end{align}

Finally,
plugging estimates \eqref{I1t1-esti}, \eqref{I1t2t3-esti},  \eqref{I1t4-esti-case2},
\eqref{I155.0}, \eqref{I1t6-esti} and \eqref{I1t7-esti}
into \eqref{equa-I1t}
we obtain \eqref{I1t-case1}
and finish the proof of Lemma \ref{Lem-I1t}.
\hfill $\square$

\begin{lemma} (Control of $\scri^{(2)}$) \label{Lem-I2t}
Consider the situations as in Theorem \ref{Thm-I-mono}.
Then, there exists $C>0$ such that for all $t\in[t^*,T_*]$,
\begin{align} \label{I2t}
 \frac{d \scri^{(2)}}{dt}
 \geq &  - \sum_{k=1}^{K} \frac{\gamma_k}{4\lambda_k^4}\int \Delta^2\chi_A(y)|\ve_k|^2 dy
           + \sum_{k=1}^{K} \frac{\gamma_k}{\lambda_k^4} {\rm Re} \int \nabla^2\chi_A(y)(\nabla \ve_k,\nabla \ol{\ve_k}) dy   \nonumber \\
&+ \sum_{k=1}^{K}\frac{\gamma_k}{\lambda^4_k}{\rm Re}\int
   \nabla\chi_A (y)\cdot\nabla \ol{Q_k}
   \ f''(Q_k)\cdot \ve_k dy
   - C A \cale_r',
\end{align}
where
\begin{align}  \label{caler'-def}
    \cale_r' =  \( \frac{Mod}{(T-t)^{3}} + \a^* (T-t)^{m-2+\frac d2} +  (T-t)^{\upsilon_*-2}\) D
                +  \frac{D^2}{(T-t)^{2}}
       + e^{-\frac{\delta}{T-t}} .
\end{align}
\end{lemma}

{\bf Proof.}
We compute as in \cite[(5.48)]{SZ20},
\begin{align}\label{equa-I2t}
   \frac{d \scri^{(2)}}{dt}
   =&-\sum_{k=1}^K \frac{\dot{\lambda}_{k}\gamma_{k}-\lambda_{k}\dot{\gamma}_{k}}{2\lambda_{k}^2}
           {\rm Im} \< \nabla\chi_A(\frac{x-\alpha_{k}}{\lambda_{k}})\cdot\nabla R, R_k\> \nonumber  \\
   &+\sum_{k=1}^K\frac{\gamma_{k}}{2\lambda_{k}}{\rm Im}
       \< \partial_t (\nabla\chi_A(\frac{x-\alpha_{k}}{\lambda_{k}}))\cdot\nabla R, {R}_k\>
    +\sum_{k=1}^K \frac{\gamma_{k}}{2\lambda_{k}^2}   {\rm Im}
       \< \Delta\chi_A(\frac{x-\alpha_{k}}{\lambda_{k}})R_k, \pa_t {R} \>  \nonumber \\
   &+   \sum_{k=1}^K \frac{\gamma_{k}}{2\lambda_{k}}    {\rm Im}
       \< \nabla\chi_A(\frac{x-\alpha_{k}}{\lambda_{k}})\cdot ( \nabla R_k + \na R \Phi_k), \partial_t  R \> \nonumber \\
    =:&  \sum\limits_{k=1}^K
         \(\scri_{t,k1}^{(2)} +\scri_{t,k2}^{(2)} + \scri_{t,k3}^{(2)} +\scri_{t,k4}^{(2)}\).
\end{align}

$(i)$ {\it Estimate of $ \scri^{(2)}_{t,k1}$ and $\scri^{(2)}_{t,k2}$.}
Since
$|\frac{\dot{\lambda}_k\gamma_k-\lambda_k\dot{\gamma}_k}{\lambda_k^2}|
\leq C \lambda_k^{-3} {Mod_k}$
and $|\partial_t (\na \chi_A(\frac{x-\a_k}{\lbb_k}))| \leq CA \lbb_k^{-2} (Mod+P)$,
by \eqref{Mod-t2},
\begin{align}  \label{esti-I2t1t2}
 |\scri^{(2)}_{t,k1} + \scri^{(2)}_{t,k2}|
 \leq& CA  \lbb_k^{-3}(Mod + P^2) \|\na R\|_{L^2} \|R\|_{L^2} \nonumber \\
 \leq& CA  (\lbb_k^{-4}Mod D^2 +\lbb_k^{-2} D^2)
 \leq CA  (T-t)^{-2}D^2 .
\end{align}

$(ii)$ {\it Estimate of $ \scri^{(2)}_{t,k3}$.}
We claim that
\begin{align} \label{I2t3.1-esti}
 \scri^{(2)}_{t,k3}
 =& -\frac{\gamma_k}{4\lambda_k^4}{\rm Re}\int \Delta^2\chi_A (\frac{x-\alpha_k}{\lambda_k} )|R_k|^2 dx+\frac{\gamma_k}{2\lambda_k^2}{\rm Re}\int \Delta\chi_A(\frac{x-\alpha_k}{\lambda_k})|\nabla R_k|^2 dx  \nonumber \\
&-\frac{\gamma_{k}}{2\lambda_{k}^2}
   {\rm Re} \< \Delta\chi_A(\frac{x-\alpha_{k}}{\lambda_{k}})R_k, f'(U_k)\cdot R_k \> \nonumber \\
  & + \calo\(A (T-t)^{-2} D^2
               +  \((T-t)^{-3} Mod + \a^* (T-t)^{m-2+\frac d2} + (T-t)^{\upsilon_*-2}\)  D  + e^{-\frac{\delta}{T-t}}\).
\end{align}

For this purpose,
by \eqref{equa-R} and  \eqref{g-gzz-expan},
\begin{align} \label{Deltachi-esti}
 \mathscr{I}^{(2)}_{t,k3}
      =& - \frac{\gamma_{k}}{2\lambda_{k}^2}
    {\rm Re}  \< \Delta\chi_A(\frac{x-\alpha_{k}}{\lambda_{k}})R_k,
       \Delta R+f'(U+z)\cdot R +f''(U+z,R)\cdot R^2
       + a_1 \cdot \na R+ a_0 R +\eta\>.
\end{align}

First, we have from \cite[(3.67)]{CSZ21} that
\begin{align} \label{Dchi-Rj-DR}
    -\frac{\g_k}{2\lbb_k^2} {\rm Re} \< \Delta\chi_A(\frac{x-\alpha_{k}}{\lambda_{k}})R_k, \Delta R\>
 =& - \frac{\g_k}{4\lbb_k^4} {\rm Re} \int \Delta^2 \chi_A(\frac{x - \a_k}{\lbb_k}) |R_k|^2 dx  \nonumber  \\
  & + \frac{\g_k}{2\lbb_k^2} {\rm Re} \int \Delta \chi_A(\frac{x - \a_k}{\lbb_k}) |\na R_k|^2 dx
   + \calo(A\|R\|_{H^1}^2).
\end{align}
Let us mention that,
an extra factor $T-t$ is gained here from the
decay properties of the cut-off function,
i.e., for $|y|\geq 2A$,
\begin{align}
    |\na \Delta\chi_A(y)| \leq C A |y|^{-2},\ \
    |\partial_{x_kx_l} \chi_A(y)| \leq C A |y|^{-1}, \ \ 1\leq k,l\leq d.
\end{align}

Moreover, rewrite
\begin{align*} \label{Dchi-Rj-f'UR}
     {\rm Re} \< \Delta\chi_A(\frac{x-\alpha_{k}}{\lambda_{k}})R_k, f'(U+z)\cdot R \>
   =& {\rm Re} \< \Delta\chi_A(\frac{x-\alpha_{k}}{\lambda_{k}})R_k, f'(U)\cdot R \>
      + \wt er,
\end{align*}
where the difference
\begin{align}
    \wt er :=& {\rm Re} \< \Delta\chi_A(\frac{x-\alpha_{k}}{\lambda_{k}})R_k,  f'(U+z) \cdot R - f'(U)\cdot R\> .
\end{align}
Then, by the bound $\|\Delta\chi_A\|_{L^\9} \leq C$, \eqref{z-LtHm-a*} and \eqref{f'v1-f'v2},
\begin{align*}
   |\wt er|
   \leq& C \int (|U|^{\frac 4d-1} + |z|^{\frac 4d-1})|z||R|^2 dx \\
   \leq& C \( (T-t)^{-\frac d2(\frac 4d-1)} + \|z\|_{L^\9}^{\frac 4d-1}\) \|z\|_{L^\9} \|R\|_{L^2}^2 \\
   \leq& C (T-t)^{-2} D^2.
\end{align*}
This along with Lemma \ref{Lem-decoup-U}  yields that
\begin{align} \label{I2jt3-DeltachiA}
   {\rm Re} \<\Delta \chi_A (\frac{x-\a_k}{\lbb_k}) R_k, f'(U+z)\cdot R\>
  = {\rm Re} \<\Delta \chi_A (\frac{x-\a_k}{\lbb_k}) R_k, f'(U_k)\cdot R_k\>
    + \calo\((T-t)^{-2} D^2 + e^{-\frac{\delta}{T-t}}\).
\end{align}

It also follows from \eqref{f''-bdd}, \eqref{U-Lp}, \eqref{R-Lp-D}
and $D=\calo((T-t)^2)$ that
\begin{align} \label{Dchi-Rj-f''UR2}
        \bigg|\frac{\g_k}{2\lbb_k^2} {\rm Re} \< \Delta\chi_A(\frac{x-\alpha_{k}}{\lambda_{k}})R_k, f''(U+z,R)\cdot R^2\> \bigg|
 \leq& C \lbb_k^{-1} (\|R\|_{L^6}^3 \|U+z\|_{L^{2(\frac 4d-1)}}^{\frac 4d-1} + \|R\|_{L^{\frac 4d+2}}^{\frac 4d+2}) \nonumber \\
 \leq& C (T-t)^{-1} \((T-t)^{-2} D^3 + (T-t)^{-2} D^{\frac 4d+2}\)  \nonumber \\
 \leq& C (T-t)^{-2}   D^2.
\end{align}

Furthermore,
by \eqref{eta-L2},
\begin{align} \label{Dchi-Rj-bcR-eta}
       & \bigg|\frac{\g_k}{2\lbb_k^2} {\rm Re}
       \< \Delta\chi_A(\frac{x-\alpha_{k}}{\lambda_{k}})R_k,  a_1\cdot \na R+a_0 R +\eta \> \bigg|   \nonumber \\
  \leq & C \lbb_k^{-1}(\|R\|_{L^2}\|\nabla R\|_{L^2}+\|R\|_{L^2}^2)
        +C \lbb_k^{-1} \|R\|_{L^2} \|\eta\|_{L^2}   \nonumber \\
  \leq& C \((T-t)^{-2} D^2 + \((T-t)^{-3} Mod + \a^* (T-t)^{m-2+\frac d2} + (T-t)^{\upsilon_*-2}\) D \)
\end{align}

Hence,
plugging \eqref{Dchi-Rj-DR}, \eqref{I2jt3-DeltachiA},
\eqref{Dchi-Rj-f''UR2} and \eqref{Dchi-Rj-bcR-eta} into \eqref{Deltachi-esti}
we obtain \eqref{I2t3.1-esti}, as claimed.

$(iii)$ {\it Estimate of $ \scri^{(2)}_{t,k4}$.}
The estimate of $ \scri^{(2)}_{t,k4}$ is similar to that of $ \scri^{(2)}_{t,k3}$.
We claim that
\begin{align}   \label{I2t3.2-esti}
  \scri^{(2)}_{t,k4} =&\frac{\gamma_k}{\lambda_k^2}{\rm Re}\int \nabla^2\chi_A(\frac{x-\alpha_k}{\lambda_k})(\nabla R_k,\nabla \ol{R_k}) dx
-\frac{\gamma_k}{2\lambda_k^2}{\rm Re}\int \Delta\chi_A(\frac{x-\alpha_k}{\lambda_k})|\nabla R_k|^2 dx \nonumber  \\
&-\frac{\gamma_{k}}{\lambda_{k}}\< \nabla\chi_A(\frac{x-\alpha_{k}}{\lambda_{k}})\cdot
         \nabla R_k,f^\prime(U_k)\cdot R_k\> \nonumber \\
     &+ \calo\(A(T-t)^{-2} D^2 + A\((T-t)^{-3} Mod+ \a^* (T-t)^{m-2+\frac d2} + (T-t)^{\upsilon_*-2}\) D\) .
\end{align}

For this purpose,
using  \eqref{equa-R} again and \eqref{g-gzz-expan} we derive
\begin{align} \label{I2t3.2-esti.0}
 \scri^{(2)}_{t,k4} =& -\frac{\gamma_{k}}{2\lambda_{k}} {\rm Re}
     \< \nabla\chi_A(\frac{x-\alpha_{k}}{\lambda_{k}})\cdot
         (\nabla R_k + \nabla R\Phi_k), \nonumber \\
 & \qquad \qquad \ \ \Delta R+f^\prime(U+z)\cdot R_k
+f''(U+z,R)\cdot R^2+(a_1\cdot \na +a_0) R+\eta\> .
\end{align}

Similarly to \eqref{Dchi-Rj-DR}, we have (see \cite[(3.73)]{CSZ21})
\begin{align} \label{nachi-naR-DR}
   & - \frac{\g_k}{2\lbb_k} {\rm Re}\< \nabla\chi_A(\frac{x-\alpha_{k}}{\lambda_{k}})\cdot
        (\nabla R_k + \nabla R \Phi_k ),  \Delta {R} \>    \nonumber \\
 =&   \frac{\g_k}{ \lbb^2_k}  {\rm Re} \int   \na^2 \chi_A (\frac{x-\alpha_{k}}{\lambda_{k}}) (\na R_k, \na \ol{R_k})
     -  \frac{\g_k}{2\lbb^2_k}  \Delta \chi_A (\frac{x-\alpha_{k}}{\lambda_{k}}) |\na R_k|^2 dx
     + \calo\(A(T-t)^{-2}D^2\).
\end{align}
We note that, the second order terms $\partial_{x_kx_l} \ol{R}$
are cancelled by the integration by parts formula.
For the detailed computations we refer to \cite[(5.66)]{SZ20}.

Moreover, as in \eqref{Dchi-Rj-f'UR},   rewrite
\begin{align} \label{nachiR-f'Uz}
     &\frac{\gamma_{k}}{2\lambda_{k}} {\rm Re}
     \< \nabla\chi_A(\frac{x-\alpha_{k}}{\lambda_{k}})\cdot
         (\nabla R_k + \nabla R\Phi_k), f'(U+z)\cdot R\>  \nonumber \\
    =& \frac{\gamma_{k}}{\lambda_{k}} {\rm Re}
     \< \nabla\chi_A(\frac{x-\alpha_{k}}{\lambda_{k}})\cdot
         \nabla R_k, f'(U_k)\cdot R_k\>
        + \wh er + \calo(A e^{-\frac{\delta}{T-t}}),
\end{align}
where the last step is due to Lemma \ref{Lem-decoup-U} and
the error term is of form
\begin{align*}
    \wh er :=& \frac{\gamma_{k}}{\lambda_{k}} {\rm Re}
     \< \nabla\chi_A(\frac{x-\alpha_{k}}{\lambda_{k}})\cdot
          (\na R_k+ \nabla R\Phi_k), f'(U+z)\cdot R - f'(U) \cdot R\>.
\end{align*}
We use the bound $\|\na \chi_A\|_{L^\9} \leq CA$, \eqref{z-LtHm-a*}, \eqref{eyvez-bdd} and \eqref{f'v1-f'v2} to bound
\begin{align} \label{wher}
   |\wh er|
  \leq &C A\int (|\na R|+|R|)(|U|^{\frac 4d-1} + |z|^{\frac 4d-1}) |z| |R| dx \nonumber \\
   \leq& CA \( \sum\limits_{k=1}^K (T-t)^{-2} \|e^{-\delta |y|} \ve_{z,k}\|_{L^\9} + \|z\|^{\frac 4d}_{L^\9}\)
           \int (|\na R|+|R|)|R| dx + CAe^{-\frac{\delta}{T-t}} \nonumber \\
   \leq& C A\(\a^*(T-t)^{m-1+\frac d2} + \a^*\) (T-t)^{-1} D^2 + CAe^{-\frac{\delta}{T-t}} \nonumber \\
   \leq& CA \( (T-t)^{-1} D^2 + e^{-\frac{\delta}{T-t}} \).
\end{align}
Plugging this into \eqref{nachiR-f'Uz} yields that
\begin{align} \label{nachiR-f'Uz*}
     &- \frac{\gamma_{k}}{2\lambda_{k}} {\rm Re}
     \< \nabla\chi_A(\frac{x-\alpha_{k}}{\lambda_{k}})\cdot
         (\nabla R_k + \nabla R\Phi_k), f'(U+z)\cdot R\>  \nonumber \\
    =& - \frac{\gamma_{k}}{\lambda_{k}} {\rm Re}
     \< \nabla\chi_A(\frac{x-\alpha_{k}}{\lambda_{k}})\cdot
         \nabla R_k, f'(U_k)\cdot R_k\>
         + \calo\( A(T-t)^{-1}D^2 +A e^{-\frac{\delta}{T-t}}\).
\end{align}

For the remaining inner products in \eqref{I2t3.2-esti.0},
by  \eqref{U-Lp} and \eqref{R-Lp-D},
\begin{align} \label{nachi-naR-f''UR2-a1a0R}
 & \bigg|\frac{\gamma_{k}}{2\lambda_{k}} {\rm Re}
     \< \nabla\chi_A(\frac{x-\alpha_{k}}{\lambda_{k}})\cdot
         (\nabla R_k + \nabla R\Phi_k),
         f''(U+z,R)\cdot R^2 +  a_1\cdot \na R +a_0 R\> \bigg|  \nonumber \\
  \leq& C A\int (|\na R|+|R|) (|U+z|^{\frac 4d-1} + |R|^{\frac 4d-1})|R|^2dx
         + C A \int (|\na R|+|R|)^2 dx \nonumber \\
  \leq& CA \( (T-t)^{-\frac d2(\frac 4d-1)} + \|z\|_{L^\9}^{\frac 4d-1} \)
          \int (|\na R|+|R|)|R|^2 dx
       + CA \int (|\na R|+|R|)|R|^{\frac 4d+1} dx + CA \|R\|_{H^1}^2 \nonumber \\
  \leq& CA (T-t)^{-2+\frac d2} \(\|\na R\|_{L^2} \|R\|_{L^4}^2 + \|R\|_{L^3}^3\)
        + CA \(\|\na R\|_{L^2} \|R\|_{L^{2(\frac 4d+1)}}^{\frac 4d+1} + \|R\|_{L^{\frac 4d+2}}^{\frac 4d+2}\)
        + C A\|R\|_{H^1}^2.
\end{align}
Then, by \eqref{R-Lp-D}, the R.H.S. above can be bounded by,
up to a universal constant $CA$,
\begin{align}\label{nachi-naR-f''UR2-a1a0R}
   & (T-t)^{-2+\frac d2} \( (T-t)^{-\frac d2-1} + (T-t)^{-\frac d2} \)D^3
   + \( (T-t)^{-3} + (T-t)^{-2} \)D^{\frac 4d+2}
   + (T-t)^{-2} D^2 \nonumber \\
   \leq&  (T-t)^{-3} D^3 + (T-t)^{-3} D^{2+\frac 4d}
         + (T-t)^{-2}D^2 \nonumber \\
   \leq&  (T-t)^{-2} D^2.
\end{align}

Finally, the last inner product involving $\eta$ can be bounded easier than the previous $\mathscr{I}^{(1)}_{t,7}$
in $\scri^{(1)}$.
By \eqref{eta-L2},
\begin{align} \label{nachi-naR-eta}
   & \bigg|\frac{\gamma_{k}}{2\lambda_{k}} {\rm Re}
     \< \nabla\chi_A(\frac{x-\alpha_{k}}{\lambda_{k}})\cdot
         (\nabla R_k + \nabla R\Phi_k), \eta\> \bigg| \nonumber \\
  \leq& C A\|R\|_{H^1} \|\eta\|_{L^2} \nonumber \\
   \leq&  C A\((T-t)^{-3} Mod + \a^*(T-t)^{m-2+\frac d2} + (T-t)^{\upsilon_*-2}\) D .
\end{align}

Hence, we conclude from \eqref{nachi-naR-DR}, \eqref{nachiR-f'Uz*}, \eqref{nachi-naR-f''UR2-a1a0R}
and \eqref{nachi-naR-eta} that \eqref{I2t3.2-esti} holds.

Now,
putting the estimates \eqref{esti-I2t1t2},
\eqref{I2t3.1-esti} and \eqref{I2t3.2-esti} altogether
and using the renormalized variable $\ve_k$ in \eqref{Rj-ej}
we arrive at
\begin{align*}
  \frac{d \scri^{(2)}}{dt}
 =&- \sum_{k=1}^K \frac{\gamma_k}{4\lambda_k^4}\int \Delta^2\chi_A (y)|\ve_k|^2 dy
  +\sum_{k=1}^K\frac{\gamma_k}{\lambda_k^4}{\rm Re}\int \nabla^2\chi_A(y)(\nabla \ve_k,\nabla \ol{\ve_k}) dy \nonumber \\
& -\sum_{k=1}^K{\rm Re} \<\frac{\gamma_k}{2\lambda_k^4}\Delta\chi_A(y) \ve_k
                + \frac{\gamma_k}{\lambda^4_k}\nabla\chi_A(y)\cdot  \nabla \ve_k,
                  f'(Q_k) \cdot \ve_k \> + \calo(A\cale_r'),
\end{align*}
where $\cale_r'$ is given by \eqref{caler'-def}.
Taking into account the identity
\begin{align} \label{esti-I2.31}
   &  -\sum_{k=1}^K{\rm Re} \<\frac{\gamma_k}{2\lambda_k^4}\Delta\chi_A(\frac{x-\alpha_k}{\lambda_k}) \ve_k
                + \frac{\gamma_k}{\lambda^4_k}\nabla\chi_A(\frac{x-\alpha_j}{\lambda_k})\cdot  \nabla \ve_k,
                  f'(Q_k) \cdot \ve_k \>  \nonumber \\
  =& \sum_{k=1}^{K}\frac{\gamma_k}{\lambda^4_k}{\rm Re}\int \nabla\chi_A (y)\cdot\nabla \ol{Q_k}
   \ f''(Q_k)\cdot \ve_k^2 dy,
\end{align}
we thus obtain \eqref{I2t},
thereby finishing the proof of Lemma \ref{Lem-I2t}.
\hfill $\square$

We are now in position to prove Theorem \ref{Thm-I-mono}.

{\bf Proof of Theorem \ref{Thm-I-mono}.}
Combining \eqref{I1t-case1} and \eqref{I2t} altogether
and then using
the renormalized variable $\ve_k$ in \eqref{Rj-ej}
we obtain
\begin{align} \label{dI}
\frac{d\scri}{dt}
   \geq&
   \sum\limits_{k=1}^K
   \frac{\gamma_k}{\lambda_k^4}
   \bigg(\int \nabla^2\chi_A (y)(\nabla \varepsilon_k,\nabla \ol{\varepsilon_k}) dy
          + \int |\varepsilon_k|^2dy  \nonumber \\
    &\qquad \qquad  -\int (1+\frac 4d)|Q_k|^\frac{4}{d}\varepsilon_{k,1}^2
           + |Q_k|^{\frac{4}{d}-2} \ol{Q_k}^2 \varepsilon_{k,2}^2  dy
    -\frac{1}{4 }\int \Delta^2\chi_A (y) |\varepsilon_k|^2 dy \nonumber \\
& \qquad \qquad  +   {\rm Re} \int (\na \chi_A (y) -y)\cdot \na \ol{Q_k}
        \ (f''(Q_k) \cdot \ve_k^2) dy \bigg)
    -CA \cale_r,
\end{align}
where $\ve_{k,1}={\rm Re}\ve_k$ and $\ve_{k,2}={\rm Im}\ve_k$, $1\leq k\leq K$.

Then, arguing as in the proof of \cite[(3.83)]{CSZ21}
we obtain that for $A$ large enough,
\begin{align}
   \frac{d \scri}{dt}
   \geq C \sum\limits_{k=1}^K \frac{\g_k}{\lbb_k^4} \int |\na \ve_k|^2 e^{-\frac{|y|}{A}} + |\ve_k|^2 dy
        + \calo\(\sum\limits_{k=1}^K \frac{\g_k}{\lbb_k^4} Scal(\ve_k) + \cale_r\).
\end{align}

Thus, using the inequality
(see \cite[(3.87)]{CSZ21})
\begin{align} \label{Scal-vek}
   \sum\limits_{k=1}^K \  Scal(\ve_k)
   \leq& C \sum\limits_{k=1}^K \( M_k^2 + \|R\|_{L^2}^4 + P^2 \|R\|_{L^2}^2 + e^{-\frac{\delta}{T-t}} \)
\end{align}
we arrive at \eqref{dIt-mono-case2}.
Therefore, the proof is complete.   \hfill $\square$

Below we will fix a large constant $A$  such that Theorem \ref{Thm-I-mono} is valid.

\section{Construction of multi-bubble Bourgain-Wang solutions}  \label{Sec-Const-BW}

In this section we construct the multi-bubble Bourgain-Wang solutions to \eqref{equa-NLS-perturb}
and derive several properties which will be used
in the conditional uniqueness part in Section \ref{Sec-Uniq-BW} later.

Throughout this section, we will take $\ve, \a^*$ sufficiently small
and $t$ close to $T$ such that
\begin{align} \label{a-ve-t-small}
   C \( \( \ve + \a^* \)^{\frac 12} + (1+\max\limits_{1\leq k\leq K} |x_k|) (T-t)^\frac{d}{8+4d} \)  \leq \frac 12,
\end{align}
where $C$ is a universal constant, independent of $\ve, \a^*$
and larger than the constants in the estimates in this section.
Let us mention that,
the exponent ${d}/{(8+4d)}$ is used in the derivation of \eqref{Rn-H23-bdd} below.
For the construction of blow-up solutions,
it will be sufficient to take the exponent $1/4$.

Let us start with the bootstrap estimates of the remainder and geometrical parameters,
which are the key towards the derivation of uniform estimates of solutions.

\subsection{Bootstrap estimates}  \label{Subsec-Boot}

Given any $\upsilon_* \geq 5$,
$m\geq 3$ if $d=2$ and $m\geq 4$ if $d=1$,
set
\begin{align} \label{kappa-def}
   \kappa:= (m+\frac d2 -1)\wedge (\upsilon_*-2).
\end{align}
Note that $\kappa \geq 3$.

\begin{proposition} (Bootstrap estimates) \label{Prop-u-Boot}
Suppose that there exists $t^*\in(0,T_*)$ such that
$u$ admits the unique geometrical decomposition \eqref{u-dec} on $[t^*,T_*]$
and  the following estimates hold:

$(i)$ For the remainder,
\begin{align} \label{D-Tt}
   \|R(t)\|_{L^2}\leq (T-t)^{\kappa+1},\ \ \|\na R(t)\|_{L^2} \leq (T-t)^{\kappa}.
\end{align}
$(ii)$ For the modulation parameters,  $1\leq k\leq K$,
\begin{align}
&|\la_{k}(t) - w_k (T-t) | + |\gamma_{k}(t)  - w_k^2 (T-t) |\leq (T-t)^{\kappa},   \label{lbbn-Tt} \\
&|\al_{k}(t)-x_k|+|\beta_{k}(t)|\leq (T-t)^{\frac{\kappa}{2}+\frac 12},  \label{anbn-Tt}\\
&|\theta_{k}(t) - (w_k^{-2}(T-t)^{-1} + \vartheta_k)| \leq (T-t)^{\kappa-2}. \label{thetan-Tt}
\end{align}

Then,
there exists $t_*\in [0,t^*)$ such that
the decomposition \eqref{u-dec}
and the following improved estimates hold on the larger interval $[t_*, T_*]$:
for $1\leq k\leq K$,
\begin{align}
&\| R(t)\|_{L^2}\leq \frac 12 (T-t)^{\kappa+1},\quad\|\nabla R(t)\|_{L^2}\leq  \frac 12 (T-t)^{\kappa},  \label{D-Tt-boot-2} \\
&\lf|\la_{k}(t) - w_k (T-t) \rt| +  |\gamma_{k}(t) - w_k^2(T-t) |\leq \half (T-t)^{\kappa}, \label{lbbn-Tt-boot-2} \\
&|\al_{k}(t)-x_k|+|\beta_{k}(t)|\leq \frac{1}{2} (T-t)^{\frac{\kappa}{2}+ \frac 12},  \label{anbn-Tt-boot-2} \\
&  |\theta_{k}(t) - (w_k^{-2}(T-t)^{-1} + \vartheta_k) |\leq  \half (T-t)^{\kappa-2}. \label{thetan-Tt-boot-2}
\end{align}
\end{proposition}

\begin{remark}
Since $\kappa \geq 3$,
\begin{align} \label{lbb-g-t}
    \lbb_k, \g_k, P \thickapprox (T-t),\ \ |\beta_k|+|\a_k-x_k|+ D =\calo((T-t)^2),
\end{align}
where the implicit constants are independent of $\ve,\a^*$.
Hence,
the results in the previous Sections \ref{Sec-Gem-Mod} and \ref{Sec-Mass-Energy} are all valid.
Moreover,
since $\kappa = (m+\frac d2-1)\wedge (\upsilon_*-2)$,
we have
\begin{align}  \label{Tt-m-nu-kappa}
   (T-t)^{m+\frac d2} + (T-t)^{\upsilon_*-1} \leq  C (T-t)^{\kappa+1}.
\end{align}
\end{remark}

In order to prove Proposition \ref{Prop-u-Boot},
by the continuity of Jacobian matrix,
the local well-posedness theory of \eqref{equa-NLS-perturb} and $C^1$-regularity of modulation parameters,
we may take $t_*(<t^*)$ sufficiently close to $t^*$,
such that the geometrical decomposition \eqref{u-dec}
and the following estimates hold
on the larger interval  $[t_*, T_*]$:
\begin{align}
&\| R(t)\|_{L^2}\leq 2  (T-t)^{\kappa+1}, \ \ \|\nabla R(t)\|_{L^2}\leq 2 (T-t)^{\kappa}, \label{R-Tt2} \\
&\lf|\la_k(t)- w_k(T-t)\rt| + |\gamma_k(t) - w_k^2(T-t)|\leq 2  (T-t)^{\kappa}, \label{lbb-ga-Tt2}   \\
&|\al_k(t) - x_k| + |\beta_k(t)|\leq 2(T-t)^{\frac {\kappa}{ 2}+\frac 12}, \label{ab-Tt2}  \\
&   |\theta_k(t)- (w_k^{-2}(T-t)^{-1} + \vartheta_k) |\leq 2  (T-t)^{\kappa-2}. \label{theta-Tt2}
\end{align}
By virtue of Theorem \ref{Thm-Mod-bdd}, \ref{Thm-mass-local}, \ref{Thm-Energy} and \ref{Thm-I-mono}
we obtain
\begin{lemma} \label{Lem-Mod-boot}
There exists $C>0$ such that for any $t\in [t_*, T_*]$,
\begin{align}
  & M_k \leq C\a^*  (T-t)^{\kappa+1},  \label{Mk-lbb-boot} \\
  &  Mod  \leq  C\a^*  (T-t)^{\kappa+1},   \label{Mod-lbb-boot}\\
   & |\lbb_k \dot{\lbb}_k + \g_k| \leq C   (T-t)^{\kappa+2},  \label{lbb-lbb-boot}
\end{align}
and for the errors $Er$ and $\cale_r$ in \eqref{Er} and \eqref{calEr-def}, respectively,
\begin{align}
   & |Er| \leq C \a^* (T-t)^{\kappa-1}, \label{Er-lbb-boot} \\
   & |\cale_r| \leq C (\ve+\a^*) (T-t)^{2\kappa-1}+C (T-t)^{2\kappa}. \label{caler-lbb-boot}
\end{align}
\end{lemma}

\begin{remark}
In comparison with \eqref{Mod-lbb-boot},
one more factor $(T-t)$ is gained in \eqref{lbb-lbb-boot} for the
particular modulation equation $\lbb_k \dot{\lbb}_k + \g_k$.
This fact is important to derive \eqref{caler-lbb-boot}
and to close the bootstrap estimates of remainder.
\end{remark}

We are now in position to prove the bootstrap estimates in Proposition \ref{Prop-u-Boot}.

{\bf Proof of Proposition \ref{Prop-u-Boot}}.
{\it (i) Estimate of $R$.}
On one hand,
similarly to \eqref{f''v1-f''v2},
\begin{align}
  |F''(U+z, R) \cdot R^2 - F''(U,R)\cdot R^2|
  \leq C \(|U|^{\frac 4d-1} + |R|^{\frac 4d-1} + |z|^{\frac 4d-1}\) |z| |R|^2 ,
\end{align}
we see that
\begin{align}
   & \bigg| \int F''(U+z, R) \cdot R^2 dx - \int F''(U,R)\cdot R^2 dx \bigg|  \nonumber \\
   \leq& C \int \(|U|^{\frac 4d-1} + |R|^{\frac 4d-1} + |z|^{\frac 4d-1}\)  |z|  |R|^2dx   \nonumber \\
   \leq& C \( (T-t)^{-\frac d2(\frac 4d-1)} \|z\|_{L^\9} \|R\|_{L^2}^2
           +  \|z\|_{L^\9}^{\frac 4d} \|R\|_{L^2}^2
           +  \|z\|_{L^\9} \|R\|_{L^{1+\frac 4d}}^{1+\frac 4d} \) \nonumber \\
   \leq& C \(\a^* (T-t)^{-2+\frac d2} D^2 + \a^* (T-t)^{-d(\frac 2d - \frac 12)}D^{1+\frac 4d} \)  \nonumber \\
   =& o\((T-t)^{-2}D^2\).
\end{align}
Taking into account $F''(U,R)\cdot R^2 = F(U+R) - F(U) - {\rm Re} (f(U)\ol{R})$
we thus get
\begin{align}
   \scri =& \frac 12 \int |\na R|^2 + \frac 12 \sum\limits_{k=1}^K \int \frac{1}{\lbb_k^2} |R|^2 \Phi_k dx
            - {\rm Re} \int F(U+R) - F(U) - f(U) \ol{R} dx \nonumber \\
          & + \sum\limits_{k=1}^K \frac{\g_k}{2\lbb_k} {\rm Im} \int \na \chi_A(\frac{x-\a_k}{\lbb_k}) \cdot \na R \ol{R} \Phi_k dx
          + o\((T-t)^{-2}D^2\).
\end{align}
Then, we  use the expansion
\begin{align*}
  F(U+R)- F(U) - {\rm Re} (f(U) \ol R)
        =& \frac 12 (1+\frac 2d) |U|^{\frac 4d} |R|^2
         + \frac{1}{d} |U|^{\frac 4d -2} {\rm Re} (U^2 \ol{R}^2)
          + \calo\(\sum\limits_{j=3}^{2+\frac 4d} |U|^{2+\frac 4d -j} |R|^j\)
\end{align*}
to derive
\begin{align}  \label{scri-esti}
  \scri = & \frac{1}{2} {\rm Re}  \int |\nabla R|^2
         + \sum_{k=1}^K  \frac{1}{\lambda_{k}^2}|R|^2\Phi_k
         - (1+\frac 2d) |U|^{\frac 4d} |R|^2
         - \frac{2}{d}|U|^{\frac 4d -2} U^2 \ol{R}^2 dx \nonumber \\
  & +\calo\( \sum\limits_{j=3}^{ 2+ \frac 4d} \int |U|^{\frac 4d+2-j} |R|^j dx
         +  \|R\|_{L^2} \|\na R\|_{L^2}\)
         + o\( (T-t)^{-2}D^2 \).
\end{align}
Note that,  the last second line on the R.H.S. above is of order $o((T-t)^{-2}D^2)$, see \cite[(4.29)]{CSZ21},
while for the quadratic terms the following  coercivity type estimate holds  (see \cite[(3.39)]{CSZ21}):
\begin{align} \label{naR2-R2-coer}
   & \frac{1}{2}{\rm Re} \int |\nabla R|^2
         + \sum_{k=1}^K  \frac{1}{\lambda_{k}^2}|R|^2\Phi_k
         - (1+\frac 2d) |U|^{\frac 4d} |R|^2
         - \frac{2}{d}|U|^{\frac 4d -2} U^2 \ol{R}^2 dx  \nonumber \\
   \geq& C\frac{D^2(t)}{(T-t)^2}
        + \calo\(\sum\limits_{k=1}^K \frac{M_k^2}{(T-t)^2}+e^{-\frac{\delta}{T-t}}\),
\end{align}
where $C>0$.
Thus,
by \eqref{scri-esti} and \eqref{naR2-R2-coer},
for $t$ close to $T$,
\begin{align} \label{I-lowbdd}
   \scri \geq& \frac C 2  \frac{D^2}{(T-t)^2}
          - C\(\sum\limits_{k=1}^K \frac{M_k^2}{(T-t)^{2} }
          + e^{-\frac{\delta}{T-t}}\).
\end{align}

On the other hand,
Theorem \ref{Thm-I-mono} yields that for any $t\in[t_*,T_*]$,
\begin{align}  \label{dIdt-lowbdd}
\frac{d \scri}{dt}
\geq   - C\cale_r.
\end{align}

Thus, we infer from \eqref{I-lowbdd}, \eqref{dIdt-lowbdd}
and the boundary condition $\scri(T_*) =0$
that for any $t\in[t_*,T_*]$,
\begin{align*}
    \frac{D^2}{(T-t)^2}
    \leq C\(\int_t^{T_*} |\cale_r| ds + \sum\limits_{k=1}^K \frac{M_k^2}{ (T-t)^2}  + e^{-\frac{\delta}{T-t}}\),
\end{align*}
or, equivalently,
\begin{align} \label{D-caler-Mk}
  D \leq C \( (T-t) \(\int_t^T |\cale_r|ds\)^{\frac 12}
              + \sum\limits_{k=1}^K |M_k| + e^{-\frac{\delta}{T-t}} \).
\end{align}
Taking into account \eqref{Mk-lbb-boot} and \eqref{caler-lbb-boot}
we then obtain
\begin{align}
   D \leq& C \( (\ve+\a^*)^\frac 12 (T-t)^{\kappa+1} + (T-t)^{\kappa+\frac 32} + \a^* (T-t)^{\kappa+1} \) \nonumber \\
     \leq& C \( \ve + \a^* + (T-t) \)^\frac 12 (T-t)^{\kappa+1},
\end{align}
which along with \eqref{a-ve-t-small} yields
\begin{align}
  D\leq \frac 12  (T-t)^{\kappa+1}.
\end{align}
Thus, estimate \eqref{D-Tt-boot-2} is verified.

{\it $(ii)$ Estimates of $\lambda_k$ and $\g_k$.}
By \eqref{Mod-lbb-boot},
\begin{align}
  \bigg|\frac{d}{dt}\(\frac{\gamma_{k}}{\lambda_{k}}\)\bigg|
=\frac{|\lambda_{k}^2\dot{\gamma}_{k}-\lambda_{k}\dot{\lambda}_{k}\gamma_{k}|}{\lambda_{k}^3}
\leq 2 \frac{Mod}{\lambda_{k}^3}
\leq C \a^* (T-t)^{\kappa-2},
\end{align}
which along with the boundary condition  $(\frac {\gamma_{k}}{ \lambda_{k}})(T_*)=w_k$
yields that
\begin{align}  \label{gamlbb-1}
 \bigg| \frac{\gamma_k}{\lambda_k} - w_k  \bigg|
\leq\int_{t}^{T_*}  \bigg| \frac{d}{dr} \(\frac{\gamma_{k}}{\lambda_{k}} \) \bigg|dr
\leq C \a^* (T-t)^{\kappa-1}.
\end{align}
This in turn yields that
\begin{align*}
  \bigg|\frac{d}{dt}(\lambda_{k} - w_k (T-t))\bigg|
 =\bigg|\dot{\lambda}_{k}+\frac{\gamma_{k}}{\lambda_{k}}+ w_k-\frac{\gamma_{k}}{\lambda_{k}}\bigg|
\leq\frac{Mod}{\lambda_{k}}+ C \a^* (T-t)^{\kappa-1}
\leq C \a^* (T-t)^{\kappa-1},
\end{align*}
and thus, by \eqref{a-ve-t-small},
\begin{align} \label{lbb-Tt*}
  |\lambda_{k} - w_k(T-t) |
\leq& \int_{t}^{T_*}  \bigg|\frac{d}{dr}(\lambda_{k}- w_k(T-r)) \bigg|dr \nonumber \\
\leq& C \a^*(T-t)^\kappa
\leq  \frac{1}{2} (T-t)^{\kappa}.
\end{align}
Hence, we prove the estimate of $\lbb_k$ in \eqref{lbbn-Tt-boot-2}.

Regarding $\g_k$,
by \eqref{Mod-lbb-boot} and \eqref{gamlbb-1},
\begin{align*}
   \bigg|\frac{d}{dt}(\g_k - w_k^2 (T-t))\bigg|
   =  \bigg|\dot\g_{k} + \frac{\g_k^2}{\lbb_k^2} + w_k^2 - \frac{\g_k^2}{\lbb_k^2} \bigg|
   \leq  \frac{Mod}{\lbb_k^2} + C\bigg|w_k - \frac{\g_k}{\lbb_k}\bigg|
   \leq C \a^* (T-t)^{\kappa-1}.
\end{align*}
Thus, taking into account $\g_k(T_*) = \omega^2_k (T-T_*)$ and \eqref{a-ve-t-small}  we get
\begin{align}  \label{g-Tt*}
   |\g_k(t) - w_k^2 (T-t)|
   \leq \int_t^{T_*}  \bigg|\frac{d}{dr} (\g_k(r) - w_k^2 (T-r)) \bigg| dr
   \leq C \a^* (T-t)^{\kappa}
   \leq \frac 12 (T-t)^{\kappa}.
\end{align}
This gives the estimate of $\g_k$ in \eqref{lbbn-Tt-boot-2}.

{\it $(iii)$ Estimates of $\beta_k$ and $\alpha_k$.}
By the improved estimate \eqref{energy-esti}, \eqref{Er-lbb-boot} and \eqref{gamlbb-1},
\begin{align} \label{b-Tt-a}
   \frac{|\beta_k|^2}{\lbb_k^2}
   \leq C \( \sum\limits_{k=1}^K \bigg|w_k - \frac{\g_k}{\lbb_k} \bigg|   + Er\)
   \leq C \a^*(T-t)^{\kappa-1},
\end{align}
which along with \eqref{a-ve-t-small} yields that
\begin{align} \label{b-Tt*}
   |\beta_k| \leq C (\a^*)^\frac 12 (T-t)^{\frac \kappa2 + \frac 12}
   \leq \frac 12 (T-t)^{\frac \kappa2 + \frac 12}.
\end{align}

Moreover, by \eqref{Mod-lbb-boot} and \eqref{b-Tt*},
\begin{align}
|\dot{\alpha}_{k}|=
\bigg|\frac{\lambda_k\dot{\alpha}_{k}-2\beta_{k}}{\lambda_{k}}+\frac{2\beta_{k}}{\lambda_{k}}\bigg|
\leq\frac{Mod}{\lambda_{k}}+\frac{2|\beta_{k}|}{\lambda_{k}}
\leq C \a^* (T-t)^{\frac{\kappa}{2}-\frac 12}.
\end{align}
Integrating both sides and using \eqref{a-ve-t-small}
and the boundary condition $\a_k(T_*) = x_k$ we get
\begin{align}
|\alpha_{k}(t)-x_k|\leq\int_{t}^{T_*}|\dot{\alpha}_{k}(r)|dr
\leq C \a^* (T-t)^{\frac \kappa 2 + \frac 12}
\leq \frac{1}{2} (T-t)^{\frac{\kappa}{2}+ \frac 12},
\end{align}
thereby proving the estimate of $\a_k$ in \eqref{anbn-Tt-boot-2}.

{\it $(iv)$ Estimate of $\theta_k$.}
It remains to estimate $\theta_k$.
By \eqref{lbb-g-t}, \eqref{Mod-lbb-boot}, \eqref{lbb-Tt*} and \eqref{b-Tt-a},
\begin{align}  \label{esti-thetan.0*}
   \bigg|\frac{d}{dt}(\theta_{k} - {w_k^{-2}(T-t)^{-1}}+ \vartheta_k)\bigg|
=& \bigg|\frac{\lbb_{k}^2 \dot{\theta}_{k} -1 -|\beta_{k}|^2}{\lbb_{k}^2}
  +\frac{|\beta_{k}|^2}{\lbb_{k}^2}
  + \frac{1}{\lbb_{k}^2}
  - \frac{1}{w_k^2(T-t)^2} \bigg| \nonumber \\
\leq& \frac{Mod}{\lbb_k^2} + \frac{|\beta_k|^2}{\lbb_k^2} +
     \frac{|\lbb_{k}-w_k(T-t)||\lbb_{k}+w_k(T-t)|}{w_k^2\lbb_{k}^2(T-t)^2 } \nonumber \\
   \leq& C \a^* (T-t)^{\kappa-3},
\end{align}
which along with \eqref{a-ve-t-small}
and the boundary $\theta_k(T_*) = w_k^{-2} (T-T_*)^{-1} + \vartheta_k$ yields that
\begin{align}
   |\theta_{k}- ({w_k^{-2}(T-t)^{-1}} + \vartheta_k) |
 \leq& \int_{t}^{T_*} \bigg| \frac{d}{dr}(\theta -  {w_k^{-2}(T-r)^{-1}} +\vartheta_k) \bigg|dr \nonumber \\
 \leq& C \a^* (T-t)^{ \kappa-2}
\leq \frac{1}{2} (T-t)^{\kappa-2}.
\end{align}
Hence, the estimate \eqref{thetan-Tt-boot-2} is verified.
Therefore, the proof of Proposition \ref{Prop-u-Boot} is complete.
\hfill $\square$

\subsection{Proof of existence}   \label{Subsec-Proof-Const-BW}

We are now in position to prove the existence part in Theorem \ref{Thm-BW-RNLS}.
Consider the approximating solutions $v_n$
satisfying the equation
\be    \label{equa-u-t}
\left\{ \begin{aligned}
 &i\partial_t v_n+\Delta v_n +  a_1 \cdot \nabla  v_n +a_0 v_n  +|v_n|^{\frac 4d}v_n =0,   \\
 &v_n(t_n)=\sum_{k=1}^{K}S_k(t_n)+z(t_n),
\end{aligned}\right.
\ee
where $\{t_n\}$ is any increasing sequence converging to $T$,
the coefficients $a_1,a_0$ are given by \eqref{a1-loworder} and \eqref{a0-loworder}, respectively,
$\{S_k\}$
are the pseudo-conformal blow-up solutions defined in \eqref{Sj-blowup},
and $z$ solves equation \eqref{equa-z}.

As a consequence of bootstrap estimates,
we have the key uniform estimates below.

\begin{lemma} (Uniform estimates)    \label{Lem-u-Unibdd}
There exists $t_*\in [0,T)$ such that for $n$ large enough,
$v_n$ admits the unique
geometrical decomposition $v_n=U_n+ z+R_n$ as in \eqref{u-dec},
with the parameters $\calp_{n,k}:=(\lbb_{n,k}, \a_{n,k}, \beta_{n,k}, \g_{n,k}, \theta_{n,k})$,
$1\leq k\leq K$,
and the estimates \eqref{D-Tt}-\eqref{thetan-Tt} hold on $[t_*,t_n]$.
Moreover,
there exists $C>0$ such that
\begin{align}\label{R-Sigma}
  \sup\limits_n \|R_n(t)\|_{\Sigma}\leq C(T-t)^{\kappa},
\end{align}
and \begin{align}\label{u-sigma}
  \sup\limits_{n} \|xv_n\|_{C([t_*,t_n];L^2)}\leq C (1+\max\limits_{1\leq k\leq K}|x_k|)^2,
\end{align}
\end{lemma}

{\bf Proof.}
The proof of the existence of a universal time $t_*$ and
uniform estimates \eqref{D-Tt}-\eqref{thetan-Tt} is similar to that of \cite[Theorem 5.1]{SZ19},
mainly based on the bootstrap estimates in Proposition \ref{Prop-u-Boot}
and bootstrap arguments (see, e.g., \cite[Proposition 1.21]{T06}).
Thus, the details are omitted here for simplicity.
Below let us mainly prove estimates \eqref{R-Sigma} and \eqref{u-sigma}.

Let $M:=1+\max_{1\leq k\leq K}|x_k|$.
Let $\varphi(x)\in C^1(\R^d,\R)$ be a radial cutoff function such that
$\varphi(x)=0$ for $|x|\leq r$,
and $\varphi(x)=(|x|-r)^2$ for $|x|>r$,
where $r=2\max_{1\leq k\leq K}\{|x_k|,1\}$.
Note that, $|\nabla \varphi|\leq C\varphi^{\half}$
for a universal constant $C>0$.

Let $w_n:= U_n + R_n$, $n\geq 1$.
Then, $v_n = w_n +z$.
By equations \eqref{equa-u-t} and \eqref{equa-z},
$w_n$ solves equation
\be     \label{equa-wn}
\left\{ \begin{aligned}
   & i\partial_t w_n+\Delta w_n + a_1 \cdot \nabla w_n+ a_0 w_n+f(v_n)-f(z)=0, \\
   & w_n(t_n)  = \sum\limits_{k=1}^K S_k(t_n)\ (=: S(t_n)).
\end{aligned}\right.
\ee
Then, by the integration by parts formula
and ${\rm Im}\ w_n \ol{f(w_n)} =0$,
\begin{align} \label{dt-wnL2}
\frac{d}{dt}\int |w_n|^2\varphi dx
 &={\rm Im}\int (2\ol{w_n} \nabla w_n+ a_1|w_n|^2)\cdot  \nabla\varphi
     + 2 w_n( \ol{f(v_n)}- \ol{f(w_n)} - \ol{f(z)}) \varphi dx \nonumber\\
 &={\rm Im}\int (2\ol{w_n} \nabla w_n+ a_1|w_n|^2)\cdot  \nabla\varphi dx
   + \calo\(\sum_{j=1}^{4/d}\int |w_n^{2+\frac 4d-j}z^{j}\varphi|dx\).
\end{align}

In order to estimate the R.H.S. of \eqref{dt-wnL2},
we note that
\begin{align}   \label{1}
 & \bigg|{\rm Im}\int (2\ol{w_n} \nabla w_n+ a_1|w_n|^2)\cdot  \nabla\varphi dx \bigg| \nonumber \\
\leq& C\int_{ |x-x_k|\geq 1,1\leq k\leq K}\(|w_n||\nabla w_n|+|w_n|^2\)\varphi^{\half} dx  \nonumber  \\
\leq& C\(\(\int_{ |x-x_k|\geq 1,1\leq k\leq K}|\nabla w_n|^2dx\)^{\half}+\(\int_{ |x-x_k|\geq 1,1\leq k\leq K}|w_n|^2dx\)^{\half}\)
   \(\int|w_n|^2\varphi dx\)^\half,
\end{align}
where  $C>0$ is  independent of $n$.
By  \eqref{u-dec},  \eqref{D-Tt}  and \eqref{Q-decay},
\begin{align}
    \bigg|\int_{ |x-x_k|\geq 1,1\leq k\leq K} |w_n(t)|^2+ |\nabla w_n(t)|^2dx\bigg|
\leq C(\|R_n(t)\|_{H^1}^2  +e^{-\frac{\delta}{T-t}})
\leq C (T-t)^{2\kappa}.
\end{align}
This yields that for a universal constant $C>0$,
\begin{align} \label{dt-wnL2-esti1}
  \bigg| {\rm Im}\int (2\ol{w_n} \nabla w_n+ a_1|w_n|^2)\cdot  \nabla\varphi dx \bigg|
  \leq C(T-t)^\kappa \(\int |w_n|^2 \vf dx \)^\frac 12.
\end{align}
Moreover, for $1\leq j\leq \frac 4d$,
since ${\rm supp} \vf \subseteq \{x: |x-x_k|\geq 1, 1\leq k\leq K\}$,
\be\ba
\int |w_n^{2+\frac 4d-j} z^{j}\varphi|dx
&\leq \(\int |w_n|^2\varphi dx\)^{\half}
\(\int |w_n|^{2+\frac 8d-2j}|z|^{2j}\varphi dx\)^\half\nonumber\\
&\leq C\(\int |w_n|^2\varphi dx\)^{\half}\(\int \(|U_n|^{2+\frac 8d-2j}+|R_n|^{2+\frac 8d-2j}\)|z|^{2j}\varphi dx\)^\half\nonumber\\
&\leq C\(\int |w_n|^2\varphi dx\)^{\half}\(\int |R_n|^{2+\frac 8d-2j}|z|^{2j}\varphi dx+ M^2 e^{-\frac{\delta}{T-t}}\)^\half\nonumber\\
&\leq C\(\int |w_n|^2\varphi dx\)^{\half}\(\|R_n\|_{L^{2(2+\frac 8d -2j)}}^{2+\frac 8d -2j} \|xz\|^2_{L^4}\|z\|_{L^\9}^{2j-2}
          + M^2 e^{-\frac{\delta}{T-t}}\),
\ea\ee
which along with \eqref{z-LtHm-a*}, \eqref{xz-H1-small}, \eqref{a-ve-t-small} and \eqref{D-Tt-boot-2}
yields that for a universal constant $C>0$,
\begin{align} \label{dt-wnL2-esti2}
\int |w_n^{2+\frac 4d-j} z^{j}\varphi|dx\leq C(T-t)^{2\kappa}\(\int |w_n|^2\varphi dx\)^{\half}.
\end{align}
Hence, plugging \eqref{dt-wnL2-esti1} and \eqref{dt-wnL2-esti2} into \eqref{dt-wnL2}
we get
\begin{align} \label{dt-wnL2-esti}
       \bigg|\frac{d}{dt}\int |w_n(t)|^2\varphi dx\bigg|
 \leq& C(T-t)^{\kappa}\(\int|w_n(t)|^2\varphi dx\)^\half.
\end{align}

Thus,
integrating \eqref{dt-wnL2-esti} from $t$ to $t_n$,
using \eqref{a-ve-t-small} and the boundary estimate
\begin{align}
  \int |w_n(t_n)|^2\varphi dx
  = \int \bigg|\sum\limits_{k=1}^KS_k(t_n) \bigg|^2\varphi dx
  \leq C M^2 e^{-\frac{\delta}{T-t_n}}
  \leq C M^2 e^{-\frac{\delta}{T-t}}
\end{align}
we obtain for $t\in[0,t_n]$,
\begin{align} \label{wL2-Tt}
\int|w_n(t)|^2\varphi dx
\leq C (T-t)^{2\kappa+2}.
\end{align}

In particular, this yields that
\begin{align}   \label{Rn-H2vf-esti}
   \int|R_n(t)|^2\varphi dx
   \leq& C \(\int|U_n(t)|^2\varphi dx+\int|w_n(t)|^2\varphi dx\)
    \leq C (T-t)^{2\kappa+2} .
\end{align}
Since  $\varphi(x)\geq \frac 14 |x|^2$ for $|x| \geq 4M $,
by \eqref{a-ve-t-small}, \eqref{D-Tt-boot-2} and \eqref{Rn-H2vf-esti},
\begin{align}
\int|xR_n(t)|^2 dx\leq C\(\int|R_n(t)|^2\varphi dx+ M^2 \int|R_n(t)|^2 dx\)
\leq C M^2(T-t)^{2\kappa+2}
\leq C (T-t)^{2\kappa+1},
\end{align}
where $C$ is independent of $n$ and $M$.
This along with \eqref{D-Tt} yields \eqref{R-Sigma}.

Similarly, we  derive that
\begin{align}
\int|x v_n(t)|^2 dx
  \leq& C \(\int |xw_n|^2 dx + \|xz\|_{L^2} \)^2 \nonumber \\
  \leq& C\(\int|w_n(t)|^2\varphi dx+ M^2 \|w_n\|^2_{L^2} + \|xz\|_{L^2}^2 \)  \nonumber \\
  \leq& C \((T-t)^{2\kappa+2}  + K M^2 \|Q\|_{L^2}^2 + M^2(T-t)^{2\kappa+2} + \|xz\|_{L^2}^2 \)
  \leq C M^2,
\end{align}
where the last step is due to \eqref{xz-H1-small}, \eqref{D-Tt}, \eqref{wL2-Tt}
and the conservation law of mass,
and $C$ is independent of $n$.
This yields \eqref{u-sigma}.
Therefore, the proof of Lemma \ref{Lem-u-Unibdd} is complete.
\hfill $\square$

{\bf Proof of existence part in Theorem \ref{Thm-BW-RNLS}.}
Let $\a^*, \ve$ be small enough such that \eqref{a-ve-t-small} holds
and  $t_n$ as in Lemma \ref{Lem-u-Unibdd}.
Let $M:=1+\max_{1\leq k\leq K}|x_k|$.
By Lemma \ref{Lem-u-Unibdd},
$\{v_n(t_*)\}$ are uniformly bounded in $\Sigma$,
and thus up to a subsequence (still denoted by $\{n\}$),
$v_n(t_*)$ converges weakly to some $v_* \in \Sigma$.
The weak convergence indeed can be enhanced to the strong one in the space $L^2$, i.e.,
\begin{align} \label{unt0-u0-L2}
    v_n(t_*) \to v_*,\ \  in\ L^2,\ as\ n \to \9.
\end{align}
This is due to the uniform integrability of $\{v_n(t_*)\}$
implied by the uniform estimate \eqref{u-sigma}:
\begin{align}  \label{un-uninteg-L2}
 \sup\limits_{n\geq 1} \|v_n(t_*)\|_{L^2(|x|>A)}
  \leq \frac{1}{A} \sup\limits_{n\geq 1}\|x v_n(t_*)\|_{L^2(|x|>A)}
  \leq \frac{C M^2}{A} \to 0, \ \ as\ A\to \9.
\end{align}

Thus,
the $L^2$ local well-posedness theory (see, e.g. \cite{BRZ14})
yields a unique $L^2$-solution $v_c$ to \eqref{equa-NLS-perturb} on $[t_*, T)$,
satisfying that $v_c(t_*)=v_*$
and
\begin{align} \label{un-u-0-L2}
\lim_{n\rightarrow \infty}\|v_n-v_c\|_{C([t_*,t];L^2)}=0,\ \ t\in [t_*, T).
\end{align}
Moreover,
since $v_* \in \Sigma$,
the local well-posedness result also yields
$v_c\in C([t_*,t]; \Sigma)$ for  $t\in(t_*,T)$.

Next, we show that $v_c$ is the desired multi-bubble Bourgain-Wang solution to \eqref{equa-NLS-perturb}.

As a matter of fact,
let
\begin{align*}
   & (\lbb_{0,k}, \a_{0,k}, \beta_{0,k}, \g_{0,k}, \theta_{0,k})
           : =  (w_k(T-t), x_k, 0, w_k^2(T-t), w_k^{-2} (T-t)^{-1} +\vartheta_k)
\end{align*}
and $\calp_{n,k}:= (\la_{n,k},\alpha_{n,k},\beta_{n,k},\gamma_{n,k},\theta_{n,k})$ be the parameters
corresponding to the geometrical decomposition of $v_n$.
Then, analogous computations as in \cite{SZ19}
and estimates \eqref{D-Tt}-\eqref{thetan-Tt}  yield, for $\kappa \geq 1$,
\begin{align}  \label{Un-S-L2-Tt}
   \|U_n - S\|_{L^2}
   \leq& C \sum\limits_{j=1}^K
        \bigg( \bigg| \frac{\lbb_{0,k}}{\lbb_{n,k}} -1\bigg|
        + \bigg|\frac{\a_{n,k}-\a_{0,k}}{\lbb_{n,k}}\bigg|
        + \bigg|\beta_{n,k} -\beta_{0,k}\bigg|
        +  \bigg| \g_{n,k} - \g_{0,k}\bigg| \nonumber  \\
       & \qquad  \ \
         + \bigg|\frac{\lbb_{n,k}^\frac d2 - \lbb_{0,k}^\frac d2}{\lbb_{0,k}^\frac d2}\bigg|
         + \bigg|\theta_{n,k} - \theta_{0,k}\bigg|  \bigg)   \nonumber \\
   \leq& C \( (T-t)^{\frac 12 \kappa -\frac 12} +   (T-t)^{\kappa-2} \) \nonumber \\
   \leq& C (T-t)^{\frac 12 (\kappa -1)},
\end{align}
which along with \eqref{un-u-0-L2} yields that
\begin{align} \label{v-S-z-L2-esti}
    \|v_c(t) - S(t) -z(t)\|_{L^2}
    \leq \lim\limits_{n\to \9} \(\|U_n(t) - S(t)\|_{L^2} + \|R_n(t)\|_{L^2}\)
    \leq C (T-t)^{\frac 12 (\kappa -1)}.
\end{align}

Moreover, as in \cite{SZ20},
\begin{align}   \label{Un-S-H1-Tt}
    \|U_n - S\|_{\Sigma}
   \leq& C M \sum\limits_{k=1}^K
        \bigg(  \frac{1}{\lbb_{0,k}} \bigg| \frac{\lbb_{0,k}}{\lbb_{n,k}} -1\bigg|
             + \bigg|\frac{\a_{n,k}-\a_{0,k}}{\lbb_{0,k}\lbb_{n,k}}\bigg|
        + \bigg| \frac{\beta_{n,k} -\beta_{0,k}}{\lbb_{0,k}} \bigg|
         +  \bigg|\frac{\g_{n,k} - \g_{0,k}}{\lbb_{0,k}} \bigg|  \nonumber  \\
       & \qquad   + \bigg|\frac{\lbb_{n,k}^{1+\frac d2} - \lbb_{0,k}^{1+\frac d2}}{\lbb_{n,k} \lbb_{0,k}^{1+\frac d2}} \bigg|
           + \bigg|\frac{\theta_{n,k} - \theta_{0,k}}{\lbb_{0,k}} \bigg|   \bigg)  \nonumber  \\
   \leq& C M (T-t)^{\frac 1 2(\kappa -3)},
\end{align}
which, via \eqref{R-Sigma}, yields that
\begin{align}
   \|v_n(t) - S(t) - z(t)\|_{\Sigma}
   \leq \|U_n(t) - S(t)\|_{\Sigma} + \|R_n(t)\|_{\Sigma}
   \leq C M (T-t)^{\frac{1}{2} (\kappa -3)}.
\end{align}
Hence,
possibly selecting a further subsequence (still denoted by $\{n\}$)
and using \eqref{un-u-0-L2}
we obtain
\begin{align*}
   v_n(t) -S(t) -z(t)
   \rightharpoonup v(t) - S(t) -z(t),\ \ weakly\ in\ \Sigma,\ as\ n\to \9,
\end{align*}
which yields that
\begin{align*}
    \|v_c(t) - S(t) -z(t) \|_{\Sigma}
    \leq \liminf\limits_{n\to \9}  \|v_n(t) - S_T(t) - z(t) \|_{\Sigma}
    \leq C M (T-t)^{\frac{1}{2}(\kappa-3)}.
\end{align*}

Therefore, the proof of existence part in Theorem \ref{Thm-BW-RNLS} is complete.
\hfill $\square$

\subsection{Further properties}   \label{Subsec-Further}

We close this section with  further properties of
the constructed multi-bubble Bourgain-Wang solutions in Theorem \ref{Thm-BW-RNLS},
which will be used in Section \ref{Sec-Uniq-BW} later.

\begin{proposition} ($H^\frac 32$ boundedness) \label{Prop-Rn-H23-bdd}
Consider the situations as in Theorem \ref{Thm-BW-RNLS}.
Then,
\begin{align} \label{Rn-H23-bdd}
\|R_n(t)\|_{{H}^\frac{3}{2}}\leq (T-t)^{\kappa-2}, \ \  t\in[t_*,t_n),
\end{align}
where $\kappa:= (m+\frac d2-1)\wedge (\upsilon_*-2)$.
\end{proposition}

{\bf Proof.}
Set $M:= 1+\max_{1\leq j\leq K}|x_j|$.
Rewrite equation \eqref{equa-R}:
\begin{align}   \label{equa-Rn}
  & i\partial_t R_n +\Delta R_n+ (a_1 \cdot \nabla +a_0)  R_n
    =  -\eta_n -f(R_n) -(f(v_n)-f(U_n+z)-f(R_n)),
\end{align}
where $ R_n(t_n)=0$ and $\eta_n$ is given by \eqref{etan-Rn}.
Applying  $\<\na\>^{\frac 32}$ to both sides of \eqref{equa-Rn}
yields
\begin{align}   \label{equa-na32Rn}
  & i\partial_t (\<\na\>^\frac 32 R_n) +\Delta (\<\na\>^\frac 32 R_n) +(a_1 \cdot \nabla +a_0) (\<\na\>^\frac 32 R_n)  \nonumber \\
   =& [a_1\cdot \na + a_0, \<\na\>^{\frac 32}] R_n
     -\<\na\>^\frac 32 \eta_n
     -\<\na\>^\frac 32f(R_n)
    -\<\na\>^\frac 32 \(f(v_n)-f(U_n+z)-f(R_n)\),
\end{align}
where $[a_1\cdot \na + a_0, \<\na\>^{\frac 32}]$ is the commutator
$(a_1\cdot \na + a_0)\<\na\>^{\frac 32} - \<\na\>^{\frac 32} (a_1\cdot \na + a_0)$.
Then, the Strichartz  and local smoothing estimates yield
\begin{align}\label{Rn-Est-32}
      \|R_n\|_{C([t,t_n]; {H}^\frac{3}{2})}
  \leq& C \bigg(\|[a_1\cdot \na +a_0,\<\na\>^\frac 32] R_n\|_{L^2(t,t_n; H^{-\frac 12}_1)}
           +\|\<\nabla\>^{\frac{3}{2}}(f(R_n))\|_{L^{\frac{4+2d}{4+d}}(t,t_n;L^{\frac{4+2d}{4+d}})} \nonumber \\
      &\ \ \ \  + \|\<\nabla\>^{\frac{3}{2}}\eta_n\|_{L^{\frac{4+2d}{4+d}}(t,t_n;L^{\frac{4+2d}{4+d}})}
          +\|\<\na\>^{\frac 32}(f(u_n)-f(U_n+z)-f(R_n))\|_{L^2(t,t_n;H_1^{-\frac 12})} \bigg) \nonumber \\
     =&: \sum\limits_{l=1}^4 J_l.
\end{align}

To estimate the R.H.S. above,
by the calculus of pseudo-differential operators,  \eqref{decay}
and \eqref{D-Tt-boot-2},
\begin{align} \label{na-bnac-commu}
   J_1 \leq C \|R_n\|_{L^2(t,t_n; H^{1}_{-1})}
   \leq  C (T-t)^\frac 12 \|R_n\|_{C([t,t_n];H^1)}
   \leq C (T-t)^{\kappa+\frac 12}.
\end{align}

Moreover, similarly to \cite[(7.8)]{SZ20},
by the product rule,
Sobolev's embedding and \eqref{D-Tt},
\begin{align}
      \|\<\nabla\>^{\frac{3}{2}}(f(R_n))\|_{L^{\frac{4+2d}{4+d}}}
\leq  C \|R_n\|_{H^1}^{\frac 4d} \|R_n\|_{H^{\frac 32}}
 \leq C(T-t)^{\frac 4d \kappa}  \|R_n\|_{H^{\frac 32}},
\end{align}
which yields that
\begin{align} \label{Rn-Est-32-2}
  J_2 \leq C(T-t)^{\frac{4}{d} \kappa + \frac{4+d}{4+2d}}\|R_n\|_{C([t,t_n]; {H}^\frac{3}{2})}.
\end{align}

Regarding $J_3$, we use the decomposition $\eta_n = \sum_{l=1}^4 \eta_l$ as in \eqref{eta1-eta4}
to derive that for
$p:= \frac{4+2d}{4+d}$
and any multi-index $|\upsilon|\leq 2$,
by \eqref{eta1-def} and \eqref{Mod-lbb-boot},
\begin{align}  \label{pxeta1}
  \|\partial_x^\upsilon \eta_1\|_{L^p(t,t_n; L^p)}
  \leq C \sum\limits_{k=1}^K \lbb_k^\frac 1p \lbb_k^{-2-|\upsilon|+d(\frac 1p - \frac 12)} Mod
  \leq C \a^* (T-t)^{\kappa-2+\frac{d}{4+2d}}.
\end{align}
Moreover, by \eqref{eta2-def},
\begin{align}
    \|\partial_x^\upsilon \eta_2 \|_{L^p}
    \leq C (T-t)^{-4-\frac d2 + \frac dp} \sum\limits_{|\upsilon|\leq 2} \sum\limits_{k=1}^K \|e^{-\delta|y|} \partial_y^\upsilon \ve_{z,k}\|_{L^\9}
        + C e^{-\frac{\delta}{T-t}},
\end{align}
which along with \eqref{eyvez-bdd} and \eqref{Tt-m-nu-kappa}  yields
\begin{align}  \label{pxeta3}
   \|\partial_x^\upsilon \eta_2 \|_{L^p(t,t_n; L^p)}
   \leq& C (T-t)^{\frac 1p - 4 + \frac{d}{2+d}}
        \sum\limits_{|\upsilon|\leq 2} \sum\limits_{k=1}^K
        \|e^{-\delta|y|} \partial_y^\upsilon \ve_{z,k}\|_{L^\9}
        + C e^{-\frac{\delta}{T-t}} \nonumber \\
   \leq& C \a^* (T-t)^{m-2+\frac 2d+\frac{d}{4+2d}}
   \leq C(T-t)^{\kappa-1}.
\end{align}
Note that,
because $\eta_3$ contains the interactions between different blow-up profiles,
by  Lemma \ref{Lem-decoup-U},
\begin{align}  \label{pxeta4}
   \|\partial_x^\upsilon \eta_3 \|_{L^p(t,t_n; L^p)}
   \leq C e^{-\frac{\delta}{T-t}}.
\end{align}
At last, by \eqref{eta4-def},
\begin{align*}
   \partial_x^\upsilon \eta_4 = \sum\limits_{k=1}^K
   \lbb_k^{-\frac d2 - |\upsilon|} \partial_y^\upsilon
   \(\wt a_{1,k} \lbb_k^{-1}\cdot \na Q_k + \wt a_{0,k} Q_k\)(\frac{x-\a_k}{\lbb_k}),
\end{align*}
which by \eqref{a1-Pnu}, \eqref{a0-Pnu} and \eqref{Tt-m-nu-kappa} yields that
\begin{align}  \label{pxeta2}
   \|\partial_x^\upsilon \eta_4\|_{L^p(t,t_n; L^p)}
    \leq& C (T-t)^{\upsilon_*-2+\frac{d}{4+2d}}
   \leq C (T-t)^{\kappa+ \frac{d}{4+2d}}.
\end{align}

Hence, we conclude that
\begin{align}\label{Rn-Est-32-1}
 J_3
 \leq \| \eta_n\|_{L^{\frac{4+2d}{4+d}}(t,t_n;H^{2,\frac{4+2d}{4+d}})}
 \leq C(T-t)^{\kappa-2+\frac{d}{4+2d}}.
\end{align}

It remains to estimate the last term $J_4$.
We estimate
\begin{align}  \label{x-fun-fUn-fRn}
 J_4 \leq & C \|\langle x\rangle(f(v_n)-f(U_n+z)-f(R_n))\|_{L^2(t,t_n;H^{1})} \nonumber \\
\leq&  C \sum\limits_{j=1}^{4/d}
       \bigg( \|\<x\> (|U_n|^{1+\frac 4d- j}+|z|^{1+\frac 4d- j}) |R_n|^j\|_{L^2(t,t_n; L^2)}
     + \|\<x\> (|\na U_n|+|\na z|)(|U_n |^{\frac 4d- j}+|z |^{\frac 4d- j}) |R_n|^j\|_{L^2(t,t_n; L^2)}  \nonumber \\
     & \qquad \ \  + \|\<x\>(|U_n|^{1+\frac 4d- j} +|z|^{1+\frac 4d- j}) |\na R_n|  |R_n|^{j-1}\|_{L^2(t,t_n; L^2)}  \bigg).
\end{align}
Note that, by \eqref{z-LtHm-a*} and \eqref{xz-H1-small},
\begin{align*}
    \|\<x\> |z|^{1+\frac 4d- j} |R_n|^j\|_{L^2(t,t_n; L^2)}
    \leq  (T-t)^\frac 12 \|\<x\>z\|_{L^2} \|z\|_{L^\9}^{\frac 4d-j}
         \|R\|_{C([t,t_n];H^1)}^j
    \leq C (T-t)^{\kappa j+\frac 12}.
\end{align*}
Since $\|\<x\> \na U_n\|_{L^\9} \leq CM (T-t)^{-\frac d2 -1}$
and $\|\na U_n\|_{L^\9} \leq  C (T-t)^{-\frac d2}$,
by \eqref{xz-H1-small} and \eqref{D-Tt},
\begin{align*}
   & \|\<x\>  \(|\na U_n| |z |^{\frac 4d- j} + |\na z| | U_n|^{\frac 4d- j} +  |\na z| |  z|^{\frac 4d- j}  \) |R_n|^j\|_{L^2(t,t_n; L^2)}  \nonumber \\
   \leq& C\bigg( M (T-t)^{-\frac d2 -\frac 12} \|z\|_{L^\9(t,t_n;L^\9)}^{\frac 4d-j} \|R\|_{C([t,t_n];H^1)}^j
         +  (T-t)^{-\frac d2(\frac 4d-j)+\frac 12} \|\<x\> \na z\|_{{L^\9}(t,t_n;H^1)}\|R\|_{C([t,t_n];H^1)}^j \nonumber \\
       &\qquad +   (T-t)^\frac 12 \|\<x\> \na z\|_{{L^\9}(t,t_n;H^1)} \|z\|_{L^\9(t,t_n;L^\9)}^{\frac 4d-j} \|R\|_{C([t,t_n];H^1)}^j \bigg)  \nonumber \\
   \leq& C \(M(T-t)^{-\frac d2 - \frac 12 +\kappa j} + (T-t)^{-2+\frac d2j + \frac 12 +\kappa j}
             + (T-t)^{\frac 12+ \kappa j} \) \nonumber \\
   \leq&  C M (T-t)^{\kappa-\frac 32}.
\end{align*}
Moreover,
\begin{align*}
   \|\<x\> |z|^{1+\frac 4d- j}  |\na R_n|  |R_n|^{j-1}\|_{L^2(t,t_n; L^2)}
   \leq& C(T-t)^{\frac 12}  \|\<x\>z\|_{{L^\9}(t,t_n;H^1)} \|z\|_{{L^\9}(t,t_n;L^\9)}^{\frac 4d-j}  \|R\|_{C([t,t_n];H^1)}^j \nonumber \\
    \leq& C (T-t)^{\frac 12 +\kappa}.
\end{align*}
The remaining terms in \eqref{x-fun-fUn-fRn} only  involves $U_n$ and $R_n$ and can be bounded by, as in \cite{SZ20},
\begin{align*}
   C M\((T-t)^{\kappa-\frac 32} + (T-t)^{\kappa-\frac 32} \|R\|_{C([t,t_n]; H^\frac 32)} \).
\end{align*}
Thus, we conclude that
\begin{align} \label{Rn-Est-32-3}
  J_4 \leq  C M \((T-t)^{\kappa-\frac 32} + (T-t)^{\kappa-\frac 32} \|R_n\|_{C([t,t_n]; {H}^\frac{3}{2})}\).
\end{align}

Therefore, estimates \eqref{Rn-Est-32},
\eqref{na-bnac-commu}, \eqref{Rn-Est-32-2}, \eqref{Rn-Est-32-1} and \eqref{Rn-Est-32-3}
altogether yield that
\begin{align} \label{Rn-CtH23}
   \|R_n\|_{C([t,t_n]; {H}^\frac{3}{2})}
   \leq C M \((T-t)^{\kappa- 2+ \frac{d}{4+2d}} + (T-t)^{\kappa-\frac 32} \|R_n\|_{C([t,t_n]; {H}^\frac{3}{2})} \),
\end{align}
which along with \eqref{a-ve-t-small} yields \eqref{Rn-H23-bdd}.
\hfill $\square$

As a consequence of Proposition \ref{Prop-Rn-H23-bdd} and
the uniform estimates \eqref{D-Tt}-\eqref{thetan-Tt},
the asymptotic behavior \eqref{v-Sj-Sigma-RNLS} can be taken in the more regular space $\dot{H}^{\frac 32}$.
Since the proof is similar to that of \cite[Proposition 7.2]{SZ20},
it is omitted here for simplicity.

\begin{corollary} \label{Cor-unST-H23}
Consider the situation as in Proposition \ref{Prop-Rn-H23-bdd} with $\upsilon_*\geq 6$,
$m\geq 4$ if $d=2$ and $m\geq 5$ if $d=1$.
Then, we have
\begin{align} \label{un-ST-H23}
    \|v_n(t) - S(t)  -z(t) \|_{\dot H^\frac{3}{2}} \leq C (T-t)^{\frac \kappa 2 -2},
\end{align}
where $\kappa = (m+\frac d2-1)\wedge (\upsilon_*-2)$.
In particular,
for the blow-up solution $v_c$  constructed in Theorem \ref{Thm-BW-RNLS},
we have
\begin{align} \label{u-ST-H23}
    \|v_c(t) - S(t)- z(t) \|_{\dot H^\frac{3}{2}} \leq C (T-t)^{\frac \kappa 2 - 2},
\end{align}
and the strong $H^1$ convergence holds:
for any $t\in (t_*,T)$,
\begin{align} \label{un-u-0-H1}
    \|v_n - v_c\|_{C([t_*,t]; H^1)} \to 0,\ \ as\ n\to \9.
\end{align}
\end{corollary}

The constructed blow-up solution $v_c$ actually admits the geometrical decomposition on the
existing time interval $[t_*,T)$,
namely,
\begin{align} \label{vc-dec}
   v_c(t,x) =  \sum_{k=1}^{K}\lbb_{k}^{-\frac d2} Q_{k} (t,\frac{x-\a_{k}}{\lbb_{k}}) e^{i\theta_{k}}
              + z(t,x) + R(t,x) \ \
         \(:= U(t,x) + z(t,x)  + R(t,x)\)
\end{align}
with
\begin{align} \label{Qj-Q-Uniq}
   Q_{k}(t,y) := Q(y) e^{i(\beta_{k}(t) \cdot y - \frac 14 \g_{k}(t) {|y|^2})},
\end{align}
the parameters $\mathcal{P}:=\{\lbb, \a, \beta, \g, \theta\}$
are $C^1$ functions
and the following  orthogonality conditions hold on $[t_*,T)$:
for $1\leq k\leq K$,
\be\ba\label{ortho-cond-Rn-wn2}
&{\rm Re}\int (x-\a_{k}) U_{k}\ol{R}dx=0,\ \
{\rm Re} \int |x-\a_{k}|^2 U_{k} \ol{R}dx=0,\\
&{\rm Im}\int \nabla U_{k} \ol{R}dx=0,\ \
{\rm Im}\int \Lambda U_{k} \ol{R}dx=0,\ \
{\rm Im}\int \varrho_{k} \ol{R}dx=0.
\ea\ee
This fact is mainly due to the uniform estimate \eqref{Mod-lbb-boot} of modulation equation,
which ensures the equicontinuity of geometrical parameters $\{\calp_n\}$ on every $[t_*,t_n]$,
$n\geq 1$, and thus permits to take the limit procedure via the Arzel\`{a}-Ascoli Theorem.
We refer to \cite{SZ20} for more details.

Hence, taking the limit $n\to \9$ in the uniform estimates \eqref{D-Tt}-\eqref{thetan-Tt} and \eqref{Rn-H23-bdd}
we get the following estimates on $[t_*,T)$: for $1\leq k\leq K$,
\begin{align}
 & \|R(t)\|_{L^2}\leq (T-t)^{\kappa+1},\ \
 \|R(t)\|_{H^1}\leq (T-t)^{\kappa}, \ \
 \|  R(t)\|_{H^{\frac 32}}\leq (T-t)^{\kappa-2},  \label{R-Tt-Uniq} \\
&\lf|\la_{k}(t) -w_k(T-t) \rt| + \ \lf|\gamma_{k}(t)  - w_k^2(T-t) \rt|\leq (T-t)^{\kappa},   \label{lbb-Tt-Uniq} \\
&|\al_{k}(t)-x_k|+|\beta_{k}(t)|\leq (T-t)^{\frac \kappa 2 + \frac 12 },  \label{anb-Tt-Uniq}\\
&|\theta_{k}(t) - (w_k^{-2}(T-t)^{-1} +\vartheta_k)|\leq (T-t)^{\kappa-2}. \label{theta-Tt-Uniq}
\end{align}
As a consequence,
for any $t\in [t_*,T)$,
$\lbb_k, \g_k, P$ are comparable to $T-t$:
\begin{align} \label{lbbj-gj-P-Tt-Uniq}
    \lbb_k(t), \g_k(t), P(t)  \thickapprox T-t,
\end{align}
and
\begin{align}
 &Mod(t)\leq C \a^*(T-t)^{\kappa+1} , \label{Mod-bdd-Uniq} \\
 & \|\partial_x^\upsilon \eta\|_{L^2} \leq C \a^* (T-t)^{\kappa-1-|\upsilon|} , \ \   t\in [t_*,T),\ |\upsilon|\leq 2. \label{eta-bdd-Uniq}
\end{align}
where $C>0$ is a universal constant independent of $\ve, \a^*$ and $t$.

\section{Conditional uniqueness of multi-bubble Bourgain-Wang solutions}  \label{Sec-Uniq-BW}

\subsection{Control of the difference}

In this subsection we assume Hypothesis $(H1)$ with $m\geq 10$, $\upsilon_*\geq 12$.
Set $\kappa:= (m+\frac d2-1)\wedge (\upsilon_*-2)$. Note that,
$\kappa \geq 9+\frac d2$.

Let $v_c$ be the constructed multi-bubble Bourgain-Wang solution in Theorem \ref{Thm-BW-RNLS},
with the corresponding parameters
$\calp=(\lbb, \a, \beta, \g, \theta)$.
Let $v$ be any blow-up solution to \eqref{equa-NLS-perturb}
satisfying
\begin{align} \label{v-Sj-t}
 \|v(t)-\sum_{k=1}^KS_k(t)-z(t)\|_{L^2} + (T-t)\|v(t)-\sum_{k=1}^KS_k(t)-z(t)\|_{H^1}\leq C(T-t)^{4+\zeta},\ \ t\in [t_*,T),
\end{align}
where $\zeta$ is any positive constant close to $0$.
Set
\begin{align} \label{w-wj}
   w:=v-v_c =\sum_{k=1}^{K}w_k,\  \ w_k :=w\Phi_k,\ 1\leq k\leq K,
\end{align}
where $\{\Phi_k\}$ are   given by
\eqref{phi-local}.
Define the renormalized variable $\epsilon_k$ by
\begin{align} \label{wj-vej-def}
w_k(t,x) :=\lbb_{k}(t)^{-\frac d2} \epsilon_{k} (t,\frac{x-\a_{k}(t)}{\lbb_{k}(t)}) e^{i\theta_{k}(t)}, \ \ 1\leq k\leq K.
\end{align}
Note that, $\epsilon_k$ is different from $\ve_k$ defined in \eqref{Rj-ej}.
Similarly to \eqref{D-def}, set
\begin{align} \label{wtD-def}
    \wt D(t):= \| w(t)\|_{L^2} + (T-t) \|\nabla w(t)\|_{L^2},
\end{align}

Then, by  \eqref{R-Tt-Uniq} and  \eqref{v-Sj-t},
\begin{align}
    \|R(t)\|_{L^2}\leq (T-t)^{\kappa+1}, \
    \|R(t)\|_{H^1}&  \leq (T-t)^{\kappa}, \
 \|  R(t)\|_{H^{\frac 32}}\leq (T-t)^{\kappa -2}\ \ with\ \kappa \geq 9, \label{D-Tt-Uniq}\\
   & \wt D(t) \leq C(T-t)^{4+\zeta},  \label{wtD-Tt-Uniq} \ \
\end{align}
and
\begin{align}
   \|w(t) \|_{L^p}^p \leq C (T-t)^{-d(\frac p2-1)} \wt D^p.   \label{wtD-Lp-Tt}
\end{align}

Moreover, by equations \eqref{equa-NLS-perturb} and \eqref{wtD-Tt-Uniq},
$w$ satisfies the equation
\be   \label{equa-w}
\left\{ \begin{aligned}
   & i\partial_t w +\Delta w + a_1 \cdot \nabla w+ a_0 w +f(v_c+w)-f(v_c)=0, \ \ t\in (t_*,T), \\
   & \lim_{t\to T} \|w(t)\|_{H^1} =0.
\end{aligned}\right.
\ee

The crucial ingredient in the  uniqueness proof is
the following Lyapunov type functional,
which is similar to the generalized energy $\scri$ in \eqref{I-def},
\begin{align} \label{def-wtI}
 \wt {\scri} := &\frac{1}{2}\int |\nabla w|^2dx
            +\frac{1}{2}\sum_{k=1}^K \frac{1}{\lambda_{k}^2} \int  |w|^2 \Phi_k dx
           -{\rm Re}\int F(v_c+w)-F(v_c)-f(v_c)\ol{w}dx \nonumber \\
&+\sum_{k=1}^K\frac{\gamma_{k}}{2\lambda_{k}}{\rm Im} \int (\nabla\chi_A) \(\frac{x-\alpha_{k}}{\lambda_{k}}\)\cdot\nabla w \ol{w}\Phi_kdx.
\end{align}

\begin{lemma}    \label{Lem-wtI-Nt}
There exist $C_1, C_2, C_3>0$ such that for $t\in[t_*,T)$,
\begin{align}  \label{wtI-Nt}
C_1 (T-t)^{-2} \wt D^2- C_2\sum_{k=1}^{K}\frac{Scal_k}{\lbb_k^2}
\leq  \wt {\scri}
\leq C_3 A (T-t)^{-2} \wt D^2,
\end{align}
where
\begin{align} \label{Scaj-def}
Scal_k(t):= \<\epsilon_{k,1},Q\>^2+\<\epsilon_{k,1},yQ\>^2+\<\epsilon_{k,1},|y|^2Q\>^2
           +\<\epsilon_{k,2},\nabla Q\>^2+\<\epsilon_{k,2},\Lambda Q\>^2+\<\epsilon_{k,2},\rho\>^2,
\end{align}
and $\epsilon_{k,1}, \epsilon_{k,2}$ are the real and imaginary parts of $\epsilon_k$, respectively.
\end{lemma}

{\bf Proof.}
We first show that,
the constructed blow-up solution $v_c$ in \eqref{def-wtI}
can be replaced by the blow-up profile $U$ given by \eqref{vc-dec},
up to the error $\calo((T-t)^{-1}\wt D^2)$, i.e.,
\begin{align} \label{Fuw-Fu-fuw}
   {\rm Re}\int F( v_c+w)-F(v_c)-f(v_c)\ol{w}dx
= & {\rm Re}\int F(U+w)-F(U)-f(U)\ol{w}dx \nonumber \\
  & + o\((T-t)^{-1} \wt D^2\).
\end{align}
To this end, we note that
\begin{align}
   |F''(v_c,w)\cdot w^2 - F''(U,w)\cdot w^2|
   \leq C\( |U|^{\frac 4d-1} +  |w|^{\frac 4d-1} +  |z|^{\frac 4d-1} +  |R|^{\frac 4d-1} \) |z+R||w|^2.
\end{align}
By \eqref{eyvez-bdd}, \eqref{R-Tt-Uniq} and \eqref{wtD-Tt-Uniq},
\begin{align} \label{U-Rw2-esti}
    \int |U|^{\frac 4d-1} |z+R||w|^2 dx
    &\leq C  (T-t)^{-2} \sum\limits_{k=1}^K \|e^{-\delta|y|}\ve_{z,k}\|_{L^\9} \|w\|_{L^2}^2
          +C(T-t)^{-\frac d2(\frac 4d-1)} \|R\|_{L^2} \|w\|^2_{H^1}
          + C e^{-\frac{\delta}{T-t}} \|w\|_{L^2}^2  \nonumber  \\
    &\leq C \( \a^* (T-t)^{m-1+\frac d2} \wt D^2 + (T-t)^{\kappa-3+\frac d2} \wt D^2 + e^{-\frac{\delta}{T-t}} \wt D^2 \) \nonumber \\
    &= o\((T-t)^{-1} \wt D^2 \).
\end{align}
Moreover,
by  \eqref{z-LtHm-a*}, \eqref{R-Tt-Uniq} and \eqref{wtD-Tt-Uniq},
\begin{align} \label{w-R-Rw2}
     \int(|w|^{\frac 4d-1} + |z|^{\frac 4d-1}+ |R|^{\frac 4d-1} ) |z+R||w|^2dx
   \leq& C \bigg(\|z\|_{L^\9}\|w\|_{H^1}^{1+\frac 4d}
             + \|R\|_{L^2} \|w\|_{H^1}^{1+\frac 4d}
          + \|z\|_{L^\9}^{\frac 4d} \|w\|_{L^2}^2 \nonumber \\
        &  + \|z\|_{L^\9}^{\frac 4d-1} \|R\|_{H^1} \|w\|_{H^1}^2
          + \|z\|_{L^\9} \|R\|_{H^1}^{\frac 4d-1} \|w\|_{H^1}^2
          + \|R\|^{\frac 4d}_{H^1} \|w\|_{H^1}^{2}\bigg)  \nonumber  \\
     =& o\( (T-t)^{-1} \wt D^2\).
\end{align}
Hence, \eqref{Fuw-Fu-fuw} follows from \eqref{U-Rw2-esti} and \eqref{w-R-Rw2}, as claimed.

Next, for the R.H.S. of \eqref{Fuw-Fu-fuw},   note that,
\begin{align} \label{FUw-expan2}
  {\rm Re} (F(U+w) - F(U) - f(U)\ol w)
  =&   \frac 12 (1+\frac 2d) |U|^{\frac 4d} |w|^2
     + \frac{1}{d} |U|^{\frac 4d -2} {\rm Re} (U^2 \ol{w}^2)  \nonumber  \\
   &  + \calo \((|U|^{\frac 4d-1} + |w|^{\frac 4d-1})|w|^3 \).
\end{align}
The error term above can be bounded by,
via \eqref{wtD-Tt-Uniq},
\begin{align}  \label{U4d-w3}
   \int (|U|^{\frac 4d-1} + |w|^{\frac 4d-1})|w|^3 dx
    \leq& C\((T-t)^{-2+\frac d2} \|w\|^3_{H^1} + \|w\|_{H^1}^{2+\frac 4d} \)
   \leq C (T-t)^{-1} \wt D^2.
\end{align}
Moreover, for the Morawetz type functional in \eqref{def-wtI},
\begin{align} \label{wtI-last-bdd}
   \bigg| \frac{\gamma_{k}}{2\lambda_{k}}{\rm Im} \int (\nabla\chi_A) \(\frac{x-\alpha_{k}}{\lambda_{k}}\)\cdot\nabla w \ol{w}\Phi_kdx \bigg|
   \leq C A \|w\|_{L^2} \|\na w\|_{L^2}
   \leq C A (T-t)^{-1} \wt D^2(t).
\end{align}

Thus, we conclude from \eqref{Fuw-Fu-fuw}, \eqref{FUw-expan2}-\eqref{wtI-last-bdd} that
\begin{align} \label{wtI-expan2nd}
   \wt{\scri}
  =   \frac{1}{2} {\rm Re}
       \int |\nabla w|^2
         + \sum_{k=1}^K  \frac{1}{\lambda_{k}^2}|w|^2\Phi_k
       - (1+\frac 2d) |U|^{\frac 4d} |w|^2
         - \frac{2}{d}|U|^{\frac 4d -2} U^2 \ol{w}^2 dx
     + \calo(A(T-t)^{-1} \wt D^2).
\end{align}

Now, on one hand,
by H\"older's inequality, \eqref{wtI-expan2nd} and \eqref{wtD-def},
\begin{align}
   |\wt \scri| \leq C\( \|w\|^2_{H^1} + (T-t)^{-2} \|w\|^2_{L^2} + (T-t)^{-1} \wt D^2 \)
   \leq C A (T-t)^{-2} \wt D^2,
\end{align}
which yields the second inequality in \eqref{wtI-Nt}.

On the other hand,
the first inequality in \eqref{wtI-Nt}  mainly follows from the coercivity type estimate below,
which is similar to \eqref{naR2-R2-coer}  mainly due to the local coercivity of
linearized operators,
\begin{align*}
  \wt \scri
  \geq  C_1 (T-t)^{-2} \wt D^2
       - C_2 \( A (T-t)^{-1} \wt D^2 + \sum\limits_{k=1}^K \lbb_k^{-2} Scal_k + e^{-\frac{\delta}{T-t}} \wt D^2 \).
\end{align*}
Hence, for $t$ close to $T$
such that
$C_2\(A (T-t)+e^{-\frac{\delta}{T-t}}\) \leq \frac 12 C_1$,
it leads to
\begin{align*}
  \wt \scri
  \geq  \frac 12 C_1 (T-t)^{-2} \wt D^2
         - C_2 \sum_{k=1}^{K} \lbb_k^{-2} Scal_k.
\end{align*}
This verifies the first inequality in \eqref{wtI-Nt}.
Therefore, the proof is complete.
\hfill $\square$

The following monotonicity property of $\wt \scri$ is crucial in the derivation of uniqueness.

\begin{theorem} (Monotonicity of $\wt \scri$)  \label{Thm-wtI-mono}
There exist $C_1, C_2>0$ such that
for $A$ large enough and $t$ close to $T$,
\begin{align}  \label{dIt-mono-new}
\frac{d \wt {\scri}}{dt}
\geq& C_1 \sum_{k=1}^{K}
      \int (\frac{1}{\lambda_{k}} |\nabla w_{k}|^2+\frac{1}{\lbb_{k}^3} |w_{k}|^2 )
      e^{-\frac{|x-\alpha_{k}|}{A\lbb_k}}dx  - C_2 A \wt \cale_r
\end{align}
where
\begin{align} \label{wtcalEr-def}
    \wt \cale_r = \frac{\wt D^2}{(T-t)^2} + \ve \frac{\wt D^2}{(T-t)^3} + \sum_{k=1}^{K}\frac{Scal_k(t)}{\lbb_k^3(t)}.
\end{align}
\end{theorem}

\begin{remark}
Comparing with the error $\mathcal{E}_r$ in \eqref{calEr-def},
we see that $\wt \cale_r$ in \eqref{wtcalEr-def} only
contains the orders of $D$  higher than one.
This fact is important in the derivation of uniqueness.
\end{remark}

{\bf Proof. }
Using equation \eqref{equa-w} we compute
\begin{align} \label{equa-wtI}
 \frac{d \wt{\scri}}{dt}
=&-\sum_{k=1}^K\frac{\dot{\lbb}_{k}}{\lbb_{k}^3} {\rm Im}\int|w|^2 \Phi_kdx
-\sum_{k=1}^K\frac{1}{\lbb_{k}^2}{\rm Im} \< f^\prime(v_c)\cdot w, {w}_k\>
   -{\rm Re} \<f''(v_c,w)\cdot w^2, \pa_t v_c\> \nonumber\\
& - \sum_{k=1}^K\frac{1}{\lbb_{k}^2}{\rm Im} \< w\nabla \Phi_k, { \nabla w}\>
 -\sum_{k=1}^K\frac{1}{\lbb_{k}^2} {\rm Im} \< f''(v_c,w)\cdot w^2, {w}_k\> \nonumber\\
& - {\rm Im}
      \<\Delta w - \sum\limits_{k=1}^K \frac{1}{\lbb_k^2} w_k +f(v_c+w)-f(v_c), {a_1\cdot \nabla w+a_0 w} \>   \nonumber\\
&-\sum_{k=1}^K\frac{\dot{\lambda}_{k}\gamma_{k}-\lambda_{k}\dot{\gamma}_{k}}{2\lambda_{k}^2}
        {\rm Im}  \< \nabla\chi_A(\frac{x-\alpha_{k}}{\lambda_{k}})\cdot\nabla w, {w}_k\>
  +\sum_{k=1}^K\frac{\gamma_{k}}{2\lambda_{k}}{\rm Im}
        \< \partial_t (\nabla\chi_A(\frac{x-\alpha_{k}}{\lambda_{k}}))\cdot\nabla w, {w}_k\> \nonumber \\
&+\sum_{k=1}^K {\rm Im}  \< \frac{\gamma_{k}}{2\lambda_{k}^2} \Delta\chi_A(\frac{x-\alpha_{k}}{\lambda_{k}})w_k
+\frac{\gamma_{k}}{2\lambda_{k}}  \nabla\chi_A(\frac{x-\alpha_{k}}{\lambda_{k}})\cdot ( \nabla w_k + \na w \Phi_k),  \pa_t {w} \> \nonumber \\
=:&\sum_{l=1}^{9}\wt{\scri}_{t,l}.
\end{align}
In order to reduce the analysis of \eqref{equa-wtI} to
the previous case in \eqref{equa-I1t} and \eqref{equa-I2t},
we show that the reference solution $v_c$ in $\wt{\scri}_{t,2}$, $\wt{\scri}_{t,3}$, $\wt{\scri}_{t,5}$,
$\wt{\scri}_{t,6}$ and $\wt{\scri}_{t,9}$
can be replaced by $U+z$,
up to the acceptable error $(T-t)^{-2} \wt D^2$.

{\it $(i)$ Estimate of  $\wt \scri_{t,2}$}.
By \eqref{f'v1-f'v2},
\begin{align*}
 & \bigg|\wt{\scri}_{t,2}+\sum_{k=1}^K\frac{1}{\lbb_{k}^2}{\rm Im}
       \< f^\prime(U+z)\cdot w, {w}_k\>  \bigg| \nonumber \\
\leq& C \sum_{k=1}^K\frac{1}{\lbb_{k}^2}
        \int (|U|^{\frac{4}{d}-1}+ |z|^{\frac 4d-1}  +|R|^{\frac 4d-1} )|R||w|^2dx \nonumber \\
\leq& C\( (T-t)^{-4+\frac d2}  \|R\|_{L^2} \|w\|_{L^4}^2
         + (T-t)^{-2} \|z\|_{L^\9}^{\frac 4d-1} \|R\|_{L^2} \|w\|^2_{L^4}
         + (T-t)^{-2} \|R\|_{H^1}^{\frac 4d}  \|w\|^2_{L^4}  \).
\end{align*}
Then, by \eqref{z-LtHm-a*}, \eqref{R-Tt-Uniq} and  \eqref{wtD-Lp-Tt},
\begin{align}  \label{wtIt3-U-Nt}
    \bigg|\wt{\scri}_{t,2}+\sum_{k=1}^K\frac{1}{\lbb_{k}^2}{\rm Im}
       \< f^\prime(U+z)\cdot w, {w}_k\>  \bigg|
   \leq& C\( (T-t)^{\kappa-3} + (T-t)^{\frac 4d \kappa - \frac d2 -2} \) \wt D^2   \nonumber \\
   \leq& C (T-t)^{-2} \wt D^2.
\end{align}

{\it $(ii)$ Estimate of $\wt{\scri}_{t,3}$}.
By the decomposition \eqref{vc-dec},
\begin{align*}
   {\rm Re} \< f''(v_c,w)\cdot w^2,  \partial_t {v_c} \>
   = {\rm Re} \< f''(v_c,w)\cdot w^2,  \partial_t ({U}+z) \>
     + {\rm Re} \< f''(v_c,w)\cdot w^2,  \partial_t {R} \>.
\end{align*}
Let us treat the two terms on the R.H.S. above separately.

First by \eqref{f''v1-f''v2},
\begin{align} \label{f''vc-f''Uz-esti0}
&|{\rm Re} \< f''(v_c,w)\cdot w^2,  \partial_t (U+z) \>
   -{\rm Re} \< f''(U+z,w)\cdot w^2, \partial_t  (U+z)\>| \nonumber \\
\leq& C \|\partial_t  (U+z) \|_{L^\infty}
    \int \(|U|^{\frac{4}{d}-2}+ |z|^{\frac{4}{d}-2} +|R|^{\frac{4}{d}-2}+ |w|^{\frac{4}{d}-2}\) |R||w|^2dx.
\end{align}
Since by \eqref{eyvez-bdd} and \eqref{pt-Uj-tx},
$\|\partial_t (U+z)\|_{L^\9} \leq C (T-t)^{-2-\frac d2}$.
Then, by \eqref{R-Tt-Uniq} and \eqref{wtD-Tt-Uniq},
the R.H.S. above can be bounded by, up to a universal constant,
\begin{align} \label{f''uww2-ptU}
      & (T-t)^{-2-\frac d2} \bigg(\((T-t)^{-2+d}+ \|z\|_{L^\9}^{\frac 4d-2}\) \|R\|_{L^2} \|w\|_{L^4}^2
      +   \|R\|_{H^1}^{\frac 4d-1} \|w\|_{L^4}^2
      +  \|R\|_{L^2} \|w\|_{L^{\frac 8d}}^{\frac 4d} \bigg) \nonumber \\
    \leq& (T-t)^{\kappa-3} \wt D^2
    \leq  (T-t)^{-2} \wt D^2.
\end{align}

Next we show that
\begin{align} \label{f''uww2-ptR}
     {\rm Re} \< f''(v_c,w)\cdot w^2,  \partial_t {R} \>
   = \calo((T-t)^{-2} \wt D^2).
\end{align}
To this end, by equation (\ref{equa-R}),
\begin{align*}
    |{\rm Re} \< f''(v_c,w)\cdot w^2,   \partial_t {R} \> |
= |{\rm Im} \< f''(v_c,w)\cdot w^2 , \Delta R+f(v_c)-f(U+z)+ (a_1\cdot \nabla +a_0)R+\eta \>|.
\end{align*}
Note that, by \eqref{Rn-H23-bdd},
\begin{align} \label{f''uw2DR}
  |{\rm Im} \< f''(v_c,w)\cdot w^2, \Delta R \>|
 \leq& C\|R\|_{\dot{H}^{\frac{3}{2}}}\|f''(u,w)\cdot w^2\|_{\dot{H}^{\frac{1}{2}}}
 \leq  C (T-t)^{\kappa-2} \|f''(u,w)\cdot w^2\|_{\dot{H}^{\frac{1}{2}}} .
\end{align}
Then, by \eqref{z-LtHm-a*}, \eqref{vc-dec}, \eqref{R-Tt-Uniq}, \eqref{wtD-Tt-Uniq}
and $\|U(t)\|_{H^1} \leq C(T-t)^{-1}$,
\begin{align} \label{f''uw2-H12}
   \|f''(v_c,w)\cdot w^2  \|_{\dot{H}^{\frac{1}{2}}}
   \leq& C \sum\limits_{j=2}^{1+\frac 4d}
          \|v_c\|_{H^1}^{1+\frac 4d -j} \|w\|_{H^1}^j   \nonumber \\
   \leq& C \sum\limits_{j=2}^{1+\frac 4d}
          (\|U\|_{H^1}^{1+\frac 4d -j} + \|z\|_{H^1}^{1+\frac 4d -j} + \|R\|_{H^1}^{1+\frac 4d -j} ) \|w\|_{H^1}^j \nonumber \\
   \leq&  C  \sum\limits_{j=2}^{1+\frac 4d} (T-t)^{-(1+\frac 4d -j)+(3+\zeta)(j-2)} \|w\|_{H^1}^2
   \leq C (T-t)^{-3} \|w\|_{H^1}^2.
\end{align}
Plugging this into \eqref{f''uw2DR} and using $\kappa \geq 5$
we obtain
\begin{align} \label{f''uww2-DR}
     |{\rm Im} \< f''(v_c,w)\cdot w^2,  \Delta {R} \>|
    \leq C (T-t)^{\kappa-5} \|w\|_{H^1}^2
   \leq C (T-t)^{-2} \wt D^2(t).
\end{align}
Moreover,
since by \eqref{f''-bdd},
\begin{align} \label{f''uw-w2-bdd}
   |f''(v_c,w)\cdot w^2|
   \leq& C (|U|^{\frac 4d-1} +  |z|^{\frac 4d-1} + |R|^{\frac 4d -1} + |w|^{\frac 4d -1}) |w|^2 \nonumber \\
   \leq& C  \((T-t)^{-2+\frac d2} + |R|^{\frac 4d-1} + |w|^{\frac 4d-1} \) |w|^2,
\end{align}
and
\begin{align} \label{fu-fU-bdd}
   |f(v_c)-f(U+z)|
   \leq C(|U|^{\frac 4d}+|R|^{\frac 4d}+ |z|^{\frac 4d})|R|
   \leq C \((T-t)^{-2} + |R|^\frac 4d \) |R|,
\end{align}
taking into account \eqref{z-LtHm-a*}, \eqref{R-Tt-Uniq} and \eqref{wtD-Tt-Uniq} we get
\begin{align} \label{f''uww2-fufUbcR}
   & \bigg|{\rm Im} \< f''(v_c,w)\cdot w^2,  f(v_c)-f(U+z)+ (a_1\cdot \nabla +a_0) R \> \bigg| \nonumber \\
   \leq& C \int  \((T-t)^{-2+\frac d2} + |R|^{\frac 4d-1} + |w|^{\frac 4d-1} \) |w|^2
                  \((T-t)^{-2}  |R| + |R|^{1+\frac 4d} + |\na R| + |R| \)   dx  \nonumber \\
   \leq& C   \((T-t)^{-2+\frac d2} + \|R\|_{H^1}^{\frac 4d-1} + \|w\|_{H^1}^{\frac 4d-1} \) \|w\|_{H^1}^2
                  \((T-t)^{-2}  \|R\|_{L^2} + \|R\|^{1+\frac 4d}_{H^1} + \|R\|_{H^1} \)  \nonumber \\
    \leq& C (T-t)^{\kappa-5+\frac d2} \wt D^2
    \leq  C (T-t)^{-2} \wt D^2.
\end{align}
Furthermore,
by  \eqref{eta-bdd-Uniq} and \eqref{f''uw-w2-bdd},
\begin{align} \label{f''uww2-esta}
  |{\rm Im} \< f''(v_c,w)\cdot w^2, \eta \>|
  \leq& C  (T-t)^{-2+\frac d2}   \|\eta\|_{L^2} \|w\|_{H^1}^2  \nonumber \\
  \leq&  C  (T-t)^{\kappa-3} \|w\|_{H^1}^2   \nonumber \\
  \leq&  C (T-t)^{-2} \wt D^2.
\end{align}

Thus, estimates \eqref{f''uww2-DR}, \eqref{f''uww2-fufUbcR} and \eqref{f''uww2-esta}
together yield  \eqref{f''uww2-ptR}, as claimed.

Therefore, we conclude from \eqref{f''vc-f''Uz-esti0}, \eqref{f''uww2-ptU} and \eqref{f''uww2-ptR} that
\begin{align}
   \wt \scri_{t,3}
    = {\rm Re} \<f''(U+z,w)\cdot w^2, \pa_t (U+z)\> + \calo\((T-t)^{-2} \wt D^2\).
\end{align}

{\it $(iii)$ Estimate of $\wt \scri_{t,5}$}.
By \eqref{D-Tt-Uniq}, \eqref{wtD-Tt-Uniq} and \eqref{f''v1-f''v2},
\begin{align*}
   & \bigg|\wt{\scri}_{t,5}+\sum_{k=1}^K\frac{1}{\lbb_{k}^2}{\rm Im} \< f''(U+z,w)\cdot w^2, {w}_k\> \bigg| \nonumber \\
\leq& C (T-t)^{-2} \int \(|U|^{\frac 4d-2} + |z|^{\frac 4d-2} + |R|^{\frac 4d-2} +  |w|^{\frac 4d-2}\) |R||w|^3 dx \nonumber \\
\leq& C(T-t)^{-4+d} \|R\|_{H^1} \|w\|_{H^1}^3
\leq C (T-t)^{-2} \wt D^2.
\end{align*}
This yields that
\begin{align} \label{wtI5-U-N}
   \wt \scri_{t,5}
   = - \sum\limits_{j=1}^K \frac{1}{\lbb_k^2}
     {\rm Im}\< f''(U,w)\cdot w^2, w_k\>
     + \calo\((T-t)^{-2} \wt D^2\).
\end{align}

{\it $(iv)$ Estimate of $\wt \scri_{t,6}$}.
Since by \eqref{f'R},
\begin{align} \label{fuwfu-fUwfU}
     & |f(v_c+w)-f(v_c) -(f(U+z+w)-f(U+z))|    \nonumber \\
 = & |f'(v_c,w)\cdot w - f'(U+z, w)\cdot w| \nonumber \\
\leq& C (|U|^{\frac 4d-1} +|z|^{\frac 4d-1} +  |R|^{\frac 4d-1} + |w|^{\frac 4d-1} )|R||w|,
\end{align}
we infer from \eqref{R-Tt-Uniq} that
\begin{align*}
   & \bigg| {\rm Im} \< f(v_c+w)-f(v_c) -(f(U+z+w)-f(U+z)),{ a_1\cdot \na w +a_0 w} \>  \bigg| \\
   \leq& C (T-t)^{-\frac d2(\frac 4d-1)} \|R\|_{H^1} \|w\|_{H^1}^2 \nonumber \\
   \leq& C (T-t)^{-2} \wt D^2 .
\end{align*}
This yields that
\begin{align} \label{wtIt6-U-Nt}
   \wt \scri_{t,6}
   = - {\rm Im}
     \<\Delta w - \sum\limits_{k=1}^K  \frac{1}{\lbb_k^2} w_k +f(U+z+w)-f(U+z), { a_1\cdot \nabla w +a_0 w  }  \>
    + \calo((T-t)^{-2} \wt D^2).
\end{align}

{\it $(v)$ Estimate of $\wt \scri_{t,9}$}.
By equation \eqref{equa-w},
\begin{align} \label{wtI9-esti0}
  \wt \scri_{t,9}
  =& \sum\limits_{k=1}^K {\rm Im}  \big\< \frac{\gamma_{k}}{2\lambda_{k}^2} \Delta\chi_A(\frac{x-\alpha_{k}}{\lambda_{k}})w_k
+\frac{\gamma_{k}}{2\lambda_{k}}  \nabla\chi_A(\frac{x-\alpha_{k}}{\lambda_{k}})\cdot ( \nabla w_k + \na w \Phi_k), \nonumber  \\
   &\qquad  i\Delta w + i ( a_1 \cdot \na w + a_0  w ) + i(f(v_c+w)-f(v_c))  \big\> .
\end{align}
Note that, unlike in \eqref{equa-R}, we have $\eta=0$ here.

Then, in view of \eqref{fuwfu-fUwfU}, we see that
\begin{align*}
  & \bigg|\< \frac{\gamma_{k}}{2\lambda_{k}^2} \Delta\chi_A(\frac{x-\alpha_{k}}{\lambda_{k}})w_k
      +\frac{\gamma_{k}}{2\lambda_{k}}  \nabla\chi_A(\frac{x-\alpha_{k}}{\lambda_{k}})\cdot ( \nabla w_k+ \na w \Phi_k),
      i(f(v_c+w)-f(v_c)) \>    \\
  & \ - \< \frac{\gamma_{k}}{2\lambda_{k}^2} \Delta\chi_A(\frac{x-\alpha_{k}}{\lambda_{k}})w_k
         +\frac{\gamma_{k}}{2\lambda_{k}}  \nabla\chi_A(\frac{x-\alpha_{k}}{\lambda_{k}})\cdot ( \nabla w_k + \na w \Phi_k),
         i(f(U+z+w)-f(U+z)) \> \bigg| \\
  \leq&C A\int \((T-t)^{-1} |w| + |\na w|\)
              \(|U|^{\frac 4d-1} + |z|^{\frac 4d-1} +|w|^{\frac 4d-1} +|R|^{\frac 4d-1}\) |R||w| dx \nonumber \\
  \leq& CA (T-t)^{-3+\frac d2} \|R\|_{H^1} \|w\|_{H^1}^2
  \leq CA(T-t)^{-2} \wt D^2.
\end{align*}
This yields that
\begin{align} \label{wtIt9-U-Nt}
   \wt \scri_{t,9}
   =&  \sum\limits_{k=1}^K  {\rm Im} \big\< \frac{\gamma_{k}}{2\lambda_{k}^2} \Delta\chi_A(\frac{x-\alpha_{k}}{\lambda_{j}})w_k
+\frac{\gamma_{k}}{2\lambda_{k}}  \nabla\chi_A(\frac{x-\alpha_{k}}{\lambda_{k}})\cdot ( \nabla w_k + \na w \Phi_k), \nonumber  \\
   &\qquad  i\Delta w + i(f(U+z+w)-f(U+z)) + i (a_1\cdot \na +a_0)w \big\>
     + \mathcal{O}\(A(T-t)^{-2} \wt D^2\).
\end{align}

Now,
the reference solution $v_c$ in \eqref{equa-wtI}
has been replaced by $U+z$
up to the order $\calo((T-t)^{-2} \wt D^2)$.
Note that,
by \eqref{R-Tt-Uniq}-\eqref{theta-Tt-Uniq} and \eqref{wtD-Tt-Uniq},
the conditions in Theorem \ref{Thm-I-mono} are verified.
Hence, arguing as in the proof of Theorem \ref{Thm-I-mono}
with $w$ replacing $R$ and using \eqref{wtD-Tt-Uniq}
we obtain \eqref{dIt-mono-new}.

As mentioned below \eqref{wtI9-esti0},
because for the difference $w$ we have $\eta =0$,
the errors involving $M_k$ and the linear terms of $D$
in \eqref{calEr-def} do not appear here,
only the higher order terms of $D$ remain.

Therefore,  the proof is complete.
\hfill $\square$

As a consequence of Lemma \ref{Lem-wtI-Nt} and  Theorem \ref{Thm-wtI-mono}
we have

\begin{corollary} \label{Cor-Nt-Scal}
For $t$ close to $T$, set
\begin{align} \label{wtN-def}
   \wt N(t) := \sup\limits_{t\leq s<T} \frac{\wt D^2(s)}{(T-s)^2}.
\end{align}
Then, there exists $C>0$ such that
\begin{align} \label{Nt-Scal}
    \wt N(t)
\leq C \(\sum_{k=1}^{K} \sup\limits_{t\leq s < T} \frac{Scal_k(s)}{\lbb_k^2(s)}
+ \int_{t}^{T}\sum_{k=1}^{K}\frac{Scal_k(s)}{\lbb_k^3(s)}  + \ve \frac{\wt N(s)}{T-s}ds \).
\end{align}
\end{corollary}

{\bf Proof.}
By Lemma \ref{Lem-wtI-Nt} and Theorem \ref{Thm-wtI-mono},
for $t<\tilde{t} <T$,
\begin{align*}
     C_1 \frac{\wt D^2(t)}{(T-t)^2}
&\leq  \wt {\scri}(t) +C_2\sum_{k=1}^{K}\frac{Scal_k(t)}{\lbb_k^2(t)} \nonumber \\
   & =  \wt \scri(\tilde{t})   +C_2\sum_{k=1}^{K}\frac{Scal_k(t)}{\lbb_k^2(t)}
-\int_{t}^{\tilde{t}}\frac{d \wt{\scri}}{ds}(s) ds \nonumber \\
&\leq  CA \( \frac{\wt D^2(\tilde{t})}{(T-\tilde{t})^{2}} + \sum_{k=1}^{K}\frac{Scal_k(t)}{\lbb_k^2(t)}
       + \int_{t}^{\tilde{t}} \frac{\wt D^2(s)}{(T-s)^2}  + \sum_{k=1}^{K}\frac{Scal_k(s)}{\lbb_k^3(s)}
       +  \ve \frac{\wt D^2(s)}{(T-s)^3}ds \),
\end{align*}
which yields that
\begin{align*}
   \sup\limits_{t\leq s\leq \tilde{t}} \frac{\wt D^2(s)}{(T-s)^2}
    \leq C A \( \wt N(\tilde{t})
             + \sum_{k=1}^{K} \sup\limits_{t\leq s\leq \wt t} \frac{Scal_k(s)}{\lbb_k^2(s)}
              + (\tilde{t}- t)  \wt N(t)
         +\int_t^{\tilde{t}} \sum\limits_{k=1}^K \frac{Scal_k(s)}{\lbb_k^3(s)}
         + \ve  \frac{\wt N(s)}{T-s}ds\).
\end{align*}
Since by \eqref{wtD-Tt-Uniq},
$\wt N(\tilde{t}) \to 0$ as $\tilde{t} \to T$,
taking $\tilde{t} \to T$ and
$t$ close to $T$
we obtain \eqref{Nt-Scal}.
\hfill $\square$

\subsection{Control of the null space}  \label{Subsec-Cont-Null}

In this subsection we derive the control of  scalar $Scal_k$.
The main result is formulated in Theorem \ref{Thm-Scal} below.
The arguments follow the lines in the proof of \cite[Theorem 7.7]{SZ20},
mainly based on algebraic identities.
For the reader's convenience,
let us sketch the main arguments below.

For every $1\leq k\leq K$, define the renormalized variables $\wt e_k$ and $e_k$ by
\begin{align} \label{w-wtej-ej-def}
w(t,x)=\lbb_{k}(t)^{-\frac d2} \wt {e}_{k} (t,\frac{x-\a_{k}(t)}{\lbb_{k}(t)}) e^{i\theta_{k}(t)}, \ \
with\ \wt e_k(t,y)=e_k(t,y)e^{i(\b_k(t)\cdot y- \frac 14 \g_k(t) |y|^2)}.
\end{align}
Note that,
the renormalized variable $e_k$ is
different from the previous one $\epsilon_k$  in \eqref{wj-vej-def}.

We use  \eqref{g-gzz-expan} to expand
\begin{align}
f(v_c+w)-f(v_c)= \partial_z f(v_c) w + \partial_{\ol z} f(v_c) \ol{w} + f''(v_c,w)\cdot w^2.
\end{align}
Then, using \eqref{g-gzz-expan} again to further expand $ \partial_z f(v_c)$ and $\partial_{\ol z} f(v_c)$
around the profile $U$
we get
\begin{align}   \label{f-G1}
   f(v_c+w) - f(v_c) = f'(U)\cdot w + G_1,
\end{align}
where
\begin{align} \label{G1-equaw}
   G_1:=
    w(\partial_z f)'(U, z+R)\cdot (z+R)
    + \ol w (\partial_{\ol z} f)'(U, z+R)\cdot (z+R)
    + f''(v_c, w)\cdot w^2.
\end{align}
Decompose $f'(U)\cdot w$ into three parts
\begin{align} \label{f'w-G2G3}
   f'(U)\cdot w
   &= f'(U_k)\cdot w
      + \sum\limits_{l\not = k} f'(U_l)\cdot w
      + [f'(U) \cdot w - \sum\limits_{l=1}^K f'(U_l) \cdot w] \nonumber \\
   &=: f'(U_k)\cdot w + G_2 + G_3,
\end{align}
and set
\begin{align} \label{G4-bc}
   G_4 :=a_1\cdot \na w+ a_0 w ,
\end{align}
where $a_1,a_0$ are given by \eqref{a1-loworder} and \eqref{a0-loworder}, respectively.

Thus, by \eqref{f-G1}, \eqref{f'w-G2G3} and \eqref{G4-bc},
equation \eqref{equa-w} can be reformulated:
\begin{align} \label{equa-w-reform}
    i\partial_tw+\Delta w+f^{\prime}(U_k) \cdot w=-\sum_{l=1}^{4}G_l.
\end{align}
Plugging \eqref{w-wtej-ej-def} into \eqref{equa-w-reform}
and using  algebraic computations one has
the equation of $e_k$ below.

\begin{lemma} \label{Lem-equa-ej}
For every $1\leq k\leq K$,
$e_k$ satisfies the equation
\begin{align} \label{equa-ej}
      i\lbb_k^2\partial_t e_k+\Delta e_k-e_k+(1+\frac{2}{d})Q^{\frac{4}{d}} e_k
        +\frac{2}{d}Q^{\frac{4}{d}} \ol {e_k}
    = -\sum_{l=1}^{4} \calh_l
      + \calo\((\<y\>^2 |\wt e_k|+  \<y\> |\na \wt e_k|) Mod_k\),
\end{align}
where
\begin{align} \label{Hl-Gl}
   \calh_l(t,y)
   = \lbb_k^{2+\frac d2} e^{-i\theta_k} e^{-i(\beta_k\cdot y - \frac{1}{4}\g_k|y|^2)} G_l(t,\lbb_k y + \a_k),
   \ \ 1\leq l\leq 4.
\end{align}
\end{lemma}

The error terms $\{\calh_l\}$ in \eqref{Hl-Gl} can be controlled by Lemma \ref{Lem-Hl} below.

\begin{lemma} \label{Lem-Hl}
Let $\calk$ belong to the generalized kernels of the linearized operator $L$
given by \eqref{L+-L-} below,
i.e., $\calk\in\{Q, yQ, |y|^2 Q, \na Q, \Lambda Q, \rho\}$.
Then, there exist $C, \delta >0$ such that
\begin{align}
   & \int | \calh_1(t,y)||\calk(y)| dy\leq C(T-t)^{4+\zeta} \wt D(t),  \label{H1-E}   \\
   & \int ( |\calh_2(t,y)|+|\calh_3(t,y)|) |\calk(y)|dy \leq Ce^{-\frac{\delta}{T-t}}\|w\|_{L^2},  \label{H2H3-E}  \\
   & \bigg|\int \calh_4(t,y)\calk(y)dy \bigg| \leq C(T-t)^{\upsilon_*+1}\|w\|_{L^2}, \label{H4-E}
\end{align}
where $\upsilon_*$ is the flatness index of the spatial functions $\{\phi_l\}$ in Hypothesis $(H1)$.
\end{lemma}

{\bf Proof.}
Estimates \eqref{H2H3-E} and \eqref{H4-E} were proved in \cite[(7,95), (7.96)]{SZ20},
hence we mainly focus on the estimate \eqref{H1-E}.

Define the renormalized variable $\ve_{R,k}$ by
\begin{align} \label{R-veRk}
   R(t,x) = \lbb_k^{-\frac d2} \ve_{R,k}(t,\frac{x-\a_k}{\lbb_k}) e^{i\theta_k}.
\end{align}

By \eqref{Hl-Gl},
\begin{align}
   \int |\calh_1(t,y)| |\calk(y) | dy
   \leq C (T-t)^{2-\frac d2} \int |G_1(t,x) \calk (\frac{x-\a_k}{\lbb_k})| dx.
\end{align}
Since  $\calk(y) \leq C e^{-\delta |y|}$
and by \eqref{G1-equaw},
\begin{align*}
   |G_1 |
   &\leq C \(|U|^{\frac 4d -1} + |z+R|^{\frac 4d -1}\) |z+R| |w|
          + C \(|U|^{\frac 4d -1} + |z+R|^{\frac 4d -1} + |w|^{\frac 4d -1}\) |w|^2,
\end{align*}
taking into account Lemma \ref{Lem-decoup-U} we derive
\begin{align}
    \int |\calh_1(t,y)| |\calk(y) | dy
    \leq& C \int e^{-\delta |y|}
               \bigg( \(e^{-\delta |y|} + |\ve_{z,k}|^{\frac 4d-1} + |\ve_{R,k}|^{\frac 4d-1}\)
               |\ve_{z,k}+\ve_{R,k}| |\wt e_k|  \nonumber \\
        & \qquad \qquad \  + \( e^{-\delta |y|} + |\ve_{z,k}|^{\frac 4d-1}  + |\ve_{R,k}|^{\frac 4d-1} + |\wt e_k|^{\frac 4d-1}\)
              |\wt e_k|^2 \bigg) dy
         + C e^{-\frac{\delta}{T-t}}.
\end{align}
Then, by \eqref{z-LtHm-a*}, \eqref{eyvez-bdd}, \eqref{R-Tt-Uniq} and \eqref{wtD-Tt-Uniq},
the R.H.S. can be bounded by, up to a universal constant,
\begin{align*}
   & \(\|e^{-\delta|y|}\ve_{z,k}\|_{L^\9} + \|\ve_{R,k}\|_{L^2} +  \|e^{-\delta|y|}|\ve_{z,k}|^{\frac 4d-1}\|_{L^\9}\|\ve_{z,k}+\ve_{R,k}\|_{L^2} \) \|\wt e_k\|_{L^2}
   + \|\ve_{R,k}\|_{H^1}^{\frac 4d-1} \|\ve_{z,k}+\ve_{R,k}\|_{H^1} \|\wt e_k\|_{H^1}    \nonumber \\
   &+ (1+\|e^{-\delta|y|}|\ve_{z,k}|^{\frac 4d-1}\|_{L^\9}) \|\wt e_k\|_{L^2}^2
   + \|\ve_{R,k}\|_{H^1}^{\frac 4d-1}\|\wt e_k\|_{H^1}^2
   + \|\wt e_k\|_{H^1}^{\frac 4d+1} + e^{-\frac{\delta}{T-t}} \nonumber \\
  \leq& \(\a^*(T-t)^{m+1+\frac d2} + (T-t)^{\kappa+1}\) \wt D
        + (T-t)^{\kappa(\frac 4d-1)} \wt D
        + \wt D^2 + (T-t)^{\kappa(\frac 4d-1)} \wt D^2
        + \wt D^{\frac 4d+1}  + e^{-\frac{\delta}{T-t}}   \nonumber \\
  \leq& (T-t)^{4+\zeta} \wt D.
\end{align*}
This yields \eqref{H1-E}
and finishes the proof. \hfill $\square$

Applying Lemmas \ref{Lem-equa-ej} and \ref{Lem-Hl}
and using algebraic identities in \eqref{Q-kernel}
one has the following ODE system of the renormalized variable $e_k$
along the six directions in the null space.

\begin{proposition} \label{Prop-ej-Lkernel}
Let $e_k$ be as in \eqref{w-wtej-ej-def}
and $e_{k,1} :={\rm Re} e_{k}$,
$e_{k,2} := {\rm Im} e_{k}$,
Then, for every $1\leq k\leq K$,
\begin{align}
&\frac{d}{dt} \<e_{k,1}, Q\>
    =  \calo((T-t)^{3+\zeta}  \sqrt{\wt N}), \label{Q-ej} \\
&\frac{d}{dt}  \<e_{k,2}, \Lambda Q\>
     =  2\lbb_k^{-2}  \<e_{k,1}, Q\>
       + \calo((T-t)^{3+\zeta} \sqrt{\wt N}),  \label{LaQ-ej} \\
&\frac{d}{dt} \<e_{k,1}, |y|^2 Q\>
    = -4\lbb_k^{-2}  \<e_{j,2}, \Lambda Q\>
        + \calo((T-t)^{3+\zeta} \sqrt{\wt N}), \label{y2Q-ej} \\
&\frac{d}{dt} \<e_{k,2}, \rho\>
   = \lbb_k^{-2} \<e_{k,1}, |y|^2 Q\>
      + \calo((T-t)^{3+\zeta} \sqrt{\wt N}), \label{rho-ej} \\
&\frac{d}{dt}  \<e_{k,2}, \na Q\>
     =    \calo((T-t)^{3+\zeta} \sqrt{\wt N}),   \label{naQ-ej}  \\
&\frac{d}{dt}  \<e_{k,1}, yQ\>
     =  -2\lbb_k^{-2} \<e_{k,2}, \na Q\>
         + \calo((T-t)^{3+\zeta}\sqrt{\wt N}), \label{yQ-ej}
\end{align}
\end{proposition}

{\bf Proof.}
By \eqref{equa-ej},
\begin{align}
   \frac{d}{dt} \<e_{k,1}, Q\>
   =& - \lbb_k^{-2} {\rm Im}
     \int Q \(\(\Delta e_k - e_k + (1+\frac 2d)Q^{\frac 4d}e_k + \frac 2d Q^{\frac 4d} \ol{e_k}\) + \sum\limits_{l=1}^4 \calh_l\) dy \nonumber \\
    & + \calo\( \lbb_k^{-2} Mod_k \int Q (\<y\>^2|\wt e_k| + \<y\> |\na \wt e_k|) dy \) .
\end{align}

Note that, by the definition of $L_{-}$
and the identity $L_{-}Q=0$ in \eqref{Q-kernel},
\begin{align*}
   {\rm Im} \int Q  \(\Delta e_k - e_k + (1+\frac 2d)Q^{\frac 4d}e_k + \frac 2d Q^{\frac 4d} \ol{e_k}\) dy
   =- {\rm Im} \int Q L_{-}e_{k,2} dy
   = - {\rm Im} \int L_{-}Q e_{k,2} dy =0.
\end{align*}

Moreover, since $\wt D \leq C (T-t) \sqrt{\wt N}$,
by Lemma \ref{Lem-Hl},
\begin{align*}
   \lbb_k^{-2} \bigg|\int Q \calh_l dy\bigg|
   \leq C \lbb_k^{-2} (T-t)^{4+\zeta} \wt D
   \leq C (T-t)^{3+\zeta} \sqrt{\wt N}.
\end{align*}

It also follows from \eqref{Mod-bdd-Uniq} that for $\kappa \geq  4$,
\begin{align*}
   \lbb_k^{-2} Mod_k \int Q \(\<y\>^2|\wt e_k| + \<y\> |\na \wt e_k|\) dy
   \leq& C \lbb_k^{-2} Mod \wt  D   \nonumber \\
   \leq& C \a^* (T-t)^\kappa \sqrt{\wt N}
   \leq C \a^* (T-t)^{3+\zeta}  \sqrt{\wt N}.
\end{align*}

Hence, \eqref{Q-ej} follows from the above estimates.
The proof of \eqref{LaQ-ej}-\eqref{yQ-ej} is similar,
see also the proof of \cite[Proposition 7.12]{SZ20}.
\hfill $\square$

As a consequence, we have the control of scalar $Scal_k$ below.
The proof is similar to that of \cite[Theorem 7.7]{SZ20}
and hence is omitted here.

\begin{theorem} (Control of $Scal_k$)\label{Thm-Scal}
There exists $C>0$
such that for $t$ close to $T$ and $1\leq k \leq K$,
\begin{align} \label{Scalj-N}
Scal_k(t)\leq C(T-t)^{2+\zeta} \wt N(t).
\end{align}
\end{theorem}

\subsection{Proof of conditional uniqueness}

We are now in position to prove the conditional uniqueness part in Theorem \ref{Thm-BW-RNLS}.

Let $\ve$ be a sufficiently small constant to be specified later
and let $t$ close to $T$ such that \eqref{a-ve-t-small} holds.
By Corollary \ref{Cor-Nt-Scal},  for any $t\in [t_*,T)$,
\begin{align} \label{D-iter}
   \wt N(t)
\leq C_1 \sum_{k=1}^{K} \sup_{t\leq s<T} \frac{Scal_k(s)}{\lbb_k^2(s)}
+C_1 \int_{t}^{T} \(\sum_{k=1}^{K}\frac{Scal_k(s)}{\lbb_k^3(s)}
   +  \ve  \frac{\wt N(s)}{T-s} \) ds,
\end{align}
which along with Theorem \ref{Thm-Scal}
yields that for some $\zeta>0$,
\begin{align} \label{D-iter*}
    \wt N(t)
\leq C_2 (T-t)^{\zeta} \wt N(t)
      +  C_2  \ve  \int_t^T \frac{\wt N(s)}{T-s} ds,
\end{align}
where $C_2$ is independent of $\ve$ and $t$.
Then, taking $t$ even closer to $T$
such that $C_2 (T-t)^{\zeta} \leq \frac 12$
we obtain the Gronwall type inequality
\begin{align} \label{D-iter.0}
   \wt N(t)
\leq 2 C_2  \ve \int_t^T \frac{ \wt N(s)}{T-s} ds.
\end{align}

Moreover, by \eqref{wtD-Tt-Uniq} and \eqref{wtN-def},
\begin{align} \label{wtN-Tt-1}
   \wt N(t) \leq C_3 (T-t)^{6+\zeta},
\end{align}
where $C_3 (\geq 1)$ is independent of $\ve$ and $t$.

We claim that for any $t$ close to $T$ and for any $l\geq 1$,
\begin{align} \label{wtN-Tt-induct}
   \wt N(t) \leq \(\frac{2C_2 C_3 \ve}{6+\zeta}\)^l (T-t)^{6+\zeta}.
\end{align}

To this end,
plugging \eqref{wtN-Tt-1}  into the Gronwall type inequality \eqref{D-iter.0} we get
\begin{align}  \label{D-iter.1}
    \wt N(t)
   \leq 2 C_2 \ve \int_t^T C_3 (T-s)^{5+\zeta} ds
   \leq \(\frac{2C_2C_3\ve}{6+\zeta}\) (T-t)^{6+\zeta},
\end{align}
which verifies \eqref{wtN-Tt-induct} at the preliminary step $l=1$.
Moreover,
plugging \eqref{wtN-Tt-induct} into \eqref{D-iter.0}
we derive that \eqref{wtN-Tt-induct} is still valid with $l+1$ replacing $l$.
Thus, the induction arguments lead to \eqref{wtN-Tt-induct}.

Therefore, take $\ve$ small enough such that  $\frac{2 C_2C_3 \ve}{6+\zeta} <1$.
Then, it follows from \eqref{wtN-Tt-induct} that
\begin{align} \label{D-iter.k}
  \wt  N(t) \leq \lim\limits_{l\to \9}  \(\frac{2 C_2C_3 \ve}{6+\zeta}\)^l (T-t)^{6+\zeta}
         = 0,
\end{align}
which yields $\wt N (t) = 0$  for $t$ close to $T$,
and so $w\equiv 0$.
The proof  of Theorem \ref{Thm-BW-RNLS} is complete.
\hfill $\square$

\section{Appendix} \label{Sec-App}

This Appendix mainly contains
preliminaries of linearized operators around the ground state,
the expansion of the nonlinearity
and the proof of Theorem \ref{Thm-Mod-bdd}.

\paragraph{{\bf Coercivity of the linearized operators}}

Let $L=(L_+,L_-)$ be the linearized operator around the ground state,
defined by
\begin{align} \label{L+-L-}
     L_{+}:= -\Delta + I -(1+{\frac{4}{d}})Q^{\frac{4}{d}}, \ \
    L_{-}:= -\Delta +I -Q^{\frac{4}{d}}.
\end{align}
The generalized null space of operator $L$ is
spanned by $\{Q, xQ, |x|^2 Q, \na Q, \Lambda Q, \rho\}$,
where
$\Lambda := \frac{d}{2}I_d + x\cdot \na$,
and $\rho$ is the unique $H^1$ spherically symmetric solution to the equation
\begin{align} \label{def-rho}
L_{+}\rho= - |x|^2Q,
\end{align}
which satisfies the exponential decay property (see, e.g., \cite{K-M-R, MP18}),
i.e., for some $C,\delta>0$,
\ben
|\rho(x)|+|\nabla \rho(x)|
\leq Ce^{-\delta|x|}.
\enn
Moreover, it holds that (see, e.g., \cite[(B.1), (B.10), (B.15)]{W85})
\be \ba \label{Q-kernel}
&L_+ \na Q =0,\ \ L_+ \Lambda Q = -2 Q,\ \ L_+ \rho = -|x|^2 Q, \\
&L_{-} Q =0,\ \ L_{-} xQ = -2 \na Q,\ \ L_{-} |x|^2 Q = - 4 \Lambda Q.
\ea\ee

Lemma \ref{Lem-coer-f-local} below contains the key localized coercivity of the linearized operator.

\begin{lemma}(Localized coercivity \cite[Corollary 3.4]{SZ20})   \label{Lem-coer-f-local}
Let $\phi$ be a positive smooth radial function on $\R^d$,
such that
$\phi(x) = 1$ for $|x|\leq 1$,
$\phi(x) = e^{-|x|}$ for $|x|\geq 2$,
$0<\phi \leq 1$,
and $\lf|\frac{\nabla\phi}{\phi}\rt|\leq C$ for some $C>0$.
Set $\phi_A(x) :=\phi\lf(\frac{x}{A}\rt)$, $A>0$.
Then,
for $A$ large enough we have
\begin{align} \label{coer-f-local}
\int (|f|^2+|\nabla f|^2)\phi_A -(1+\frac 4d)Q^{\frac4d}f_1^2-Q^{\frac4d}f_2^2dx\geq C_1\int(|\nabla f|^2+|f|^2)\phi_A dx-C_2Scal(f),
\end{align}
where $C_1, C_2>0$,  $f_1, f_2$ are the real and imaginary parts of $f$, respectively,
and $Scal(f)$ denotes the scalar products along  the unstable directions
in the null space
\begin{align} \label{Scal-def}
Scal(f) :=\<f_1,Q\>^2+\<f_1,xQ\>^2+\<f_1,|x|^2Q\>^2+\<f_2,\nabla Q\>^2+\<f_2,\Lambda Q\>^2+\<f_2,\rho\>^2,
\end{align}
\end{lemma}

\paragraph{\bf Expansion of the nonlinearity.}

Let us recall the expansion that
for any continuous differentiable function $g:\mathbb{C}\to \mathbb{C}$
and for any $v, w\in \mathbb{C}$,
(see, e.g., \cite[(3.10)]{KV13})
\begin{align} \label{g-gz-expan}
   g(v+w) = g(v) + g'(v, w)\cdot w
\end{align}
with
\begin{align}
  g'(v,w) \cdot w
     :=&  w  \int_0^1 \partial_z g(v+sw) ds
           + \ol w  \int_0^1 \partial_{\ol z} g(v+sw) ds,  \label{f'R}
\end{align}
where $z=x+iy \in \mathbb{C}$,  $\partial_z g$ and $\partial_{\ol{z}} g$
are the usual complex derivatives $\partial_z g= \frac 12(\partial_x g - i \partial_y g)$,
$\partial_{\ol z} g= \frac 12(\partial_x g + i \partial_y g)$,
respectively.
Moreover,
if  $\partial_z g$ and $\partial_{\ol{z}} g$  are also continuously differentiable,
we may expand $g$ up to the second order
\begin{align} \label{g-gzz-expan}
   g(v+w)
   =& g(v)
     + g'(v)\cdot w
     + g''(v, w) \cdot w^2,
\end{align}
where
\begin{align}
   & g'(v)\cdot w
     :=  \partial_z g(v) w
     + \partial_{\ol{z}} g(v)  \ol w, \nonumber \\
   & g''(v,w) \cdot w^2
     :=   w^2 \int_0^1 t  \int_0^1 \partial_{zz}g(v+st w) dsdt
      + 2|w|^2 \int_0^1 t \int_0^1 \partial_{z \ol z}g(v+st w) dsdt \nonumber \\
      & \qquad \qquad \qquad \ \ + \ol w^2 \int_0^1 t  \int_0^1 \partial_{\ol z \ol z}g(v+st w) dsdt,    \label{f''R2}
\end{align}

In particular,
for $f(z) :=|z|^{\frac 4d} z$ with $d=1,2$, $z\in \mathbb{C}$,
one has
\begin{align} \label{f-Taylor}
   f(v+w) =&  f(v) + f'(v)\cdot R
         +  f''(v)\cdot w^2
         + \calo\(\sum\limits_{l=3}^{1+\frac 4d} |v|^{1+\frac 4d-l} |w|^l\),
\end{align}
where
\begin{align}
   f'(v)\cdot w :=&  \partial_z f(v) w+ \partial_{\ol{z}} f(v) \ol{w}
                 = (1+\frac 2d) |v|^{\frac 4d} w + \frac 2d |v|^{\frac 4d-2} v^2 \ol{w},   \label{f-linear} \\
   f''(v)\cdot w^2
                :=& \frac 12 \partial_{zz}f(v) w^2 + \partial_{z\ol{z}} f(v) |w|^2 + \frac 12 \partial_{\ol{z}\ol{z}} f(v) \ol{w}^2 \nonumber \\
                 =& \frac 1d(1+\frac 2d)|v|^{\frac 4d -2} \ol{v} w^2
                    + \frac 2d (1+\frac 2d) |v|^{\frac 4d -2} v |w|^2
                    +   \frac 1d(\frac 2d-1) |v|^{\frac 4d -4} v^3 \ol{w}^2.  \label{f-quadratic}
\end{align}
The following estimates are also useful:
\begin{align}
   & |f(v_1) - f(v_2)| \leq C (|v_1|^{\frac 4d} + |v_2|^{\frac 4d})|v_1-v_2|, \label{fv1-fv2} \\
   & |f'(v_1)\cdot w - f'(v_2)\cdot w| \leq C (|v_1|^{\frac 4d-1}+|v_2|^{\frac 4d-1})|v_1-v_2||w|,   \label{f'v1-f'v2}  \\
   & |f''(v_1,w)\cdot w^2 - f''(v_2,w)\cdot w^2| \leq C (|v_1|^{\frac 4d-2} + |v_2|^{\frac 4d-2} + |w|^{\frac 4d-2}) |v_1-v_2||w|^2, \label{f''v1-f''v2} \\
   & |f''(v,w)\cdot w^2| \leq C (|v|^{\frac 4d-1}+|w|^{\frac 4d-1})|w|^2.  \label{f''-bdd}
\end{align}

\paragraph{\bf Proof of Theorem \ref{Thm-Mod-bdd}.}

We adapt the arguments as in \cite{SZ19,SZ20}.

$(i)$ {\it Reformulation of the equation of remainder.}
By \eqref{u-dec} and \eqref{g-gzz-expan},
\begin{align} \label{fv-fUz-expan-2}
   f(v) = f(U+z) + f'(U+z)\cdot R + f''(U+z, R) \cdot R^2.
\end{align}
Plugging this into \eqref{equa-R} leads to the equation of $R$:
\begin{align}\label{eq-U-R}
   i\partial_tR
   +\sum_{k=1}^{K}
   \(\Delta R_{k}+(1+\frac 2d)|U_{k}|^{\frac{4}{d}}R_{k}
   +\frac 2d|U_{k}|^{\frac{4}{d}-2}U_{k}^2 \ol{R_{k}}
   +  i\partial_tU_{k}+\Delta U_{k}+|U_{k}|^{\frac{4}{d}}U_{k}\)
  = -\sum\limits_{l=1}^5 H_l,
\end{align}
where $H_1, H_2$ contain the interactions between different blow-up profiles $U_j$ and $U_l$, $j\not =l$,
\begin{align}
  H_1 : =&  f'(U)\cdot R  - \sum\limits_{k=1}^K f'(U_k)\cdot R_k, \label{H1-def} \\
  H_2 : =&  f(U) - \sum\limits_{k=1}^K f(U_k), \label{H2-def}
\end{align}
the terms $H_3, H_4$ contain the regular flow $z$, i.e.,
\begin{align}
  H_3 :=&  f'(U + z)\cdot R - f'(U) \cdot R  + f''(U+z,R)\cdot R^2, \label{H3-def} \\
  H_4 :=& f(U+z) - f(U) - f(z), \label{H4-def}
\end{align}
and the lower order perturbations are contained in $H_5$:
\begin{align} \label{H5-def}
   H_5:= \sum\limits_{l=1}^K \(a_1\cdot\nabla  (U_l+R_l)  +a_0  (U_l+R_l) \),
\end{align}
where $a_1, a_0$ are the coefficients of lower order perturbations given by
\eqref{a1-loworder} and \eqref{a0-loworder}, respectively.

$(ii)$ {\it Estimate of Modulation equations.}
Let us take the modulation equation $\lbb^2_k \dot{\g}_k + \g_k^2$
to illustrate the main arguments below.
As $R(T_*)=0$,
we may take $t^*$ close to $T$ such that
$\|R\|_{C([t^*, T_*];H^1)}  \leq 1$.

Taking the inner product of \eqref{eq-U-R} with
$\Lambda_k {U_{k}}$ and then taking the real part
we get
\begin{align}\label{eq-Ul-Rl}
   &-{\rm Im}\langle\partial_tR,\Lambda U_{k}\rangle
   +{\rm Re}\langle\Delta R_{k}+(1+\frac 2d)|U_{k}|^{\frac{4}{d}}R_{k}+\frac 2d|U_{k}|^{\frac{4}{d}-2}U_{k}^2\ol{{R}_{k}},\Lambda_k U_{k}\rangle  \nonumber \\
   &+{\rm Re}\langle i\partial_tU_{k}+\Delta U_{k}+|U_{k}|^{\frac{4}{d}}U_{k},\Lambda_k U_{k}\rangle  \nonumber \\
 = &-{\rm Re} \<\sum_{j\neq k} (\Delta R_{j}+(1+\frac 2d)|U_{j}|^{\frac{4}{d}}R_{j}+\frac 2d|U_{j}|^{\frac{4}{d}-2}U_{j}^2\ol{{R}_{j}}),\Lambda_k U_{k} \>   \nonumber \\
   & -{\rm Re} \< \sum_{j\neq k}( i\partial_tU_{j}+\Delta U_{j}+|U_{j}|^{\frac{4}{d}}U_{j}),\Lambda_k U_{k} \>
    -\sum\limits_{l=1}^5 {\rm Re}\langle H_l,\Lambda_k U_{k}\rangle.
\end{align}

First for the L.H.S. of \eqref{eq-U-R},
we have (see the proof of \cite[(4.38)]{SZ20}, \cite[(6.43)]{CSZ21})
\begin{align} \label{Mod-LHS-bdd}
   \lbb_k^2 \times ({\rm L.H.S.\ of}\ \eqref{eq-Ul-Rl})
  =&- \frac{1}{4}\|yQ\|_2^2(\lbb_k^2\dot{\g_k}+\g_k^2)+  M_k   \nonumber \\
  &+ \calo\( (P+\|R\|_{L^2} + e^{-\frac{\delta}{T-t}}) Mod + P^2 \|R\|_{L^2} + \|R\|_{L^2}^2 + e^{-\frac{\delta}{T-t}} \).
\end{align}

Next we show that
the R.H.S. of \eqref{eq-Ul-Rl}
contribute acceptable orders.
This is mainly due to the exponentially small interactions
between different blow-up profiles
and to the flatness of both the regular profile $z$ and lower order coefficients $a_1, a_0$
at the singularities.

To be precise,
in view of Lemma \ref{Lem-decoup-U} and \eqref{equa-Ut},
we have that for some $\delta>0$,
\begin{align}
  & \bigg| \<  \sum_{j\neq k} \(\Delta R_{j}+(1+\frac 2d)|U_{j}|^{\frac{4}{d}}R_{j}+\frac 2d|U_{j}|^{\frac{4}{d}-2}U_{j}^2 \ol{{R}_{j}} \)
              + H_1,\Lambda_k U_{k} \> \bigg|
\leq C \lbb_k^{-2}e^{-\frac{\delta}{T-t}}\|R\|_{L^2},  \label{RHS-Rj-Ul} \\
   & \bigg| \<\sum_{j\neq k} \( i\partial_tU_{j}+\Delta U_{j}+|U_{j}|^{\frac{4}{d}}U_{j} \) + H_2,\Lambda_k U_{k} \> \bigg|
\leq C  \lbb_k^{-2}e^{-\frac{\delta}{T-t}}(1+Mod)  ,  \label{RHS-Uj-Ul}
\end{align}

For the third term $H_3$,
by \eqref{f'v1-f'v2}
and Lemma \ref{Lem-decoup-U},
\begin{align} \label{RHS-f'Uz-LamU}
       & \bigg| {\rm Re} \< f'(U + z)\cdot R - f'(U) \cdot R, \Lambda_k U_k\>  \bigg| \nonumber \\
   \leq& C \( \int (|U|^{\frac 4d-1} + |z|^{\frac 4d-1}) |z||R| \Lambda_k U_k| dx
   +  e^{-\frac{\delta}{T-t}} \) \nonumber \\
   \leq&  C \(\lbb_k^{-\frac d2 (\frac 4d+1)} \lbb_k^{\frac d2} \|e^{-\delta|y|}\ve_{z,k}\|_{L^\9} \|R\|_{L^2}
          + \|z\|_{L^\9}^{\frac 4d-1}  \|e^{-\delta|y|}\ve_{z,k}\|_{L^\9} \|R\|_{L^2}
          + e^{-\frac{\delta}{T-t}} \)  \nonumber \\
   \leq& C \(   \lbb_k^{-2} \|e^{-\delta|y|}\ve_{z,k}\|_{L^\9} D
             + e^{-\frac{\delta}{T-t}}\).
\end{align}

Moreover, by \eqref{f''-bdd} and \eqref{R-Lp-D},
\begin{align} \label{RHS-Nf}
|\<f''(U+z,R)\cdot R^2, \Lambda_k U_{k} \>|
  \leq& C  \int \(|U|^{\frac 4d-1} + |z|^{\frac 4d-1} + |R|^{\frac 4d-1}\) |R|^2 |\Lambda_k U_k| dx   \nonumber \\
  \leq& C\( \lbb_k^{-2}\|R\|_{L^2}^2
              + \lbb_k^{-2} \|\ve_{z,k}\|_{L^\9}^{\frac 4d-2} \|e^{-\delta|y|} \ve_{z,k}\|_{L^\9} \|R\|_{L^2}^2
              + \lbb_k^{-\frac d2} \|R\|_{L^{1+\frac 4d}}^{1+\frac 4d} + e^{-\frac{\delta}{T-t}} \)   \nonumber \\
  \leq& C \( \lbb_k^{-2} D^2  + e^{-\frac{\delta}{T-t}} \).
\end{align}
Hence, we conclude from \eqref{RHS-f'Uz-LamU} and \eqref{RHS-Nf} that
\begin{align} \label{RHS-H3}
   {\rm Re} \<H_3, \Lambda_k U_k\>
   = \calo\(\lbb_k^{-2} \|e^{-\delta|y|} \ve_{z,k}\|_{L^\9} D  + \lbb_k^{-2} D^2 + e^{-\frac{\delta}{T-t}} \).
\end{align}

We also see that
\begin{align} \label{RHS-H4}
  |{\rm Re} \<H_4, \Lambda_k U_k\>|
  \leq  C \sum\limits_{j=1}^{4/d} \int |U|^{1+\frac 4d-j} |z|^j |\Lambda_k U_k| dx
  \leq  C \lbb_k^{-2}  \|e^{-\delta|y|} \ve_{z,k}\|_{L^\9}.
\end{align}

Regarding the $H_5$ term
on the R.H.S. of  \eqref{eq-Ul-Rl},
by Lemma \ref{Lem-decoup-U}, the change of variables and integrating by parts formula,
\begin{align} \label{bc-wtbwtc}
   {\rm Re} \<H_5, \Lambda_k U_k\>
   = & {\rm Re} \< \lbb_{k}^{-1} \wt {a}_{1,k}\cdot\nabla (Q_{k}+\varepsilon_{k})
+ \wt {a}_{0,k} (Q_{k}+\varepsilon_{k}),\Lambda Q_{k} \>
       + \calo(e^{-\frac{\delta}{T-t}}) \nonumber \\
   =&  - \lbb_k^{-1} {\rm Re} \< {\rm div} \wt {a}_{1,k}\   (Q_{k}+\varepsilon_{k}), \Lambda Q_k\>
      - \lbb_k^{-1} {\rm Re} \< Q_{k}+\varepsilon_{k}, \ol{\wt {a}_{1,k}} \cdot \na (\Lambda Q_k)\>  \nonumber \\
    &  + {\rm Re} \< \wt {a}_{0,k}(Q_{k}+\varepsilon_{k}), \Lambda Q_{k} \>
       + \calo(e^{-\frac{\delta}{T-t}})
\end{align}
where $\wt {a}_{1,k}$ and $\wt {a}_{0,k}$ are defined as in Lemma \ref{Lem-wtbwtc}.
Then, applying Lemma \ref{Lem-wtbwtc}
we obtain
\begin{align} \label{RHS-b-c}
    | {\rm Re} \<H_5, \Lambda_k U_k\> |
    \leq  C \(\lbb_k^{-2} P^{\upsilon_*+1} + e^{-\frac{\delta}{T-t}}\).
\end{align}

Hence, it follows from estimates \eqref{RHS-Rj-Ul}, \eqref{RHS-Uj-Ul},
\eqref{RHS-H3}, \eqref{RHS-H4} and \eqref{RHS-b-c} that
\begin{align} \label{Mod-RHS-bdd}
   {\rm R.H.S.\ of}\ \eqref{eq-Ul-Rl}
   \leq& C  \lbb_k^{-2}
            \( e^{-\frac{\delta}{T-t}} Mod
               + D^2
               + \|e^{-\delta|y|} \ve_{z,k}\|_{L^\9}
               +P^{\upsilon_*+1}
               + e^{-\frac{\delta}{T-t}} \).
\end{align}

Now, combining \eqref{Mod-LHS-bdd} and \eqref{Mod-RHS-bdd} together
we conclude that for each $1\leq k\leq K$,
\begin{align}
|\lbb_k^2\dot{\g_k}+\g_k^2|
   \leq& C \bigg( (P + \|R\|_{L^2}+e^{-\frac{\delta}{T-t}} )Mod  +|M_k|
          + P^2 D+ D^2  \nonumber \\
      & \qquad  + \|e^{-\delta |y|} \ve_{z,k}\|_{L^\9}
          +P^{\upsilon_*+1} + e^{-\frac{\delta}{T-t}} \bigg).
\end{align}

Similar arguments apply to the remaining four modulation equations
$|{\lbb_k \dot{\a}_k}-2\beta_k|$,
$|\lbb_k\dot{\lbb}_k + \g_k|$,
$|\lbb^2_k\dot{\beta}_k + \beta_k \g_k|$
and $|\lbb_k^2\dot{\theta}_k - 1-|\beta_k|^2|$,
by taking the inner products of equation (\ref{eq-U-R})
with $i(x-\a_k) U_k$, $i|x-\a_k|^2 U_k$, $\nabla {U_k}$, $ \varrho_{k}$,
respectively, and
then taking the real parts.
This leads to
\begin{align}\label{Mod-bdd-l}
Mod_{k}(t)
\leq& C\bigg( (P+\|R\|_{L^2}+e^{-\frac{\delta}{T-t}}) Mod   +  |M_k|
    + P^2 D + D^2
      + \|e^{-\delta |y|} \ve_{z,k}\|_{L^\9}
           +P^{\upsilon_*+1} +  e^{-\frac{\delta}{T-t}} \bigg).
\end{align}

Therefore,
taking $t^*$ even closer to $T$ such that
$$(1+C) (P(t) + \|R(t)\|_{C([t^*, T_*];H^1)} +e^{- \frac{\delta}{T-t}}) \leq \frac 12$$
and then summing over $k$ and using \eqref{eyvez-bdd}
we obtain \eqref{Mod-bdd}.

$(iii)$ {\it Improved estimate of $\lbb_k \dot{\lbb}_k + \g_k$}.
Taking the inner product of equation \eqref{eq-U-R} with $|x-\a_k|^2 U_k$, then taking the imaginary part
and arguing as in the proof of \eqref{Mod-RHS-bdd}
we have that, similarly to \eqref{eq-Ul-Rl},
\begin{align} \label{inner-R-y2Uk}
   & {\rm Re}\langle\partial_tR, |x-\a_k|^2 U_k\rangle
   +{\rm Im}\langle\Delta R_{k}+(1+\frac 2d)|U_{k}|^{\frac{4}{d}}R_{k}+\frac 2d|U_{k}|^{\frac{4}{d}-2}U_{k}^2\ol{{R}_{k}}, |x-\a_k|^2 U_k\rangle  \nonumber \\
   &+{\rm Im}\langle i\partial_tU_{k}+\Delta U_{k}+|U_{k}|^{\frac{4}{d}}U_{k}, |x-\a_k|^2 U_k \rangle  \nonumber \\
  =& \calo\(    D^2
               + \|e^{-\delta|y|} \ve_{z,k}\|_{L^\9}
               +P^{\upsilon_*+1}
               + e^{-\frac{\delta}{T-t}}   \).
\end{align}
Note that,
the bound on the R.H.S. above equals to \eqref{Mod-RHS-bdd}
multiplied by $\lbb_k^2$,
which essentially relies on the exponential decay of ground state.
We also used the fact that,
by \eqref{Mod-bdd},
$Mod = \calo(1)$,
and thus $e^{-\frac{\delta}{T-t}} Mod =\calo(e^{-\frac{\delta}{T-t}} )$.

Regarding the L.H.S. of \eqref{inner-R-y2Uk},
by the orthogonality condition (\ref{ortho-cond-Rn-wn}), \eqref{equa-Ut} and Lemma \ref{Lem-decoup-U},
\begin{align}  \label{ptR-Ul}
{\rm Re}\langle\partial_t R,|x-\a_{k}|^2 U_{k}\rangle
&=2 \dot{\a}_{k}\cdot{\rm Re}\langle R,(x-\a_{k})U_{k}\rangle-{\rm Re}\langle R,|x-\a_{k}|^2\partial_t U_{k}\rangle  \nonumber \\
&= -{\rm Re}\langle R_{k},|x-\a_{k}|^2\partial_t U_{k}\rangle
      +\calo( e^{-\frac{\delta}{T-t}}\|R\|_{L^2}).
\end{align}
Then, using \eqref{Uj-Qj-Q}, \eqref{equa-Ut}, \eqref{Rj-ej} and
the algebraic identity
\begin{align} \label{equa-Qk}
  \Delta Q_k - Q_k + |Q_k|^\frac 4d Q_k
  = |\beta_k - \frac{\g_k}{2} y |^2 Q_k
    - i \g_k \Lambda Q_k
    + 2i \beta_k \cdot \na Q_k
\end{align}
we get
\begin{align*}
  -  {\rm Re}\langle R_{k\emph{}},|x-\a_{k}|^2\partial_t U_{k}\rangle
  =& - {\rm Im} \< \ve_k, |y|^2(\Delta Q_k + |Q_k|^\frac 4d Q_k)\>
     + \calo(Mod \|\ve_k\|_{L^2}) \\
  =& - {\rm Im} \< \ve_k, |y|^2 Q_k\>
     - \g_k {\rm Re} \<\ve_k, |y|^2 \Lambda Q_k\>
     + 2 \beta_k \cdot {\rm Re} \<\ve_k, |y|^2 \na Q_k\> \\
   &  + \calo\((Mod+P^2) D \),
\end{align*}
By the integration by parts formula and the almost orthogonality \eqref{Orth-almost},
\begin{align}
   & - \g_k {\rm Re} \<\ve_k, |y|^2 \Lambda Q_k\>
     + 2 \beta_k {\rm Re} \<\ve_k, |y|^2 \na Q_k\> \nonumber \\
  =& \g_k {\rm Re} \<\Lambda \ve_k, |y|^2 Q_k\>
     - 2 \beta_k {\rm Re} \<\na \ve_k, |y|^2 Q_k\>
     + \calo(e^{-\frac{\delta}{T-t}}\|R\|_{L^2}).
\end{align}
Thus, we obtain
\begin{align} \label{ptR-x2Uj}
      {\rm Re}\langle\partial_t R,|x-\a_{k}|^2 U_{k}\rangle
   =& - {\rm Im}  \< \ve_k, |y|^2 Q_k\>
     + \g_k {\rm Re} \<\Lambda \ve_k, |y|^2 Q_k\>
     - 2 \beta_k \cdot {\rm Re} \<\na \ve_k, |y|^2 Q_k\> \nonumber \\
     & + \calo\( (Mod+P^2+e^{-\frac{\delta}{T-t}}) D\).
\end{align}

Furthermore,
using  \eqref{equa-Ut}, the identities
\begin{align}
   & \Lambda  Q_k = \(\Lambda Q + i(\beta_k \cdot y - \frac 12 \g_k |y|^2) Q\)e^{i(\beta_k\cdot y - \frac 14 \g_k |y|^2)}, \label{LaQj-LaQ} \\
   & \na  Q_k = \(\na Q + i (\beta_k - \frac 12 \g_k y) Q\)e^{i(\beta_k\cdot y - \frac 14 \g_k |y|^2)}, \label{naQj-naQ}
\end{align}
and $\<\Lambda Q, |y|^2 Q\> = - \|yQ\|_{L^2}^2$
we compute
\begin{align} \label{ptUj-DUj-fUj}
      {\rm Im}\langle i\partial_tU_{k}+\Delta U_{k}+|U_{k}|^{\frac{4}{d}}U_{k},|x-\a_{k}|^2U_{k}\rangle
   = (\lbb_k \dot \lbb_k + \g_k) \|yQ\|_{L^2}^2.
\end{align}

Therefore,
plugging \eqref{ptR-x2Uj} and \eqref{ptUj-DUj-fUj}  into \eqref{inner-R-y2Uk}
and using the change of variables,
the bound $|M_k|\leq CD$ and \eqref{Mod-bdd}
we obtain the equation for the
renormalized variable  $\ve_{k}$ below
\begin{align}  \label{equa-venl}
&{\rm Im}\langle\Delta \varepsilon_{k}-\varepsilon_{k}+(1+\frac 2d)|Q_{k}|^{\frac{4}{d}}\varepsilon_{k}
+\frac 2d|Q_{k}|^{\frac{4}{d}-2}Q_{k}^2\ol{\varepsilon_k},|y|^2Q_{k}\rangle \nonumber \\
&+\gamma_{k}{\rm Re}\langle \Lambda \varepsilon_{k},|y|^2 Q_{k}\rangle
 - 2\beta_{k}\cdot{\rm Re}\langle \na \varepsilon_{k},|y|^2 Q_{k}\rangle
 +(\lambda_{k}\dot{\lambda}_{k}+\gamma_{k}) \|yQ\|_{L^2}^2 \nonumber \\
=& \mathcal{O}\big(P^2 D + D^2 + \|e^{-\delta|y|}\ve_{z,k}\|_{L^\9} + P^{\upsilon_*+1} + e^{-\frac{\delta}{T-t}} \big).
\end{align}

Writing the L.H.S. of \eqref{equa-venl} in terms of real and imaginary parts
we see that the first three terms are exactly the first line of \cite[(4.27)]{SZ19}
and hence are of order $\calo(P^2\|R\|_{L^2})$,
due to \cite[(4.28)]{SZ19}.
This yields that
\begin{align}
   {\rm L.H.S.\ of}\ \eqref{equa-venl}
   = (\lbb_k\dot{\lbb}_k + \g_k) \|yQ\|_{L^2}^2 + \calo(P^2 D).
\end{align}
Therefore,
plugging this into \eqref{equa-venl} and using \eqref{eyvez-bdd} we obtain
the desired estimate \eqref{lbb-ga-mod}.  \hfill $\square$

\section*{Acknowledgements}

M. R\"ockner and D. Zhang thank for the financial support
by the Deutsche Forschungsgemeinschaft
(DFG, German Science Foundation) through SFB 1283/2
2021-317210226 at Bielefeld University.
Y. Su is supported by NSFC (No. 11601482).
D. Zhang  is also grateful for the support by NSFC (No. 11871337)
and Shanghai Rising-Star Program 21QA1404500.


\begin{thebibliography}{99}

\bibitem{BCIR94}
O. Bang, P.L. Christiansen, F. If, K.O. Rasmussen,
Temperature effects in a nonlinear model of monolayer Scheibe aggregates.
{\it Phys. Rev. E} {\bf 49} (1994), 4627--4636.

\bibitem{BCIRG95}
O. Bang, P.L. Christiansen, F. If, K.O. Rasmussen, Y.B. Gaididei,
White noise in the two-dimensional nonlinear Schr\"odinger equation,
{\it Appl. Anal}. {\bf 57} (1995), no. 1-2, 3--15.



\bibitem{BRZ14}
V. Barbu, M. R\"ockner, D.Zhang,
The stochastic nonlinear Schr\"odinger equation with multiplicative noise:
the rescaling aproach, {\it J. Nonlinear Sciences}, {\bf 24} (2014), 383--409.

\bibitem{BRZ16}
V. Barbu, M. R\"ockner,  D. Zhang,
Stochastic nonlinear Schr\"odinger equations.
{\it Nonlinear Anal}. {\bf 136} (2016), 168--194.

\bibitem{BRZ16.1}
V. Barbu, M. R\"ockner,  D. Zhang,
The stochastic logarithmic Schr\"odinger equation.
{\it J. Math. Pures Appl.} {\bf (9)} 107 (2017), no. 2, 123-149.

\bibitem{BRZ18}
V. Barbu, M. R\"ockner, D. Zhang,
Optimal bilinear control of nonlinear stochastic Schr\"odinger equations driven by linear multiplicative noise.
{\it Ann. Probab}. {\bf 46} (2018), no. 4, 1957--1999.

\bibitem{BG09}
A. Barchielli, M. Gregoratti,
Quantum Trajectories and Measurements in Continuous Case.
The Diffusive Case,
{\it Lecture Notes Physics} {\bf 782}, Springer Verlag,
Berlin, 2009.

\bibitem{B12}
M. Beceanu,
A critical center-stable manifold for Schr\"odinger's equation in three dimensions.
{\it Comm. Pure Appl. Math.} {\bf 65} (2012), no. 4, 431--507.

\bibitem{BW97}
J. Bourgain, W. Wang,
Construction of blowup solutions for the nonlinear Schr\"odinger equation with critical nonlinearity.
Dedicated to Ennio De Giorgi. {\it Ann. Scuola Norm. Sup. Pisa Cl. Sci.} {\bf (4) 25} (1997), no. 1-2, 197--215.

\bibitem{BM14}
Z. Brz\'{e}zniak, A. Millet, On the stochastic Strichartz estimates and the stochastic nonlinear Schr\"{o}dinger equation on a compact Riemannian manifold.
{\it Potential Anal.} {\bf 41} (2014), no. 2, 269--315.



\bibitem{CSZ21}
D. Cao, Y. Su, D. Zhang,
On uniqueness of multi-bubble blow-up solutions and multi-solitons to $L^2$-critical nonlinear Schr\"odinger equations,
\texttt{arXiv:2105.14554}.

\bibitem{C03}
T. Cazenave,  Semilinear Schr\"odinger equations.
{\it Courant Lecture Notes in Mathematics}, 10. New York University,
Courant Institute of Mathematical Sciences,
 New York; American Mathematical Society, Providence, RI, 2003. xiv+323 pp.

\bibitem{C20}
T. Cazenave, An overview of the nonlinear Schr\"odinger equation. Lecture notes,
2020. https://www.ljll.math.upmc.fr/cazenave/.


\bibitem{Co11}
V. C\^{o}mbet,
Multi-soliton solutions for the supercritical gKdV equations.
{\it Comm. Partial Differential Equations}. {\bf 36} (2011), no. 3, 380--419.

\bibitem{C06}
R. C\^ote,
Construction of solutions to the subcritical gKdV equations with a given asymptotical behavior.
{\it J. Funct. Anal.} {\bf 241} (2006), no. 1, 143--211.

\bibitem{C07}
R. C\^ote,
Construction of solutions to the $L^2$-critical KdV equation with a given asymptotic behaviour.
{\it Duke Math. J.} {\bf 138} (2007), no. 3, 487--531.

\bibitem{CF20}
R. C\^{o}te, X. Friederich,
On smoothness and uniqueness of multi-solitons of the
non-linear Schr\"odinger equations.
2020. hal-02873307v2,
to appear in {\it Commun. Partial Differ. Equ.}.

\bibitem{CKLS18}
R. C\^ote, C. Kenig, A. Lawrie, W. Schlag,
Profiles for the radial focusing 4d energy-critical wave equation.
{\it Comm. Math. Phys.} {\bf 357} (2018), no. 3, 943--1008.

\bibitem{CL11}
R. C\^ote, S. Le Coz,
High-speed excited multi-solitons in nonlinear Schr\"odinger equations.
{\it J. Math. Pures Appl.} {\bf (9)} 96 (2011), no. 2, 135--166.

\bibitem{CMM11}
R. C\^ote, Y. Martel, F. Merle,
Construction of multi-soliton solutions for the $L^2$-supercritical gKdV and NLS equations.
{\it Rev. Mat. Iberoam.} {\bf 27} (2011), no. 1, 273--302.

\bibitem{BD02}
A. de Bouard, A. Debussche,
On the effect of a noise on the solutions of the focusing supercritical nonlinear Schr\"odinger equation.
{\it Probab. Theory Related Fields} {\bf 123} (2002), no. 1, 76--96.

\bibitem{BD03}
A. de Bouard, A. Debusche,
The stochastic nonlinear Schr\"odinger equation in $H^1$,
{\it Stoch. Anal. Appl.} {\bf 21} (2003), 97--126.

\bibitem{BD05}
A. de Bouard, A. Debussche,
Blow-up for the stochastic nonlinear Schr\"odinger equation with multiplicative noise.
{\it Ann. Probab.} {\bf 33} (2005), no. 3, 1078--1110.

\bibitem{BDM02}
A. de Bouard, A. Debussche, L.D. Menza,
Theoretical and numerical aspects of stochastic nonlinear Schr\"{o}dinger equations,
{\it Journ\'{e}es "\'{E}quations aux D\'{e}riv\'{e}es Partielles"}
(Plestin-les-Gr\`{e}ves, 2001),
Exp. No. III, 13 pp., Univ. Nantes, Nantes, 2001.

\bibitem{DL02}
A. Debussche, L.D. Menza,
Numerical simulation of focusing stochastic nonlinear Schr\"{o}dinger equations,
{\it Phys. D} {\bf 162} (2002), no. 3-4, 131--154.

\bibitem{DL02.2}
A. Debussche, L.D. Menza,
Numerical resolution of stochastic focusing NLS equations.
{\it Appl. Math. Lett.} {\bf 15} (2002), no. 6, 661--669.

\bibitem{D15}
B. Dodson,
Global well-posedness and scattering for the mass critical nonlinear Schr\"odinger equation
with mass below the mass of the ground state.
{\it Adv. Math.} {\bf 285} (2015), 1589--1618.


\bibitem{DJKM17}
T. Duyckaerts, H. Jia, C. Kenig, F. Merle,
Soliton resolution along a sequence of times for the focusing energy critical wave
equation.
{\it Geom. Funct. Anal.} {\bf 27} (2017), no. 4, 798--862.

\bibitem{DKM13}
T. Duyckaerts, C. Kenig,  F. Merle,
Classification of radial solutions of the focusing, energy-critical wave equation.
{\it Camb. J. Math.} {\bf 1} (2013), no. 1, 75--144.

\bibitem{DKM19}
T. Duyckaerts, C. Kenig, F. Merle,
Soliton resolutin for the radial critical wave equation in all odd space dimensions.
\texttt{arXiv:1912.07664v1}.


\bibitem{DNPZ92}
S. Dyachenko, A.C. Newell, A. Pushkarev, V.E. Zakharov,
Optical turbulence: weak turbulence, condensates and collapsing filaments in the nonlinear Schr\"odinger equation.
{\it Phys. D} {\bf 57} (1992), no. 1--2, 96--160.


\bibitem{F17}
C.J. Fan, log-log blow up solutions blow up at exactly m points.
{\it Ann. Inst. H. Poincar\'{e} Anal. Non Lin\'{e}aire} {\bf34} (2017),
no. 6, 1429--1482.


\bibitem{FSZ20}
C.J. Fan, Y. Su, D. Zhang, A note on log-log blow up solutions for stochastic
nonlinear Schr\"odinger equations, \texttt{arXiv:2011.12171,}
to appear in {\it Stoch. Partial Differ. Equ. Anal. Comput.}.




\bibitem{FH14}
P. Friz, M. Hairer,
A Course on Rough Paths.
With an Introduction to Regularity Structures. {\it Universitext}. Springer, Cham, 2014. xiv+251 pp.

\bibitem{G04}
M. Gubinelli, Controlling rough paths. {\it J. Funct. Anal.} {\bf 216} (2004), no. 1, 86--140.


\bibitem{GS11}
S.J. Gustafson, I.M. Sigal, Mathematical Concepts of Quantum Mechanics. Second edition.
{\it Universitext}. Springer, Heidelberg, 2011. xiv+382 pp.

\bibitem{HRZ18}
S. Herr, M. R\"ockner, D. Zhang,
Scattering for stochastic nonlinear Schr\"odinger equations.
{\it Comm. Math. Phys}. {\bf 368} (2019), no. 2, 843--884.


\bibitem{KV13}
R. Killip, M. Visan,
Nonlinear Schr\"odinger equations at critical regularity. Evolution equations, 325-437,
{\it Clay Math. Proc.,} {\bf 17}, Amer. Math. Soc., Providence, RI, 2013.


\bibitem{KK20}
K. Kim, S. Kwon,
On pseudoconformal blow-up solutions to the self-dual Chern-Simons-Schr\"odinger equation:  existence, uniqueness, and instability,
\texttt{arXiv:1909.01055}, to
appear in {\it Mem. Amer. Math. Soc}.

\bibitem{K-M-R}
J. Krieger, Y. Martel, P. Rapha\"{e}l,
Two-soliton solutions to the three-dimensional gravitational Hartree equation.
{\it Comm. Pure Appl. Math.} {\bf 62} (2009), no. 11, 1501--1550.


\bibitem{KS06}
J. Krieger,  W. Schlag,
Stable manifolds for all monic supercritical focusing nonlinear Schr\"odinger equations in one dimension.
{\it J. Amer. Math. Soc.} {\bf 19} (2006), no. 4, 815--920.

\bibitem{KS09}
J. Krieger, W. Schlag, Non-generic blow-up solutions for the critical
focusing NLS in 1-d. {\it J. Eur. Math. Soc.} {\bf 11} (2009), no.1, 1--125.

\bibitem{LeLP15}
S. Le Coz, D. Li, T.P. Tsai,
Fast-moving finite and infinite trains of solitons for nonlinear Schr\"odinger equations.
{\it Proc. Roy. Soc. Edinburgh Sect. A} {\bf 145} (2015), no. 6, 1251--1282.

\bibitem{LeT14}
S. Le Coz, T.P. Tsai,
Infinite soliton and kink-soliton trains for nonlinear Schr\"odinger equations.
{\it Nonlinearity} {\bf 27} (2014), no. 11, 2689--2709.

\bibitem{Ma05}
Y. Martel,
Asymptotic N-soliton-like solutions of the subcritical and critical generalized Korteweg-de Vries equations.
{\it Amer. J. Math.} {\bf 127} (2005), no. 5, 1103--1140.

\bibitem{MM06}
Y. Martel, F. Merle, Multi solitary waves for nonlinear Schr\"odinger equations.
{\it Ann. Inst. H. Poincar\'{e} Anal. Non Lin\'{e}aire} {\bf 23} (2006), no. 6, 849--864.

\bibitem{MP18}
Y. Martel, P. Rapha\"el,
Strongly interacting blow up bubbles for the mass critical nonlinear Schr\"odinger equation.
{\it Ann. Sci. \'{E}c. Norm. Sup\'{e}r.} {\bf 51}  (4) (2018), no. 3, 701-737.

\bibitem{MMT08}
J. Marzuola, J. Metcalfe, D. Tataru, Strichartz estimates and local
smoothing estimates for asymptotically flat Schr\"{o}dinger
equations, {\it J. Funct. Anal.}, {\bf 255} (6) (2008), 1479--1553.

\bibitem{M90}
F. Merle, Construction of solutions with exactly k blow-up points for the Schr\"odinger equation with critical nonlinearity.
{\it Comm. Math. Phys.} {\bf 129} (1990), no. 2, 223--240.

\bibitem{M93}
F. Merle,
Determination of blow-up solutions with minimal mass for nonlinear Schr\"odinger equations with critical power.
{\it Duke Math. J.} {\bf 69} (1993), no. 2, 427--454.

\bibitem{MR03}
F. Merle, P. Rapha\"el, Sharp upper bound on the blow up rate for critical nonlinear
Schr\"{o}dinger equation. {\it Geom. Funct. Anal.} {\bf13} (2003), 591--642.

\bibitem{MR04}
F. Merle, P. Rapha\"el, On universality of blow up profile for $L^2$ critical nonlinear Schr\"odinger
equation. {\it Invent. Math.} {\bf 156} (2004), 565--672.

\bibitem{MR05.2}
F. Merle, P. Rapha\"el,
The blow-up dynamic and upper bound on the blow-up rate for critical nonlinear Schr\"odinger equation.
{\it Ann. of Math.} 161 (2005), no. 1, 157--222.



\bibitem{MR05}
F. Merle, P. Rapha\"el,
Profiles and quantization of the blow up mass for critical nonlinear Schr\"odinger equation.
{\it Comm. Math. Phys.} 253 (2005), no. 3, 675--704.

\bibitem{MR06}
F. Merle, P. Rapha\"el, Sharp lower bound on the blow up rate for critical nonlinear
Schr\"odinger equation. {\it J. Amer. Math. Soc.} {\bf19} (2006), no. 1, 37--90.

\bibitem{MRS13}
F. Merle, P. Rapha\"el, J. Szeftel,
The instability of Bourgain-Wang solutions for the $L^2$ critical NLS.
{\it Amer. J. Math.} {\bf 135} (2013), no. 4, 967--1017.


\bibitem{MRRY20}
A. Millet, A.D. Rodriguez, S. Roudenko, K. Yang,
Behavior of solutions to the 1d focusing stochastic nonlinear Schr\"odinger equation
with spatially correlated noise,
\texttt{arXiv:2006.10695v1.}


\bibitem{MRY20}
A. Millet, S. Roudenko, K. Yang,
Behavior of solutions to the 1d focusing stochastic $L^2$-critical and supercritical
nonlinear Schr\"odinger equation with space-time white noise,
\texttt{arXiv:2005.14266v1.}


\bibitem{P01}
G. Perelman, On the blow up phenomenon for the critical nonlinear Schr\"{o}dinger equation in 1D. {\it Ann. Henri. Poincar\'{e}} {\bf 2} (2001), 605-673.

\bibitem{RS11}
P. Rapha\"el, J. Szeftel,
Existence and uniqueness of minimal blow-up solutions to an inhomogeneous mass critical NLS.
{\it J. Amer. Math. Soc.} {\bf 24} (2011), no. 2, 471--546.


\bibitem{RGBC95}
K.O. Rasmussen, Y.B. Gaididei, O. Bang, P.L. Chrisiansen,
The influence of noise on critical collapse
in the nonlinear Schr\"odinger equation.
{\it Phys. Letters A} {\bf 204} (1995), 121--127.

\bibitem{S09}
W. Schlag,
Stable manifolds for an orbitally unstable nonlinear Schr\"odinger equation.
{\it Ann. of Math.} {\bf (2) 169} (2009), no. 1, 139--227.

\bibitem{SZ19}
Y. Su, D. Zhang,
Minimal mass blow-up solutions to rough nonlinear Schr\"odinger equations.
\texttt{arXiv:2002.09659v1.}

\bibitem{SZ20}
Y. Su, D. Zhang,
On the multi-bubble blow-up solutions to rough nonlinear Schr\"{o}dinger equations,
\texttt{arXiv:2012.14037v1.}


\bibitem{SS99}
C. Sulem, P.L. Sulem,
The Nonlinear Schr\"odinger Equation: Self-Focusing
and Wave Collapse.
{\it Applied Mathematical Sciences} {\bf 139}, Springer, New York,
1999.


\bibitem{T06}
T. Tao, Nonlinear Dispersive Equations. Local and Global Analysis.
CBMS Regional Conference Series in Mathematics, 106.
Published for the Conference Board of the Mathematical Sciences, Washington, DC;
by the American Mathematical Society, Providence, RI, 2006. xvi+373 pp.


\bibitem{W83}
M. Weinstein,
Nonlinear Schr\"odinger equations and sharp interpolation estimates.
{\it Comm. Math. Phys.} {\bf 87} (1982/83), no. 4, 567--576.

\bibitem{W85}
M. Weinstein, Modulational stability of ground states of nonlinear Schr\"{o}dinger equations.
{\it SIAM J. Math. Anal.} {\bf 16} (1985), no. 3, 472--491.

\bibitem{Z17}
D. Zhang, Strichartz and local smoothing estimates for stochastic dispersive equations.
\texttt{arXiv:1709.03812}

\bibitem{Z19}
D. Zhang, Optimal bilinear control of stochastic nonlinear Schr\"odinger equations: mass-(sub)critical case.
{\it Probab. Theory Related Fields} {\bf 178} (2020), no. 1--2, 69--120.


\end{thebibliography}
\end{document}